\documentclass{article}

\usepackage[english]{babel} % 'french', 'german', 'spanish', 'danish', etc. 
\usepackage{xcolor}
\definecolor{weblmcolor}{cmyk}{0.86,0.23,0.44,0.02}
\usepackage{mathdots}
\usepackage{graphicx}
\usepackage[export]{adjustbox}
\usepackage{graphicx}
\usepackage{multicol,multirow}
\usepackage{amsmath,amssymb,amsfonts}
\usepackage{amsthm}
\usepackage{mathrsfs}
\usepackage{mleftright}
\usepackage{url}
\usepackage{rotating}
\usepackage{appendix}
\usepackage[numbers]{natbib}
\usepackage{ifpdf}
\usepackage[T1]{fontenc}
\usepackage{tikz-cd}
\usepackage{float}
\usepackage{newtxtext}
\usepackage{txfonts}
\usepackage{authblk}
\usepackage{tabularx}
\usepackage{verbatim}
\usepackage{tabto}
\usepackage{changepage}
\usepackage{enumitem}
\usepackage{arydshln}
\usepackage{subcaption}
\usepackage{placeins}
\usepackage{blindtext}
\usepackage{booktabs}
\usepackage{MnSymbol}
\newtheorem{theorem}{Theorem}
\newtheorem{lemma}{Lemma}
\theoremstyle{definition}

\newtheorem{example}{Example}
\newtheorem{definition}{Definition}
\newtheorem{corollary}{Corollary}

\numberwithin{equation}{section}
\newlength{\colA}
\newlength{\colB}
\newlength{\colC}
\newlength{\colD}
\newlength{\colE}
\newlength{\colF}
\newlength{\colG}
\DeclareRobustCommand{\tcr}[1]{\textcolor{red}{#1}}
\DeclareRobustCommand{\tcb}[1]{\textcolor{blue}{#1}}
\newcommand{\sm}[1]{\begin{smallmatrix}#1\end{smallmatrix}}
\newcommand{\croot}{\rotatebox[origin=c]{270}{\(\textcolor{red}{\circlearrowright}\)}}
\newcommand{\cabroot}{\rotatebox[origin=c]{270}{\(\textcolor{red}{\circlearrowright}\)}\,\sm{\textcolor{red}{0}\\c=a\cdot1+b}}
\newcommand{\tok}{\rotatebox[origin=c]{270}{\(\textcolor{red}{\nwrsquigarrow}\)}}

\counterwithout{figure}{subsection}
\usepackage[colorinlistoftodos]{todonotes}

\newcommand*{\email}[1]{%
 \footnotesize \href  {#1}\par
    }

\usepackage{geometry}
\usepackage[
  colorlinks=true,
  linkcolor=weblmcolor,
  citecolor=weblmcolor,
  urlcolor=weblmcolor,
  bookmarksopen=true,
  bookmarksdepth=3
]{hyperref}
\AtBeginDocument{%
  \hypersetup{
    pdfborder={0 0 0},
    pdfborderstyle={}
  }%
}
\title{Collatz/Bařina \textit{a}n+b conjecture proofs by 16/64:\\
Cayley a\textit{n}+b tree not Cantor-Gödel-Conway-Wolfram undecidable}

\begin{document}
\author[1]{\small Jan Kleinnijenhuis}
\author[2]{\small Alissa M. Kleinnijenhuis}
\author[3]{\small Mustafa G. Aydogan}
\affil[1]{\footnotesize Vrije Universiteit Amsterdam, The Network Institute}
\affil[ ]{\ {\ {j.kleinnijenhuis@vu.nl  }}}
\affil[2]{\footnotesize Columbia University, Columbia Business School}
\affil[ ]{\ {\ {amk2400@gsb.columbia.edu}}}
\affil[3]{\footnotesize Columbia University, Department of Genetics and Development}
\affil[ ]{\ {\ {mga2157@cumc.columbia.edu}}}

\maketitle

%\begin{frontmatter}

%\author[1]{\small Alissa M. Kleinnijenhuis}
%\author[2]{\small Jan Kleinnijenhuis}
%\affil[1]{\footnotesize Mathematical Institute, University of Oxford}
%\affil[2]{\footnotesize Vrije Universiteit Amsterdam, The Network Institute}
%\affil[ ]{\ {\ { alissa.kleinnijenhuis@maths.ox.ac.uk;j.kleinnijenhuis@vu.nl}}}

% \authormark{Proof decidability \textrm{a}n+\textrm{b}, proof 3n+1 conjecture}
% \keywords{Collatz conjecture, Cayley colour graph}

%v\keywords[MSC Codes]{\codes[Primary]{05C63}; \codes[Secondary]{
% 05C05, 05C20, 05C60, 05C76, 11B50, 11F03 }}

% \keywords{Collatz conjecture, Cayley colour graph}

\abstract{
Lothar Collatz considered biperiodic iterations, and conjectured that picturing   $n\mathrm{\to3}n\mathrm{+1}$ iterations for odd numbers 
and $n\mathrm{\to}n/2$ for even numbers, and proving that all natural numbers have a root path to the `trivial' cyclic root $(c=4)\to2\to1\to c\to\dots$, would connect elementary graph theory and number theory. John Conway replied that a\textit{n}+b conjectures, except 1$n$+1, are undecidable due to a-periodic even $n$=a$n$+b successors. 
The infinite Cayley tree $T_{\ge 0} \subseteq G_{\ge 0}$ of Collatz-Conway $a$$\nmid $n$=a$n$_\text{odd}$+b numbers paired with Cantor binary tree numbers $y_{\ge0}$=1$y$+0 $\in [0,1]_2\stackrel{\textup{bin}}{=}[0,1]_{10}$ in root paths of arrow digits $\tcr{\to0}$ and $\tcb{\to1}$,
, pairs Collatz 3$n$+1 congruence classes $[4,16]_{m=18}$ ($d_{\textup{ensity}}$=$\frac{2}{18}$) and $[0,1]_{2}$($d$=$\frac{2}{2})$ in a biperiodic, 2-colour, breadth-first tape of number pairs $\sm{y\\n}$ of red and blue arrow function digits $\sm{\textcolor{red}{0}\\L=fg^p=(n-b)/a \cdot2^p}$ and $\sm{\textcolor{blue}{1}\\U=g^q=n\cdot2^q}$ with cyclic root $\cabroot$. The $E=V=F=\mathbb{N}_0$ tree $T_{\ge 0}$ has arrows ($E_{\textup{dges}}$) to wondrously ordered arrowhead nodes $a$$\nmid$n$=a$n$_\text{odd}$+b ($V_{\textup{ertices}})$ with same-numbered upwardly open sites above  ($\textup{sur}F_{\textup{aces}})$, paired with well-ordered root paths $y$, projecting to the Cartesian $y$-axis well-ordered sites $n$, paired with wondrously ordered root paths $y$. 

Stephen Wolfram pictured iterations for myriad `rules'. The solution to his (one?) $S$-combinator challenge is that one breadth-first Cayley $T_{\ge0}$ traversal alternately passes red and blue number pair nodes  $\tcr{\to}s_{\tcr{i}}\tcb{\to}S_{\tcb{j}}$ in same-generation cotrees $\tcr{\to}t_{\tcr{i}}\tcb{\to}T_{\tcb{j}}$, generated by biperodic color-shifting iterators $UL^{\tcr{i}=1,2,3,\dots}$ and $LU^{\tcb{j}=1,2,3,\dots}$. This disproves a\textit{n}+b-conjecture undecidability by Cantor diagonals, Gödel numbering, Conway unpredictability and Wolfram irreducibility, allowing proofs by the Cantor 1$n$+0 binary tree, Gödel colons, the Conway 1$n$+1 tree and the Wolfram-Turing-Leibniz tape ruliad.

Tree $T_{\ge0}$ transforms into the integer Cartesian upper half graph $C_H$ with on its $(x,y)$-coordinates in Quadrant I well-ordered root paths of blue cotree $T_{\tcb{j}}$ numbers and in Quadrant II red cotree $t_{\tcr{i}}$ numbers, providing a Hamiltonian tour from and back to the root $\cabroot$ at its origin $(0,0)$. Its automorphism graph Aut$(T_{\ge0},[L,U])$ shows 1 red $L$eftward and $a-1$ blue $U$pward subtree-cotree ladder graphs of geometric series of subtree-cotree congruence classes---infinite arithmetic progressions---$\ell m+\alpha$, for periods $\ell=0,1,2,\dots$, and residues $\alpha\in\{0,\dots,a-1\}$. The `$c$otree $c$ongruence $c$lasses' tape $T_{CCC}$ visualizes the periods $\ell=0,1,2,\dots$ `cotree density' Hamiltonian, and the period $\ell=0$ `root-trunk', proof or refutation. The Turing-Leibniz tape $T_L$ with infinite rows of $\eta_{\ge1}\in\mathbb{N}_1$ proven $a$n$_{\textup{odd}}$+b conjectures, for $\eta=0,1,2,\dots$ the Conway 1$n$+1 tree, the Collatz 3$n$+1 tree, the Bařina $7n\pm1$ tree, $\dots$. Each column shows Turing-computable roots $c=a+1$, $c_{\eta0}=2^{\eta}+1=2$ and $c_{\eta\ge1}=2^{{2^{\eta-1}+1}}$, and for $\eta\ge1$: lowest trunk nodes $c^2$, periodicities  $\frac{c}{2}(c-1)$ of yellow odd nodes $a\nmid$n$_{\textup{odd}}$, and Lean-certified `cotree density' Hamiltonian tour, and `root-trunk', proofs $c^2<\frac{c}{2}(c-1)^2$: for Collatz$_{\eta1}$ 3$n$+1 ($16<18$) and for Bařina$_{\eta2}$ $7n\pm1$ ($64<196$). Successive $a=c-1$ parameters are the short legs of primitive Pythagorean triangles given by the sum of legs in the previous one (OEIS 365577).

Supplementary proofs complete hitherto heuristic proofs: German/Klop staircase term rewriting of arbitrary numbers to the trivial root, the Turing halting proof in the trivial root rather than in diverging trajectories or non-trivial cycles, and Terence Tao's recent 'almost all' proof on the Syracuse tree of odd numbers.
}

\newpage
\renewcommand{\contentsname}{}
\vspace{-0.6cm}
\small{\tableofcontents}
\newpage
\small{\listoffigures}
\small{\listoftables}

%\end{frontmatter}

\phantomsection
\normalsize{}

\newpage
\section{Introduction: wondrously ordered Collatz numbers}

\noindent The Collatz function $C$ returns $n\to n/2$ arrows for even numbers and $n\to 3n+1$ arrows for odd numbers \cite{RN9}. The Collatz conjecture to be proven maintains that iterations of the Collatz function give each natural number a \textit{root path} to the branching number $c=3\cdot1+1=4$ in the cyclic \textit{trivial root} $\dots \to (4\to 2\to 1\to 4\to \dots )$ \cite{RN9,taoconway,taoFridman25}. 
Example \ref{ex:firstEx} shows the converging root paths of $7$, $8$, and $9$. The root path of $9$ converges to that of $7$, which  converges to that of $8$, which converges to $c=4$ in the cyclic root.

\begin{figure}[h]
\caption{The Collatz graph $G_C$ (\textit{k}=13 breadth-first levels)}
\label{fig:gcol}
\vspace{-0.3cm}
\center{\includegraphics[width=0.91\textwidth]{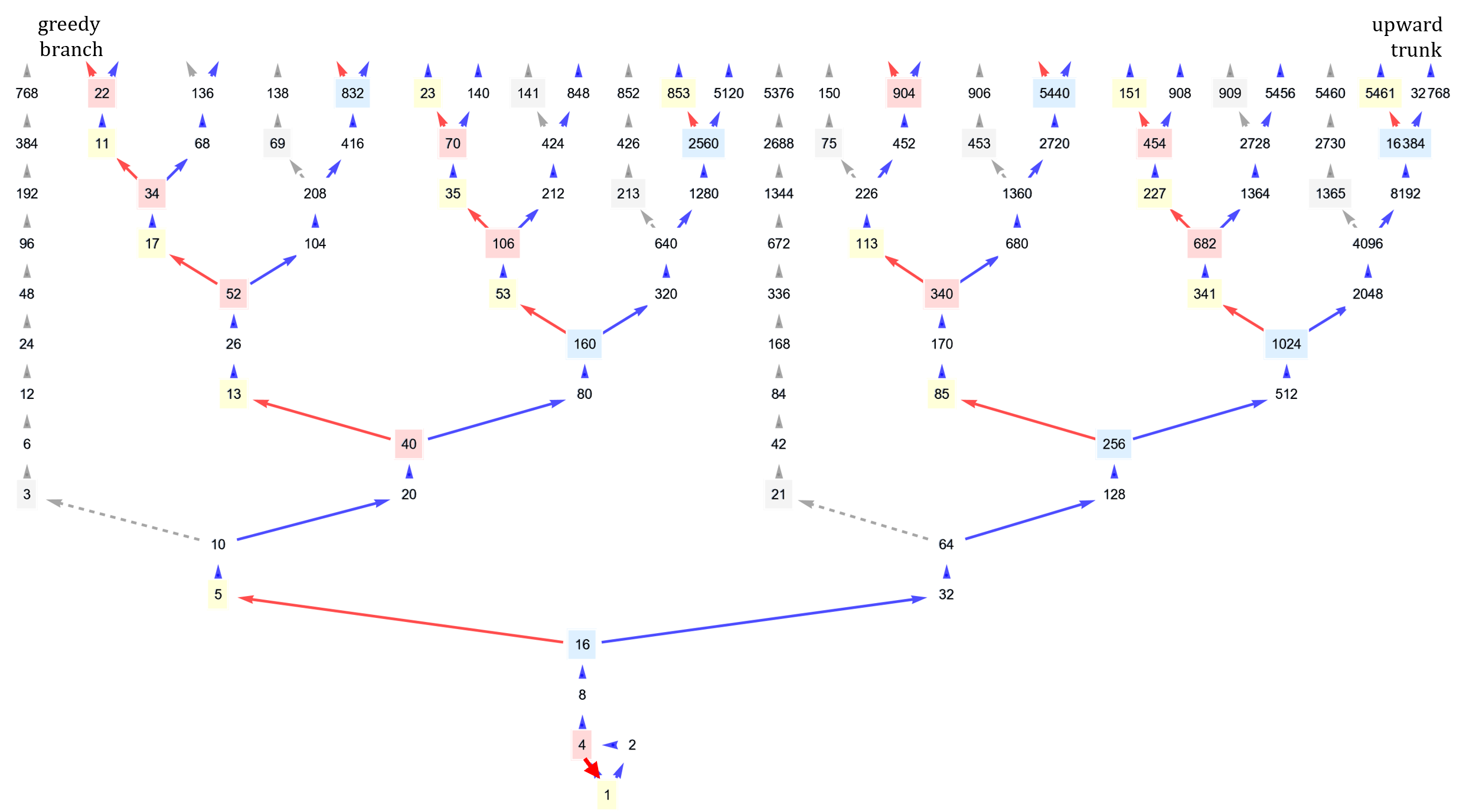}}
\end{figure}

\begin{example}
\label{ex:firstEx}
\textbf{Three root paths converging to the trivial root number} $c=3\cdot1+1=4$ \\
\par $7\textcolor{red}{\;\to\;}22\textcolor{blue}{\;\to\;} 11\textcolor{red}{\;\to\;} 34\textcolor{blue}{\;\to\;} 17\textcolor{red}{\;\to\;} 52\textcolor{blue}{\;\to\;} 26\textcolor{blue}{\;\to\;} 13\textcolor{red}{\;\to\;} 40\textcolor{blue}{\;\to\;} 20\textcolor{blue}{\;\to\;} 10\textcolor{blue}{\;\to\;} 5\textcolor{red}{\;\to\;} 16\textcolor{blue}{\;\to\;} (8\textcolor{blue}{\;\to\;} \dots )$ \par
$8\textcolor{blue}{\;\to\;} (4\textcolor{blue}{\;\to\;} 2\textcolor{blue}{\;\to\;} 1\textcolor{red}{\;\to\;} 4\textcolor{blue}{\;\to\;} \dots )$. \par
$9\textcolor{red}{\;\to\;} 28\textcolor{blue}{\;\to\;} 14\textcolor{blue}{\;\to\;} (7\textcolor{red}{\;\to\;} \dots)$
\end{example}
\noindent The brackets ($\cdots$) indicate root path convergence. The red arrows highlight $n \textcolor{red}{\to} 3n+1$ expansions for odd numbers. The term 'wondrous numbers' \cite{RN31,RN1} indicates that successive natural numbers have \textit{wondrously}  ordered, \textit{not obviously simple} ordered \cite{RN23}, root paths, just as successive numbers have very different prime factorisations, e.g. $\sm{7\\7^1},\; \sm{8\\2^3},\; \sm{9\\3^2}$. Fig.\ref{fig:gcol} of the Collatz graph $G_C$ shows $k=13$ breadth-first levels of numbers generated from the trivial root a$n$+b number $c=3\cdot1+1=4$.

\vspace{-0.3cm}
\noindent \textbf{Legend.}  The inverse functions for even numbers $g:n\to2n$ and odd numbers $f:n\to(n-1)/3$ generate the Collatz graph $G_C$. A blue \textit{upward path} $U: n\textcolor{blue}{\to U(n)}=n\cdot 2^q$ of a red or blue branching number to its blue upward child $\textcolor{blue}{U(n)}$ consists of $q$ successive $g$ arrows. E.g. $U(16)=g^4(16)=256$ is the blue upward arrow shortening $U\colon 16\textcolor{blue}{\to} 32\textcolor{blue}{\to} 64\textcolor{blue}{\to} 128\textcolor{blue}{\to 256}$ (blue). The red \textit{leftward path} $L: n\textcolor{red}{\to L(n)}=fg^p=\frac{n-1}{3} \cdot 2^p$ to its red leftward child $L(n)$ consists of an $f$ arrow to a yellow number in odd congruence classes $[1,5]_6$, followed by $p$ successive $g$-arrows. E.g. $L(16)=fg^3(16)=40$ is the red arrow shortening $L\colon 16\textcolor{red}{\to} 5\textcolor{blue}{\to} 10\textcolor{blue}{\to} 20\textcolor{blue}{\to}\textcolor{red}{40}$. Grey numbers in $[0,3]_6$ divisible by 3 do not branch by $f:n\to(n-1)/3$ (vertical line graphs). \quad\textbf{\textbar}

The $k=13$ successor generations in Fig.\ref{fig:gcol} of the lowest branching number $c=3\cdot1+1=4$ come from $k=13$ \textit{breadth-first} iterations of the functions $f:n\to(n-1)/3$ and $g:n\to2n$. These iterations generate high numbers such as $5440, 5460, 5461$, and on the rightmost path, labeled the \textit{upward trunk}, $16384$ and $32768$. In the leftmost path, labeled the \textit{greedy branch} \citep[Seq.A225570]{RN8}, they generate the low number $22$. Example \ref{ex:firstEx} shows the even lower rooted numbers $f(22)=7$ and $g^2f(7)=9$ at iterations $k=13+1$ and $k=14+3=17$.
Note the inconsistent ‘diagrammatic’ \citep[p.33]{RN3240}, ‘leftmost-innermost’ \citep[p.906]{RN23} arrow order notation to denote composite functions such as $g^2f(7)$, avoiding the suggestion of multiplication by a consistent diagrammatic, bracketless, notation $7g^2f$.

To be proven is that all root paths are included in one connected tree, labeled $T_{\ge0}$, generated from the lowest branching number $c=3\cdot1+1=4$ \citep[ch.3]{RN3239}. It takes $k=17$ iterations to generate all the numbers $n<10$ from it (Example \ref{ex:firstEx}). By 2025, the record for the number below which all numbers can be generated from it is $n<2^{71}$, held by an algorithm developed by David Bařina \cite{barina2025}.

\begin{figure}[h!]
\caption{The 3-regular 3-colour Cayley graph $T_{\ge 0}$, \textbf{ a.} as a planar graph, \textbf{ b.} as a fractal graph}
\label{fig:tge0}
\vspace{-0.5cm}
\noindent \center{\includegraphics[width=0.95\textwidth]{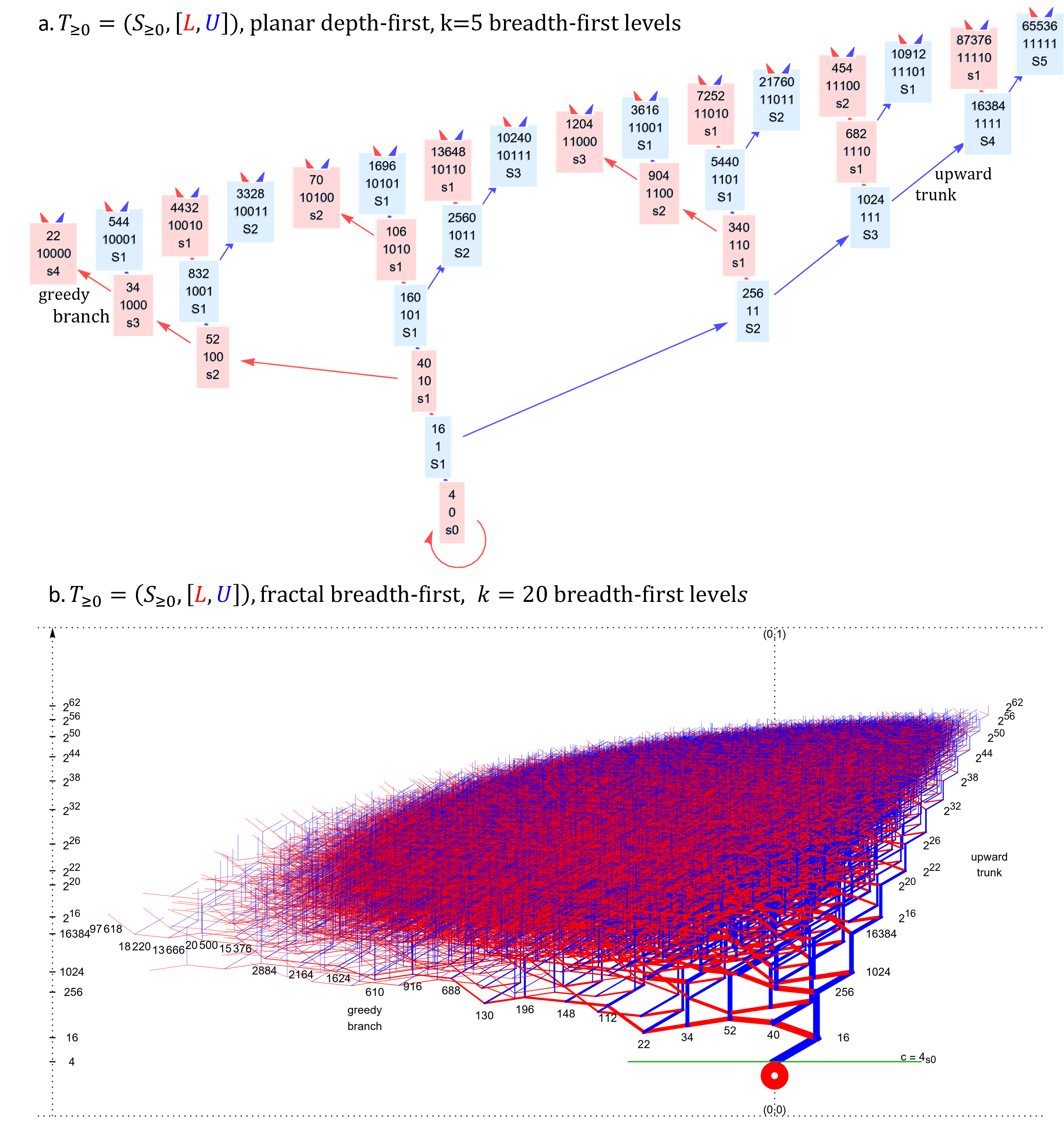}}
\end{figure}

The functions $f:n\to(n-1)/3$ and $g:n\to2n$ and their inverses $f^{-1}:n\to3n+1$ and $g^{-1}:n\to n/2$ generate \textit{rest classes} $[r]_m$ (Eqs.\ref{eq:fmf},\ref{eq:gmg}, notation \cite{RN3231}), being infinite \textit{arithmethic progressions} of numbers $\ell m+r$ in positions $r$ in periods $\ell=0,1,2,\cdots$ with periodicity $m$. Purple a$n$+b classes $S_{\ge0}\;=\;[3\cdot\colorbox{yellow}{1}+1,\;3\cdot\colorbox{yellow}{5}+1]_{3\cdot6}\;=\;[\colorbox{purple}{\textcolor{white}{4}},\colorbox{purple}{\textcolor{white}{16}}]_{18}$ comprise both red children $s_{\ge0}$ generated by $L=fg^p$ (Eq.\ref{eq:Lfgp}) and blue children $S_{\ge1}$ generated by $U=g^q$ (Eq.\ref{eq:Ugq}).
For them $f: (n-1)/3\in
[\colorbox{yellow}{\textcolor{black}{1}},
\colorbox{yellow}{\textcolor{black}{5}}]_{6}$  gives yellow odd classes indivisible by $a=3$ with `perfect' periodicity $1\cdot2\cdot3=1+2+3=6$. To be proven is that Collatz 3$n$+1 tree $T_{\ge0}$ (Figs.\ref{fig:gcol},\ref{fig:tge0}) includes all numbers in $S_{\ge0}$ (purple, Eq.\ref{eq:sge0col}) $=s_{\ge0}$ (red $L$eftward, Eq.\ref{eq:Lfgp}) $\,\cup\; S_{\ge1}$ (blue $U$pward, Eq.\ref{eq:Ugq}).
\begin{align}
f^{-1}&:n\quad\quad\in\;[\colorbox{yellow}{1},
\colorbox{lightgray}{3},
\colorbox{yellow}{5}]_6\to
&&f:n\quad\quad\quad\in[\colorbox{purple}{\textcolor{white}{4}},
\colorbox{lightgray}{10},
\colorbox{purple}{\textcolor{white}{16}}]_{m=18}\to \nonumber\\
&\;\;\; 3n+1 \in
[\colorbox{purple}{\textcolor{white}{4}},
\colorbox{lightgray}{10},
\colorbox{purple}{\textcolor{white}{16}}]_{m=18} 
&&\quad(n-1)/3 \in\;
[\colorbox{yellow}{1},
\colorbox{lightgray}{3},
\colorbox{yellow}{5}]_6 & \label{eq:fmf} \\
g^{-1}&:n\in[0,2]_4\to n/2\in[0,1]_2 
&&g:n\in[0,1]_2 \to 2n\in [0,2]_4 & \label{eq:gmg} \\
(a-1)&=2 \textrm{ branching classes }S_{\ge0}=[\colorbox{purple}{\textcolor{white}{4}},\colorbox{purple}{\textcolor{white}{16}}]_{m=18}&&\textrm{for which }f: (n-1)/3\in
[\colorbox{yellow}{\textcolor{black}{1}},
\colorbox{yellow}{\textcolor{black}{5}}]_{6} & \label{eq:sge0col}\\
&\textrm{red \textit{L}eftward arrow to red $L$-child }
&& L=fg^{\textcolor{red}{p}}:n\;\textcolor{red}{\to}\; \textcolor{red}{L(n)}\;\,=\;(n-b)/a\cdot2^{\textcolor{red}{p}}\;\;\in \;\;s_{\ge0} \label{eq:Lfgp} \\
&\textrm{blue $U$pward arrow to blue $U$-child }
&& U=\;g^{\textcolor{blue}{q}}\;\;:n\;\textcolor{blue}{\to}\; \textcolor{blue}{U(n)}\;=\;n\cdot2^{\textcolor{blue}{q}} \;\;\in \;\;S_{\ge1}\label{eq:Ugq}
\end{align}
John Conway (1937-2020) states ``in one of the earliest mathematical papers on the $3n+1$ problem" \citep[p.221]{RN9}, that the $3n+1$ function defined by Lothar Collatz (1910-1990) “prompts consideration of more general functions” \cite{RN3209}. Conway introduces a$n+$b functions named after the $f^{-1}:n\to \textrm{a}n+\textrm{b}$ subfunction for odd numbers, still assuming halving even numbers by $g^{-1}:n\to n/2$, with inverse number generating functions $f:n\to(n-b)/a$ and $g:n\to2n$ (Eqs.\ref{eq:fmf},\ref{eq:gmg}). Whether an a$n$+b function generates all natural numbers from its trivial root $c=a\cdot1+b$ is \textit{undecidable} according to Conway. For which branching numbers which powers-of-two $p$ in $L=fg^p:n\to(n-b)/a\cdot2^p$ (Eq.\ref{eq:Lfgp})  give `defined' yellow odd numbers (with purple successors) rather than `undefined' fractions would be \textit{a-periodic} and \textit{unpredictable} \cite{RN3209,fractran,endrullisGrabmayer}. 

Arthur Cayley (1821–1895), regarded as the third pioneer of graph theory \citep[3A,6A]{RN3245} after Leonhard Euler (1707–1783) \citep[1A]{RN3245} and William Hamilton (1805–1865) \citep[2C]{RN3245}, invented colour graphs. As the publishers of his days did not print colours, Cayley explains that ``the continuous lines represent black lines, and the broken lines, red lines'' \cite{RN3245}. Same-colour arrowhead nodes require however also broken nodes. The \textit{infinite} 3-colour Cayley tree $T_{\ge0}$ of a Collatz/Conway a$n$+b function \textit{periodically} excludes gray odd arrowhead numbers $[\colorbox{lightgray}{3}]_6$ divisible by $a$ and their even $3n+1$ numbers $[\colorbox{lightgray}{10}]_{m=18}$ (Fig.\ref{fig:gcol}) (Figs.\ref{fig:tge0},\ref{fig:treeP1},\ref{fig:tree7pm1},\ref{fig:tree3m1},\ref{fig:tree5p1}).

The $3$-\textit{regular} \citep[10A]{RN3245}, \textit{biperiodic} \cite{collatz74} Cayley $an + b$ graph is not Conway-undecidable. By Eq.\ref{eq:rootPairs}, each node $n$ has 1 upward child $U(n)=g^{\textcolor{blue}{q}}(n)$ (Eq.\ref{eq:Ugq}), 1 leftward child $L(n)=fg^{\textcolor{red}{p}}(n)$ (Eq.\ref{eq:Lfgp}) and 1 parent, including the cyclic root $\croot c=a\cdot1+b$ (Eq.\ref{eq:Croot}). The red and blue powers-of-two, $\textcolor{red}{p}$ and $\textcolor{blue}{q}$, show up in congruence classes graphs (Figs.\ref{fig:cc3p1},\ref{fig:cc3m1},\ref{fig:cc7pm1}) as path lengths of $g$-arrows to purple $S_{\ge0}$ classes, which are arithmetic progressions $\ell m+r$ for periods  $\ell=0,1,2,\dots$: for $3n+1$ rest classes $r\in[4,16]$ with periodicity $m=18$. To be shown is that the period $\ell=0$ \textit{root-trunk-test} proves an a$n$+b conjecture by  $(U(c)=c^2\in(\ell=0))<m$ (Eq.\ref{eq:rtt}): Collatz $3n+1$ by $U(c)=16$ and Bařina $7n\pm1$ by $U(c)=64$. 

\newpage
\vspace{-0.5cm}
\noindent \textbf{Legend } The 3-regular Cayley graph $T_{\mathrm{\ge }0}=(S_{\ge 0},[L,U])$ has arrows shortening paths (Fig.\ref{fig:gcol}), to the leftward and upward child, e.g. $L=fg^3\colon\textcolor{blue}{16}\textcolor{red}{\to}\left(5\textcolor{red}{\to} 10\textcolor{red}{\to} 20\textcolor{red}{\to}\right)\textcolor{red}{40}\,$ and $U=g^4\colon \textcolor{blue}{16}\textcolor{blue}{\to}\left(32\textcolor{blue}{\to}64\textcolor{blue}{\to}128\textcolor{blue}{\to}\right) \textcolor{blue}{256}$, as $(40-1)/3=13$ and $(256-1)/3=85$ are odd. \textbf{a}. Its planar depth-first graph shows $n$=3n$_{\textup{odd}}$+1 numbers of upwardly open sites, bordered to the left by blue arrows to upward successors
of the leftward child in generations $\textcolor{blue}{S_1,S_2,S_3,\dots}$, and to the right by red arrows to the leftward successors of the upward child in generations $\textcolor{red}{s_1,s_2,s_3,\dots}$. \textbf{b}. Its fractal breadth-first graph shows $k=20$ iterations generating $2^{20}$ wondrously ordered $n$=3n$_{\textup{odd}}$+1 numbers $\croot\tcr{4}\textcolor{blue}{\to16}\tok\tcr{\;40}\tcb{\to256}\tok\tcr{\;52}\tcb{\to160}\tcr{\to340}\tcb{\to1024}\tok\tcr{\;34}\tcb{\to}\cdots$, projected well-ordered to the Cartesian $y$-axis $4\to16\to22\to34\to52\to\cdots$ contracted to the open interval $(0,1)$.\textbf{ {\textbar}}

\section{Proof approach and proof preview}
\label{sec:proofa}

The combined insights of Collatz \cite{RN9,collatz74,collatz1985geom}, Conway \cite{RN3209,fractran}, Cantor \citep[S6]{cantor1883bin} and Wolfram \cite{RN32,RN23,combinators} provide a proof approach and proof preview in four steps deciding whether a Collatz-Conway a$n$+b conjecture is proven or refuted on its Cayley graph $T_{\ge0}$ (Figs.\ref{fig:tge0},\ref{fig:treeP1},\ref{fig:tree7pm1},\ref{fig:tree3m1},\ref{fig:tree5p1}).

The opening sentences of \textit{The Motivation and Origin of the }$3n+1\ $\textit{Problem} by Lothar Collatz (1910-1990) read: ``There are numerous connections between elementary number theory and elementary graph theory. One such connection can be made using the fact that one can picture a number theoretic function $f\left(n\right)$ with a directed graph" where each iteration is drawn with ``an arrow from $n$ to $f\left(n\right)$" \cite{RN3051}. ``Through the lectures of Edmund Landau, Oscar Perron and Issai Schur during my time as a student between 1928 and 1933, I became interested". The lectures by Issai Schur \cite{schur1928} show the prominence in German mathematics before 1933 of \textit{unpictured} Cayley graphs \cite{noether1882,peiferCayleyDehn}. Collatz was both aware of \textit{biperiodic} (“doppelt periodische”) functions \cite{collatz74,collatz1985geom} and of massive iterative calculations \cite{althofer2019}, and font of \textit{picturing} trees by 2-colour shadows showing 1-colour branches \cite{bredendiek}. In the 1980s, emerging computer graphics enabled left-right-symmetric 2-colour tessellations \cite{collatz1985geom}, as independently exploited by Collatz \cite{collatz1985geom}, Mandelbrot $\&$ Frame \cite{RN3247}, Penrose \cite{penrose1979penta} and Wolfram \cite{RN23,combinators}. \textit{Computational mathematics} programs to both compute and picture directed multi-colour number graphs \cite{rahnSultanow2021,weisNumber}, notably Wolfram \textit{Mathematica}, appeared only in the late 1980s \cite{RN32}.

Conway generalized Collatz 3$n$+1 iterations to \textit{a-periodic}, \textit{unpredictable}, iterations and undecidable a$n$+b conjectures \cite{fractran,RN3209}.  Wolfram pictured, starting from “initial conditions” \citep[index]{RN23}\cite{wolfram21rulial,wolframFridman25}, iterations of thousands of `simple', one-periodic, 2-colour, “rules” (functions), behaving as breadth-first level $\tok$-functions (Eq.\ref{eq:sorS}) towards incrementally broader node levels of a “trivariate”, 2-colour (black-white), tree / triangle in a `cellular automation', in an integer Cartesian coordinate system \cite{RN23,combinators}. Wolfram’s hypothesis of “computational irreducibility” \citep[chs.11-12]{RN23} for specific rules \citep[rules 30,110]{RN23} states the “intrinsic randomness” \citep[p315-326]{RN23}, or `a-periodicity' \cite{RN3209}, of black cells/coordinates.

For each Collatz-Conway a$n$+b function (Figs. \ref{fig:tge0} Collatz 3$n$+1, \ref{fig:treeP1} Conway 1$n$+1, \ref{fig:tree7pm1} Bařina 7$n\pm$1, \ref{fig:tree3m1} 3$n$-1, \ref{fig:tree5p1} 5$n$+1) the biperiodic, 3-colour, 3-regular Cayley tree
$T_{\ge0}\subseteq G_{\ge0}$ \citep[10A]{RN3245}\citep[p904,p1029]{RN23}\cite{petersen} of even  $a$$\nmid$$n$=an$_\textup{odd}$+b arrowhead nodes ($V_{\textup{ertices}}$) and above them their upwardly open sites ((sur)$F_{\textup{aces}}$) "should be in the first instance drawn with the arrows" ($E_{\textup{dges}}$) \citep[p181]{RN3242} to same-colour arrowhead nodes. From the cyclic root $\cabroot$ upwards, its number node pairs $\sm{y\\n}$ (Eqs.\ref{eq:yn}-\ref{eq:rootPairs}) pair even  $a$$\nmid$$n$=an$_\textup{odd}$+b root successor arrowhead number nodes and above them their same-numbered upwardly open, although narrowing, sites, to the Cantor binary number tree \citep[S6]{cantor1883bin} of root paths $y$, alternately with a terminal even red $L$eftward arrow $\tcr{0}$ and a terminal odd blue $U$pward arrow $\tcb{1}$ to same-colour arrowhead nodes. Observation shows that the Cayley tree $T_{\ge0}\subseteq G_{\ge0}$ is the $E=V=F=\mathbb{N}_0$ tree, quite unlike the, in each direction closed, $V+F-E=2$ Euler polyhedron, regarded as the birth of topology \cite{richeson2012euler}. Mathematical proof textbooks may not cover observation-based geometric, topological and graph-theoretic proofs  \cite{chartrand2013mathematical}. 

Biperiodic binary functions \cite{collatz74,collatz1985geom} defined by Georg Cantor (1845-1918) \citep[\S6]{cantor1883bin} lengthen the root path number $y$ to each, either red or blue ('purple'), node at breadth-first level $\tok k-1$ with a terminal $\tcr{0}$ or a terminal $\tcb{1}$ (Eq.\ref{eq:rootPairs}) to the root paths of its red $L$eftward child (Eq.\ref{eq:redChild}) and of its blue $U$pward child (Eq.\ref{eq:blueChild}) at the next higher $k$-digits breadth-first level $\tok k$ (Eqs.\ref{eq:sorS}-\ref{eq:bf7pm1}, Figs.\ref{fig:tge0},\ref{fig:treeP1},\ref{fig:tree7pm1},\ref{fig:tree3m1},\ref{fig:tree5p1}). Breadth-first level root successor level $\tok k$ (for $k=1,2,3,\dots$), comprises an exponentially higher number of $2^{k-1}$ root paths of $k-1$ binary digits, alternately with an even red $L$eftward terminal arrowhead $\tcr{\to0}$ and an odd blue $U$pward terminal arrowhead $\tcb{\to1}$ (Lemma \ref{lem:2comb}). A biperiodic, 3-colour, left-right-asymmetric, Cayley tree $T_{\ge0}$ shows up with at breadth-first level $k$ a leftmost root path $1\tcr{0}^{k-1}$ and rightmost root path $\tcb{1}^k$ (Figs.\ref{fig:tge0},\ref{fig:treeP1},\ref{fig:tree7pm1},\ref{fig:tree3m1},\ref{fig:tree5p1}, Eq.\ref{eq:tleg}).

The approach to prove plausible a$n$+b conjectures ($T_{\ge0}=G_{\ge0}$), and to refute others ($T_{\ge0}\subset G_{\ge0}$) is to preview in four subsections \ref{sec:31eqs}-\ref{sec:35eqs} observable equations in (transformed graphs of) Cayley tree $T_{\ge0}\subseteq G_{\ge0}$ without a terminal $\square$, and to trace in their sub-subsections lemmas and theorems on unobserved `transformations' or `metamorphoses' with a $\square$. The term `graph transformation' \cite{endrullis2024term,ehrig2015,biggs1993algebraic} (Roman origin) should not suggest inventions of otherwise non-existent transformations. The term `graph metamorphoses' (Greek origin) should not\ suggest mathematically non-verifiable observations because of debated metamorphoses reported in politics (Caesar 100-44bC, `divide and conquer'), war (Homer, 8th century bC, Trojan horse), economics (Adam Smith, 1723-1790, `invisible hand'; Karl Marx, 1818-1883, `hinter dem Rücken'), religion (reincarnation, resurrection) and literature (Franz Kafka, 1883-1924, ``Metamorphosis", ``Die Verwandlung"). 
Their introduction draws on Douglas Hofstadter’s improvisations on math dialogues by Plato ($\pm$427-347bC) \citep[ch7]{platoRepublic}\cite{platotimaeus,platoparmenides}, whose 1979 bestseller \textit{Gödel, Escher, Bach} \cite{RN31} revives them with more recent insights, such as Turing-computable numbers (Table \ref{fig:TLEE}), undecidability by Gödel numbers (Lemmas \ref{lem:cantorno}-\ref{lem:wolframno}), multiple periodic symmetries in the music of Bach, and geometric series in paintings of Escher and in the improvised dialogues on geometric 3$n$+1 iterations between Achilles and Zeno's Turtoise, invariably positioned at Achilles' infinite horizon in the perception of the Turtoise \citep[p240-245,400-402]{RN31}. 

Each of the subsections \ref{sec:31eqs}–\ref{sec:35eqs} begins with a list of observable equations without a terminal $\square$ on a transformed graph defining Cayley tree $T_{\ge0}$, places them in their historical context without explaining them separately, and presents lemmas and theorems, each with a terminal $\square$, on the unobserved transformation yielding the transformed graph, previewed in subsection \ref{sec:transmeta}. Equations without ‘=’ express one (notation of an) observation. Equations with ‘=’ express two of them, a left-hand and a right-hand side. The left-hand-side function formulas denote right-hand-side paths of serial arrows to observations, or sequences of parallel, random-access, arrows to them. Table \ref{fig:TLEE} (subsection \ref{sec:35eqs}) comprises the formula of successive binary numbers on the infinite Turing-Leibniz tape generating successive cyclic roots $c$=a+1 of an infinite series of Collatz-Conway $n$=an$_\textup{odd}$+b conjectures proven by Hamiltonian tours (Eq.\ref{eq:cdp}) and computationally reduced root-trunk tests (Eq.\ref{eq:rtt}), starting with the proofs of the first two Mersenne prime $a$n$_\textup{odd}$+b conjectures included also on the Euclid-Euler tape (1: Collatz 3$n$+1, 2: Bařina 7$n\pm$1).  

\subsection{Unobservable Cayley tree $T_{\ge0}$ transformations, observed transformed graphs} 
\label{sec:transmeta}
The prisoners in Plato's cave observe object projections/shadows to a wall \citep[ch7]{platoRepublic}.  Euclid's (323-283 bC) 2-colour planar geometry explores 1-colour projections/pictures on a 1-colour wall/paper/board/screen. The biperiodic 3-colour Cayley tree $T_{\ge0}$ comprises an infinite Turing-Leibniz tape of breadth-first ordered number pairs $\sm{y\\n}$ (Eqs.\ref{eq:yn},\ref{eq:rootPairs}) of wondrously ordered $n$=a$n_{\textup{odd}}$+b numbers paired with well-ordered, breadth-first ordered, binary numbers of root paths $y$  in Cantor's binary tree $T_{\ge0}$ (\citep[\S6]{cantor1883bin})  (Figs.\ref{fig:tge0},\ref{fig:treeP1},\ref{fig:tree7pm1},\ref{fig:tree3m1},\ref{fig:tree5p1}, Eqs.\ref{eq:sorS}-\ref{eq:bfpar}). They do not observe the Cayley tree $T_{\ge0}$ as the 3-colour infinite propositional logic graph \cite{knuthAlltrees} with $\sm{\tcr{0}\\\top,\textup{True}}$ and $\sm{\tcb{1}\\\bot,\textup{False}}$ (Lemma \ref{lem:knuth}). They observe finite shadows of the infinite Cayley tree (Fig.\ref{fig:tge0}b) projected to the open $y$ interval (0,1) in a planar Cartesian ($x,y$) coordinate system (Fig.\ref{fig:tge0}b) projected by an infinite number of geometrically contracting heights of breadth-first-levels (Fig.\ref{fig:tge0}b,Theorem \ref{the:fractal}). 

The infinite-to-finite transformation is known from Plato's dialogues on the paradox of Zeno ($\pm$ 490-430 bC). According to a slow turtoise, ahead of super-athlete Achilles, Achilles will never surpass the turtoise at his infinite horizon \cite{platoparmenides}\citep[400-402]{RN31}. Zeno's paradox involves different measures of space and time when moving in the same direction, dealt with in the special relativity theory of Albert Einstein (1879-1955) using Lorentz contraction (after Hendrik Lorentz 1853-1926) visualized in the paintings of M.C. Escher (1898-1972), which were exhibited to ICM 1954 participants, including future Noble prize laureate 2020 Roger Penrose \cite{penrose1979penta}. Escher's ``circle limit" paintings picture the infinite number of steps that Achilles has to take towards the turtoise at his horizon as a finite, surpassable, geometric series sum of an infinite number of geometrically contracting heights of breadth-first-levels. 

The infinite Turing-Leibniz tape of number pairs $\sm{y\\n}$ (Eqs.\ref{eq:yn},\ref{eq:rootPairs}) with well-ordered, breadth-first ordered root path numbers $y$ in the biperiodic, 3-colour, Cayley tree $T_{\ge0}$ projects to the $y$-axis of the Cartesian coordinate system number pairs $\sm{y\\n}$ of well-ordered $n$=a$n_{\textup{odd}}$+b numbers paired with wondrously ordered root paths $y$. The geometric series contraction to the open y-axis interval (0,1) shows a vertical line graph with wondrous node distances, decreasing toward the upper bound of (0,1) (metamorphosis 1, Section \ref{sec:31eqs}, Fig.\ref{fig:tge0}b, Eq.\ref{eq:bfpar}, Theorem \ref{the:fractal}), not revealing in a 1-colour projection its biperiodic, bicolour nature.

Instead of color branches to same-colour arrowhead nodes growing from the cyclic Cayley tree $T_{\ge0}$ root $\cabroot$ they observe series of non-congruent but still `isomorphic', contracting and lengthening, one-colour, shadows, not able to distinguish red infinite series of cotrees $t_{\tcr{i}=1,2,3,\dots}$ with red same-generation arrowhead successors $s_{\tcr{i}}$ (Figs.\ref{fig:utree}bdf) from blue cotrees $T_{\tcb{j}=1,2,3,\dots}$ with blue same-generation arrowhead successors $S_{\tcb{j}}$ (Section \ref{sec:32eqs}, Figs.\ref{fig:ltree}bdf, \ref{fig:utree}bdf, Eqs.\ref{eq:dfsi},\ref{eq:dfSj}). 
Plato considers the 3-regular \citep[11A]{RN3245},\cite{petersen} cube representing Earth with an inscribed 4-regular octahedron representing Air, having 3 respectively 4 edges at each vertex \cite{platotimaeus}.  In 1895, Röntgen invented `new light beams' visualizing the inscribed bones graph of a body circumference graph \cite{rontgen1895}. The prisoners do not see the inscribed octahedron of a cube. They do not see that the Cayley $E=V=F=\mathbb{N}_0$ tree $T_{\ge0}$ shows arrows ($E$dges) to same-colour arrowhead nodes ($V$ertices) with above each node an upwardly open, although narrowing, same-numbered site (sur$F$ace). They do not see that the upwardly unclosed site (sur$F$ace) above its cyclic root $\cabroot$ shows up as the regions to the left, to the right and below tree $T_{\ge0}$ (Figs.\ref{fig:tge0},\ref{fig:treeP1},\ref{fig:tree7pm1},\ref{fig:tree3m1},\ref{fig:tree5p1}).  Paths of serial arrows to successor nodes (Eqs.\ref{eq:dfsi},\ref{eq:dfSj}) in $T_{\ge0}$
show up in its inscribed 4-regular Cartesian half graph $C_H$ as transformed sequences of parallel, `random-access', arrows to same-numbered successor \textit{sites} (Eqs.\ref{eq:dfsipar},\ref{eq:dfSjpar},\ref{eq:xpar}). The Cartesian half graph provides a Hamiltonian tour to all $T_{\ge0}$ number pairs located on upper half coordinates ($x,y$) from and to the trivial root number pair $\cabroot$ at its origin (0,0) and cyclic root, using two out of four arrows from/to each coordinate $x\leftrightarrow x\pm1,y\leftrightarrow y\pm1$ (metamorphosis 2, Section \ref{sec:32eqs}, Figs.\ref{fig:hamiltonCo},\ref{fig:hamiltonBa},Theorem \ref{the:cartesian}).

The third metamorphosis of $T_{\ge0}$ is its automorphism graph \cite{RN11,RN20} Aut($T_{\ge0},[L,U]$). It is an infinite series of Plato's 3-regular cubes/earths and 4-regular octahedrons/airs, each with inscribed octahedrons and cubes. Each cube has a grandchild and grandparent cube, each octahedron has a grandchild and grandparent octahedron, giving as automorphism graph an infinite Jacob ladder of cube$\leftrightarrow$octahedron steps connecting the cube rail of connected grandparent$\leftrightarrow$grandchild cubes to the octahedron rail of connected grandparent$\leftrightarrow$grandchild octahedrons. Its transformed automorphism graph $\textup{Aut}(T_{\ge0},[L,U])$\; (Fig.\ref{fig:aut}), shows $1$ red $L$eftward  ladder with subtree $t_{0,\ge i}\;\tcr{\to}$ cotree $t_i$ steps and $a-1$ $U$pward ladders with subtree $t_{\ge j}\;\tcb{\to}$ cotree $T_j$ steps, defining by the isomorphism transformation formula presented by Max Dehn (1878-1952) \citep[p117]{RN19},\cite{peiferCayleyDehn}  `conjugate', `spurious influence', inner arrows in red $L$eftward and blue $U$pward cotrees as well as outer paths with a terminal conjugate arrow to $L$eftward cotrees (Figs.\ref{fig:ltree}bdf) via same-Least-Common-Multiple-periodicity (lcmp) $L$eftward subtrees (Figs.\ref{fig:ltree}a$\tcr{\to}$b, c$\tcr{\to}$d, e$\tcr{\to}$f; Eq.\ref{eq:dfsicon}) and $U$pward cotrees (Figs. \ref{fig:ltree}a$\tcb{\to}$b, c$\tcb{\to}$d, e$\tcb{\to}$f; Eq.\ref{eq:dfSjcon}) (metamorphosis 3, Sections \ref{sec:34eqs}, Figs.\ref{fig:hamiltonCo},\ref{fig:hamiltonBa},Theorem \ref{the:cartesian}).

Finite sums of infinite geometric series explain Plato's recount of Zeno's paradox: will Achilles surpass the Turtoise at his infinite horizon? Emmy Noether (1882-1935) considered multiple  \textrm{lcm}-periodicities of congruence classes of number nodes in isomorphic submodules (subtrees, cotrees) of non-commutative, non-confluent, `modules' (trees, graphs) \cite{RN18}. Multiple  \textrm{lcm}-periodicities give finite converging sums of infinite geometric progressions of periodic densities of congruence classes of same-generation successor numbers $s_{i}$ and $S_{j}$. A congruence class $[r]_m$ is an infinite arithmetic progression $\ell m+r$ for periods $\ell=0,1,2,\dots$ with periodic density $d=\frac{1}{m}$. The automorphic biperiodic Cayley tree Aut$(T_{\ge0},[L,U]$ (Fig.\ref{fig:aut}) shows one converging geometric series of $L$eftward cotree congruence class densities (Figs.\ref{fig:ltree}bdf$\dots$), and $a-1$ converging geometric series of successive upward cotree congruence class densities (Fig.\ref{fig:utree}bf$\dots$,d$\dots$). This provides the cotree density proof $d(T_{\ge0})=d(S_{\ge0})$  or refutation $d(T_{\ge0})<d(S_{\ge0})$ (Figs.\ref{fig:hamiltonCo},\ref{fig:hamiltonBa}, Eq.\ref{eq:cdp}) of a Hamiltonian tour provided by the Collatz-Conway a$n$+b conjecture  (Metamorphosis 4, Figs. \ref{fig:cc3p1},\ref{fig:cc3m1},\ref{fig:cc7pm1}). 

The `cotree congruence class tape' $T_{CCC}$ (Table \ref{fig:Gccc}) recombines equations from its automorphism graph and congruence class graph. In left-right order its $c=a+1$ columns are its infinite series of leftward subtrees (Eq.\ref{eq:lsuar}), leftward cotrees (Eq.\ref{eq:lcoar}), and $a-1$ upward cotrees (Eq.\ref{eq:ucoar}). Its red $L$eftward subtree and cotree columns share in each row $i=1,2,3,\dots$ the lcm$_{\textup{pq}}$ (Eq.\ref{eq:lcmpq}) (least common multiple periodicty/modulus)$10^i$ of  purple, either red $\tcr{\to0}$ or blue $\tcb{\to1}$ arrowhead numbers in $T_{\ge0}$: $U^{i=1,2,3,\dots}:[T_{\ge0}]_{\textup{lcm}_{pq}} \tcr{\to} [t_{0;\ge i}]_{\textup{lcm}_{pq} 10^i}$ and 
$UL^{i=1,2,3,\dots}:[T_{\ge0}]_{\textup{lcm}_{\textup{pq}}}\tcr{\to} [t_i]_{\textup{lcm}_{pq} 10^i}$ (Eq.\ref{eq:dfsi}). Its $a-2$ $U$pward cotree columns share in their rows  blue generation successors in blue cotrees:   
$LU^{j=1,\dots,a-2}:
[T_{\ge0}]_{\textup{lcm}_{\textup{pq}}}\tcb{\to} [T_j]_{\textup{lcm}_{pq} 1^{j+1}}$ (Eq.\ref{eq:dfSj}). 
The exponentially growing number of congruence classes in successive leftward cotrees reflects the parameter $A^+=a^2+a+1 $ (Eq.\ref{eq:Apm}) in the number-of-congruence-classes multiplier $(\textup{log}_2(c)-1)A^+$, amounting for Collatz $3n+1$ to $(2-1)(3^2+3+1)=13$ and for Bařina $7n\pm1$ to $(3-1)(7^2+7+1)=114$. The parameter $A^-=a^2-a-1$ (Eq.\ref{eq:Am}) in $(\textup{log}_2(c)-1)^2A^-$ gives the same number of congruence classes in successive upward cotrees, $1^2(3^2-3-1)=5$ for Collatz $3n+1$  and $(3-1)^2(7^2-7-1)=164$ for Bařina $7n\pm1$. The top row includes sums of $1+(a-1)=a$ geometric series of cotree densities, amounting to their cumulative row total $d(tT)$ of red and blue cotree densities (Eq.\ref{eq:cdp}), proving or refuting the $an+b$ conjecture by whether $d(tT)$ equals or subcedes the periodic density of purple, red or blue, branching numbers $d(S_{\ge0})$ (Eq.\ref{eq:dSge0}). 

Recombinations of equations also yield the infinite Turing-Leibniz tape $T_L$ with in its columns proven/ provable Collatz-Conway $a$n$_{\textup{odd}}$+b conjectures (Table \ref{fig:TLEE}) and in its rows equations with as left-hand side functions yielding  (transformed) Cayley tree features. 
The Turing-Leibniz tape  (Table \ref{fig:TLEE}) includes an infinite row of arrows to natural arrowhead numbers $\eta=\croot\tcr{0}\tcb{\to1}
\tok1\tcr{0}\tcb{\to}
1\tcb{1}\tok10\tcr{0}\tcb{\to}\cdots$ (Eqs.\ref{eq:sorS}-\ref{eq:bf7pm1}) generating the row of cyclic roots  $c=a\cdot1+b=2^{2^{\eta-1}+1}$ of $a\nmid$ a n{$_\textup{odd}$}+b numbers, except for $\eta=0$ generating the cyclic root $c=2$ of the Conway $1n+1$ function by $c=2^{\eta-1}=2$. Successive row numbers are denoted by $\eta$, as $n$ is already in use for provable/proven $n$=an$_\textup{odd}$+b conjectures in its columns. Columns $1,2,3,7,31,\dots$ are Euclid-Euler pairs of Mersenne prime $(1<n_x\le a) \nmid$ a$n_\textup{odd}$+b conjectures, after Marin Mersenne (1588-1648), paired with `perfect' numbers $2^{p-1}(2^p-1)=\frac{c}{2}(c-1)$, being numbers that are both the sum and the product of their divisors: for Collatz 3$n$+1: $c=a+1=4=2^2\Rightarrow\frac{c}{2}(c-1)=1\cdot2\cdot3=6=1+2+3;$ for Bařina 7$n\pm$1: $c=a+1=8=2^3\Rightarrow\frac{c}{2}(c-1)=1\cdot4\cdot7=1\cdot2\cdot14=28=1+2+4+7+14$. The periodicities $\frac{a+1}{2}a=\frac{c}{2}(c-1)$ are the periodicities of yellow odd number classes indivisible by $a$ (Figs.\ref{fig:gcol},\ref{fig:tree7pm1},\ref{fig:cc3p1},\ref{fig:cc7pm1}). Euclid ($\pm $300bC) observed pairs of primes $2^{p=a=c-1}-1$ and perfect numbers in his geometric 2-color figure \citep[Prop.36]{Euclid300}. By 2026 the search for higher primes has delivered $52$ $\sm{\textup{Mersenne prime}\\ \textup{perfect number}}$ pairs \cite{gimps}. 
Row $a=c-1$ in Table $T_L$ shows in columns $\log_2(c=\{2^1,2^2,2^3,2^7,2^{31},\dots\})=\{1,2,3,7,31,\dots\}$ Euclid-Euler cyclic roots $c=a\cdot1+1$ that are Mersenne primes $(1<n_x\le a)\nmid a\,$n$_\textup{odd}$+b, paired with perfect numbers.

The infinite rows of the Turing-Leibniz tape of provable/proven $a$n$_{\textup{odd}}$+b conjectures (Table \ref{fig:TLEE}) include a cyclic power-of-2 number root $c=a\cdot1+1=(c-1)\cdot2^q$, a lowest upward trunk number $U(c)=c^2$, an lcm-periodicity $\frac{c}{2}(c-1)$ of yellow odd classes indivisible by $a$ (a perfect number for Euclid-Euler columns $1,2,3,7,31,\dots$), an lcm-periodicity $\frac{c}{2}(c-1)^2$ of either red or blue, purple, numbers, an lcm-periodicity $(\log_2 c-1)\frac{c^3}{2}(c-1)^2$ of red and blue numbers, a `cotree density' Hamiltonian tour proof $d(tT)$ from and back to the trivial root $c$ (Eq.\ref{eq:cdp}) and a computationally reduced root-trunk proof $c^2<\frac{c}{2}(c-1)^2$ (Eq.\ref{eq:rtt}). 

Rahn, Sultanow, Henkel, Ghosh and Aberkane  \cite{rahnSultanow2021} picture dark-light projections of Collatz 3$n$+1 convergence to the trivial root. They report modular arithmetic switches between bases/periodicities $2$, $3$ and $4$ \citep[p18]{rahnSultanow2021}, interpreted as "pseudo-randomness" due to "symmetry breaking properties" by ``the chaoticity of conversions between bases 2 and 3" \citep[p27-29]{rahnSultanow2021}. The non-chaotic Turing-Leibniz-tape $T_L$  of provable/proven 3-colour $an$+b conjectures (Table \ref{fig:TLEE}) integrates $c$-dependent cyclic roots $\croot c=a\cdot1+1$ generated from the natural numbers by $2^{2^{\eta-1}+1}$, $a=c-1$ parameters, with $\psi=(c-2)$-number base trees of geometric series of cotree densities. In Hofferbert's music metaphor \cite{RN31}, biperiodic $c$-expressible parameters show up in an infinite number of $a\nmid n=\;$an$_{\textup{odd}}$+b fugues, tessellating \textit{infinite numbers} of cotree \textit{periodicities}.

\subsection{Ready for such problems}

Paul Erd\H{o}s' (1913-1996) stated:
``mathematics is not yet ready for such problems" (op.cit.\cite{lagarias00}). Collatz' expectation is met that picturing the $3n_\textup{odd}+1$ function as a directed graph connects ``elementary" graph theory and number theory \cite{RN9,RN34}. Picturing directed 3-colour graphs shows graph-theoretic Hamiltonian tour proofs and root-trunk-proofs of provable/proven Collatz-Conway $a\,$n$_\textup{odd}$+b conjectures as infinite Turing-Leibniz tape $T_L$ rows, including Euclid-Euler tape columns. ‘Number theory’ mirrors ‘graph theory’, even though at the time of writing computer (proof assistant) languages (e.g. C(++), Prolog, Mathematica, Python, R, Lean) lack ‘types’ for congruence classes as infinite arithmetic progressions, numbers as cyclic roots, or arrow-colour functions to same-colour arrowhead nodes.
Number theoretic specifications of proofs that are also graph-theoretic proofs provide Lean proof certificates for Collatz-Conway a$n$+b conjectures on the infinite Turing-Leibniz tape $T_L$, including Collatz 3$n$+1 and Barina 7$n\pm$1 (Tables \ref{fig:TLEE},  \ref{tab:cdp3p1}, \ref{tab:cdp7pm1}).

\newpage
\section{Cayley tree $T_{\ge0}$ of a Collatz-Conway $a\nmid n=an_{\textrm{odd}}+b$ function}

\subsection {Eq. series 1: Elementary Graph Theory 1736-1936}
\label{sec:31eqs}
\begin{theorem} Every Collatz-Conway a$\nmid$n=\textup{a}n$_{\textup{odd}}+\textup{b}$ function
is a:
\label{the:caytree}
\end{theorem}
\vspace{-0.8cm}
\begin{align}
3&\textrm{-regular 3-colour Cayley graph } G_{\ge0}\supseteq T_{\ge0}\textrm{ at least comprising the } \nonumber \\ 
&\textrm{left-right-asymmetric, non-confluent, non-commutative,  Cayley tree} &T_{\ge0}\subseteq G_{\ge0}\left(S_{\ge0},\left[\sm{\textcolor{red}{0}\\L},\sm{\textcolor{blue}{1}\\U}\right]\right) \label{eq:tleg}\\ 
&E=V=F=\mathbb{N}_0 \textrm{ with the `trivial' cyclic root }&\rotatebox[origin=c]{270}{\(\textcolor{red}{\circlearrowright}\)}\,\sm{\textcolor{red}{0}\\c=a\cdot1+b} \label{eq:Croot}\\ 
&\textrm{for which \textit{q} iterations of \textit{g} yield its lowest trunk number }U(c) &\;\;\textcolor{blue}{\to}\sm{\textcolor{blue}{1}\\U(c)=g^q(c)=c\cdot2^q}\label{eq:Trunk}\\
&\textrm{from which grows its crown of number pairs }\sm{y\\n}\textrm{ of:}& \sm{y \textrm{ binary root path number }\quad\quad\;\\n=a\cdot n_{odd}+b;\; a\nmid n; f(n)=\frac{n-b}{a}\in[1]_2} \label{eq:yn} \\
    &\textrm{- binary root paths $y\in(c-2)$-base number tree growing by:} \nonumber\\
&\quad\textrm{numeric functions: }\; \tcr{0}:y\textcolor{red}{\to}\,y\,\cdot10+0=y\tcr{0} \textrm{; }\;\,\tcb{1}:y\textcolor{blue}{\to}\,y\,\cdot10+1=y\textcolor{blue}{1}  \notag \\
&\quad\textrm{labeled as: }\quad\;\;\,\quad\quad g:y\textcolor{red}{\to}\,y\,\cdot10+0=y\tcr{0} \textrm{; }\;f:y\textcolor{blue}{\to}\,y\,\cdot10+1=y\textcolor{blue}{1}  \notag \\
&\textrm{- even/red/\textit{L}eftward \& odd/blue/\textit{U}pward function digits }&
\sm{\textcolor{red}{0}\\L}=
\sm{\textcolor{red}{0}^{\textcolor{white}{-}}\\fg^p\in s_{\ge0}},\;\;
\sm{\textcolor{blue}{1}_{\textcolor{white}{-}}\\U^{\textcolor{white}{-}}}=
\sm{\textcolor{blue}{1}\\g^q\in S_{\ge1}} \label{eq:rootPairs}\\
&\quad\textrm{$p$ and $q$: \;powers-of-2,\; $g$-path-lengths in transformed graph $G_{CC}$ } \notag
\end{align}

\normalsize

\noindent The first subsection \ref{sec:31eqs} on Eqs.\ref{eq:tleg}-\ref{eq:rootPairs} bears the name elementary \textit{Graph theory 1736-1936}, after the source textbook by Biggs, Lloyd $\&$ Wilson anticipating the 1976 computer algorithm proof by Appel and Haken that four colours suffice to distinguish adjacent sites \citep[6A-6D,7A-7C,9A-9D]{RN3245}. 1736 is the year in which Leonhard Euler (1707-1783), regarded as the first pioneer of graph theory  \citep[1A,ch.1]{RN3245}, published the first application of the `analysis situs', proposed by Gottfried Wilhelm Leibniz (1646-1716), to the odd number of bridges/edges between city quarters/sites of Königsberg, implying that at least one visited Königsberg site could not be finally left without using an already used bridge/edge twice. 

1936 is the year in which Alan Turing (1912-1954) conceived an $a$(utomatic)-machine, a ladder graph, of binary numbers on one of its rails paired with computable numbers \cite{turing1936} on the other, labeled a \textit{Turing tape} by his doctoral adviser Alonzo Church (1903-1995), here labeled a \textit{Turing-Leibniz} (multi)tape $T_L$, as already Leibniz construed a `calculus ratiocinator' pairing binary numbers $y$ for calculations with computable numbers $n$ in pairs $\sm{y\\n}$ (Eq.\ref{eq:yn}), employing computable decimal input and output numbers \citep[ch.4.2]{lleibniz}. 

Arthur Cayley (1821–1895), regarded as the third pioneer of graph theory \citep[ch.3,3A,6A]{RN3245}, invented trees and colour graphs \citep[3A,6A]{RN3245}\cite{RN3242} as observable, visible, `analytic forms'. As the publishers of his days did not print colours, Cayley explained that analytic forms "should be in the first instance drawn with the arrows" \citep[p181]{RN3242} and that (on white paper) ``the continuous lines represent black lines, and the broken lines, red lines''. Cayley published `arrows' with a same-colour arrowhead node ($V_{\textup{ertex}}$) as undirected `lines' ($E_{\textup{edges}}$) with black or broken arrow labels. Black dots do not discern black and broken (red) arrowhead nodes.

In 1883, Georg Cantor (1845-1918), referring to Spinoza (1632-1677), Locke (1632-1704) and Leibniz (1645-1716)  \citep[\S5]{cantor1883bin}, defines the binary number tree $T_{\ge0}$ of non-negative binary numbers $\mathbb{N}_0=0,1,2,\dots$ filling congruence classes $[0,1]_2\stackrel{\textup{bin}}{=}[0,1]_{10}$ (Figs.\ref{fig:tge0},\ref{fig:treeP1},\ref{fig:tree7pm1},\ref{fig:tree3m1},\ref{fig:tree5p1}) by iterations of two biperiodic depth-first binary functions \citep[\S6]{cantor1883bin}, here coloured red and blue, and numbered $\tcr{0}:$ and $\tcb{1}:\,$, or labeled $g:$ and $f:\,$, that lengthen a binary number  $y_{\ge0}\in\mathbb{N}_0$ with a terminal arrow function digit to its red arrowhead child $\tcr{\to0}$ (Eq.\ref{eq:redChild}) and its blue arrowhead child $\tcb\to1$  (Eq.\ref{eq:blueChild}):
\begin{align}
& \tcr{0}:(=g:) y\tcr{\to} y\cdot10+0=y\tcr{0},  \textrm{ implying redundant leading } \tcr{0} \textrm{ digits in cyclic root }\croot\tcr{0}=0\cdot10+0
   \notag\\
 & \tcb{1}:(=f:) y\tcb{\to} y\cdot10+1=y\tcb{1} \tag{\ref{eq:yn}}
\end{align}

\noindent Recursive `number as function' iterations $\tcr{0}:(=g:)$ and  $\tcb{1}:(=f:)$ generate the binary number tree of $2^{k-1}$ non-negative binary root path numbers of $k$ arrow-function digits at successive breadth-first, or `diagonal', levels 
$k=\croot\sm{1\\1+2^{1-1}=2}\tok\sm{2\\2^{2-1}=2}\tok\sm{3\\2^{3-1}=4}\tok\sm{4\\2^{4-1}=8}\tok\cdots$ (Eq.\ref{eq:sorS}). Level $k=1$ comprises not just $2^1-1=1$, but $1+(2^1-1)=2$ nodes: the cyclic root $g:\croot\tcr{0}\to10\cdot0+0=\croot\tcr{0}$ regenerates itself (Eq.\ref{eq:Croot}) and generates the lowest upward trunk node $f:\tcr{0}\to10\cdot0+1=\tcb{1}$. Cantor did not picture the resulting biperiodic, 3-colour, one-trunk, binary number tree growing from its cyclic root that is paired in Cayley tree $T_{\ge0}$ to Collatz-Conway $a\nmid$ a$n_{\textup{odd}}$+b numbers (Eqs.\ref{eq:Trunk}, Eq.\ref{eq:sorS}, Figs.\ref{fig:tge0}a,\ref{fig:treeP1},\ref{fig:tree7pm1},\ref{fig:tree3m1},\ref{fig:tree5p1}). The notation of odd and even congruence 
classes, being infinite arithmetic progressions $0,2,4,\dots$ and $1,3,5,\dots$, 
as $[0]_2$ and $[1]_2$ and of binary numbers as $[0,1]_2\stackrel{\textup{bin}}{=}[0,1]_{10}$ dates from after 1936 (e.g.\cite{RN3231}), just like pictures of an infinitely high one-trunk tree contracted to the open interval (0,1) (Fig.\ref{fig:tge0}b) \cite{RN3247} instead of as a derooted shrub. 

The 1924 article ``On the Infinite" by David Hilbert (1862–1943) \cite{RN16} explains the Cantor functions $\tcr{0}:(=f:)$ and $\tcb{1}:(=g:)$ that identify the binary numbers $N_0=0\tcb{\to1}\tok{2}\cdots$ in tree $T_{\ge0}$ (Eq.\ref{eq:rootPairs}), by their metamorphosis as the fully booked infinite Grand Hotel capable of hosting an infinite row of new guests. 

\begin{corollary}
\label{the:hilberthot}
Hilbert's fully booked grand hotel capable of hosting an infinite rows of new guests is a metamorphosis of the Cayley tree $T_{\ge0}$ with binary root paths $y$ defined by iterations of the Cantor-functions $\tcr{0}:(=f:)y\to10\cdot y+0=y\tcr{0}$ and $\tcb{1}:(=g:)y\to10\cdot y+1=y\tcb{1}$. 

\noindent\textit{Proof}: First, the infinite row of new guests pairs up with the current guests in their rooms with a root path $y$ of binary digits $[\tcr{0},\tcb{1}]$. Next, each current guest clears its room and climbs with its new guest to the next higher breadth first level $\tok$. Finally, the new guest and the current guest take the cleared rooms of respectively the current guest's blue $U$pward child $\tcb{1}:y\to10\cdot y+1=y\tcb{1}$ (Eq.\ref{eq:blueChild}) and the current guest's red $L$eftward child $\tcr{0}:y\to10\cdot y+0=y\tcr{0}$ (Eq.\ref{eq:redChild}), which recycles the hotel manager $\tcr{0}$ as $\tcr{0}:0\to10\cdot 0+0=0\tcr{0}=\tcr{0}.$ \qed
\end{corollary}

The Cayley tree $T_{\ge0}$ of even Collatz-Conway numbers $a\nmid n\,=\,a\,n_{\textup{odd}}\,+\,b$ with cyclic root $\cabroot$ (Eq.\ref{eq:Croot}) generated by combinations of the $L$eftward and $U$pward functions, $L=fg^p$ and $U=g^q$, paired with binary root path numbers $y$ of red and blue arrow function digits $\tcr{\to0}$ and $\tcb\to1$ generated by the red and blue binary Cantor-functions $g:(=\tcr{0}:)$ and $f:(\tcb{1}:)$, consists of paired number nodes $\sm{y\\n}$ (Eq.\ref{eq:yn}) of even/red/\textit{L}eftward arrow function digits
$\sm{\textcolor{red}{0}_{\textcolor{white}{-}}\\L^{\textcolor{white}{-}}}=
\sm{\textcolor{red}{0}^{\textcolor{white}{-}}\\fg^p\in s_{\ge0}}$ and
odd/blue/\textit{U}pward arrow function digits
$\sm{\textcolor{blue}{1}_{\textcolor{white}{-}}\\U^{\textcolor{white}{-}}}=
\sm{\textcolor{blue}{1}\\g^q\in S_{\ge1}}$ (Eq.\ref{eq:rootPairs}), with powers-of-2 $p$ and $q$ that are $g$-path-lengths in its transformed congruence class graph $G_{CC}$ (Figs.\ref{fig:cc3p1},\ref{fig:cc3m1},\ref{fig:cc7pm1}). The paired numbers $\sm{y\\n}$ are categorized by their red $L$eftward generation $s_i\subset s_{\ge0}$ or blue $U$pward generation $S_j\subset S_{\ge1}$ (Eq.\ref{eq:rootPairs},Figs Figs.\ref{fig:tge0},\ref{fig:treeP1},\ref{fig:tree7pm1},\ref{fig:tree3m1},\ref{fig:tree5p1}).

A Collatz-Conway a$n$+b conjecture is proven if the Cayley tree $T_{\ge0}\subseteq G_{\ge0}$ with cyclic root $\cabroot$ is the entire Cayley graph  $T_{\ge0}= G_{\ge0}$ and refuted if its only one of the trees in a Cayley forest $G_{\ge0}\supseteq T_{\ge0}$. Collatz-Conway $n$=an${_\textup{odd}}$+b functions show up as infinite Cayley graphs   (Figs.\ref{fig:tge0},\ref{fig:treeP1},\ref{fig:tree7pm1},\ref{fig:tree3m1},\ref{fig:tree5p1}) that are (Eq.\ref{eq:tleg}):
\begin{itemize}
    \item 3-colour graphs (Eq.\ref{eq:tleg}), showing, next to the colour of same-coloured `sites' from Leibniz' `analysis situs', usually the colour of Plato's wall, sand, white paper or of the blackboard/screen, two colours of arrows and arrowhead nodes, often red and blue like here, or blue and brown \cite{sanderson}. 
    \item 3-regular graphs by the definition of Petersen \citep[11A]{RN3245}\cite{RN3242}\citep[1029]{RN23} (Eq.\ref{eq:tleg}). Each node has one red arrow to its red $L$eftward arrowhead child, one blue arrow to its $U$pward arrowhead child, and one, either red or blue, parental arrow (Eq.\ref{eq:tleg}), e.g. the red, parental cyclic root arrow $\cabroot$ (Eq.\ref{eq:Croot}). 
    \item base-$(c-2)$ arrow number functions, thus Bařina $7n\pm1$ rests also on $5: n \to n\cdot10+5=y5$ (Fig.\ref{fig:tree7pm1}).
\item $E=V=F=\mathbb{N}_0$ infinite number trees $T_{\ge0}$ with same-colour arrows ($E_\textup{dges}$) and arrowhead nodes ($V_\textup{ertices}$), and same-numbered nodes and upwardly open sites above them ($\text{\footnotesize{\textup{sur}}}F_\textup{aces}$) (Eq.\ref{eq:tleg}). Sites are bordered to the left by a serial path of blue $U$pward arrows and to the right by a serial path of red $L$eftward arrows. The upwardly open site of the cyclic root number pair $\cabroot$ to the left and right and below Cayley tree $T_{\ge0}$ is bordered to the left by the blue ``upward trunk'' and to the right by the red ``greedy branch''. 
    Euler's formula $V+F-E=2$ for closed solid bodies (solida) \cite{euler1758solid} does not apply.
\end{itemize}

\noindent The 3-regular, hence odd-regular, Cayley tree $T_{\ge0}\subseteq G_{\ge0}$ prevents an `Eulerian tour' as not a single visited quarter/site can be finally left without using an already used bridge/edge. Its transformation to a Cartesian half graph $C_H$ in the next subsection \ref{sec:32eqs} provides a 'Hamiltonian tour' from and back to the trivial root $\cabroot$ of all $T_{\ge0}$ number pairs.

\subsubsection{Wondrously ordered a$n$+b numbers with well-ordered $y$-axis projections}

Theorem \ref{the:fractal} states, first, that a geometric series transformation contracts the infinite Cayley a$n$+b tree $T_{\ge0}\subseteq G_{\ge0}$ to a fractal binary tree in the open interval (0,1) (Fig.\ref{fig:tge0}b), as pioneered by Mandelbrot and Frame \cite{RN3247}.
For ever higher breadth-first levels $k$ the fractal tree with upwardly open, although narrowing, sites becomes near the upper limit $1$ of the open interval $(0,1)$ a space-filling function with infinitesimal small root path distances $(x\pm\delta x,y\pm\delta y)$ to adjacent root paths, as is to be expected already from $2^{20}=1\,048\,576$ root paths for $k=20$ breadth-first iterations (Fig.\ref{fig:tge0}b). Subsection \ref{sec:35eqs} relies on similar geometric series contractions.

Second, the theorem shows the metamorphosis of well-ordered, breadth first ordered root paths to wondrously ordered $a\nmid an{_\textup{odd}}+b$ numbers into well-ordered $a\nmid an{_\textup{odd}}+b$ numbers led to by wondrously ordered root paths. Both well-orderings are required in the upcoming subsections. Without loss of generality, the proof uses numerical examples from the Collatz $3n+1$ function.

\begin{theorem} 
\label{the:fractal}
    The breadth-first levels $k=1\tok2\tok3\tok\cdots$ in the fractal representation of the Cayley graph \textup{(Fig.\ref{fig:tge0}b)}, holding for $k=1$ the cyclic root $\sm{\tcr{0}\\c=3\cdot1+1=4}$ and the lowest upward trunk number $\sm{\tcb{1}\\U(c)=c g^q=c\cdot2^q=4\cdot2^2=16=c^2}$, and for   $k>1$ a number of $2^{k-1}$ well-ordered binary root path numbers of red and blue binary digits $\sm{\tcr{0}\\L}$ and
    $\sm{\tcb{1}\\U}$ to wondrously ordered Collatz 3$\nmid$3$n{_\textup{odd}}$+1 numbers:
\[     \rotatebox[origin=c]{270}{\(\textcolor{red}{\circlearrowright}\)}
\sm{\textcolor{red}{0}^0\\\sm{4\\\sm{L\\s_0}}}
\textcolor{blue}{\to}\sm{\textcolor{blue}{1}^1\\\sm{16\\\sm{LU\\S_1}}}
\tok
\sm{1\textcolor{red}{0}^1\\
\sm{40\\\sm{UL^1\\s_1}}}
\textcolor{blue}{\to}
\sm{1\textcolor{blue}{1}^1\\\sm{256\\\sm{LU^2\\S_2}}}
\tok\
\sm{1\textcolor{red}{0}^2\\\sm{52\\\sm{UL^2\\s_2}}}
\textcolor{blue}{\to}
    \sm{10\textcolor{blue}{1}^1\\\sm{160\\\sm{ULU^1\\S_1}}}
\textcolor{red}{\to}
\sm{11\textcolor{red}{0}^1\\\sm{340\\\sm{U^2L^1\\s_1}}}
\textcolor{blue}{\to}
\sm{\textcolor{blue}{1}^3\\\sm{1024\\\sm{LU^3\\S_3}}}
\tok
\sm{1\textcolor{red}{0}^3\\
\sm{34\\\sm{UL^3\\s_3}}}
\textcolor{blue}{\to}\cdots\quad\quad\quad 
\eqref{eq:cywell}
\]

\noindent project to the Cartesian $y$-axis well-ordered a$n$+b numbers with wondrously ordered root paths \textup{(Eqs.\ref{eq:lowb}, \ref{eq:frac0})}: $\sm{y\\n}\in\sm{[0,1]_2\\ [4,16]_{18}}=\sm{0\\4},\sm{1\\16},\sm{10^4\\22},\sm{10^3\\34},\sm{10^1\\40},\sm{10^2\\52},\cdots$, as a geometric series contracted to the open interval (0,1). 
\end{theorem}

\noindent\textit{Proof/Explanation}.
The length of the red and blue arrow  $L:4\tcr{\to}4$ and $U:4\tcb{\to}16$ at breadth-first level $k=1$, is set to $a=1/3^2$. This length is contracted at each higher breadth-first level $k=2,3,4,\cdots$ by multiplication with its complement $(1-a)=2^3/3^2$ raised to the power ${k-1}$ (Eq.\ref{eq:contraction}). This amounts to a geometric series with a sum coming arbitrarily close to the limit $\sum_{k=1}^{\infty} a (1-a)^{k-1} = 1$ (Eq.\ref{eq:contractionSum}) as the finite horizon of the open interval $(0,1)$. Eq.\ref{eq:contractionSum} shows the well-known geometric sum formula $s=a/(1-r)$ as $a/(1-a(1-a)/a)$. 
\begin{equation}
\label{eq:contraction}
\textrm{Arrow length at each breadth-first level } k=1,2,3,\cdots:\;\; a_k\;=\;a (1 - a)^{k-1},\textrm{ in Fig.\ref{fig:tge0}b:}\;\; a_1=1/3^2
\end{equation}
\begin{equation}
\label{eq:contractionSum}
\textrm{Geometric sum arrow lengths at levels }k=1,2,3,\cdots: \;\; 
\sum_{k=1}^{\infty} a (1 - a)^{k-1}  =\frac{a}{1-\frac{a(1-a)}{a}} =\frac{a}{a}=1 
\end{equation}
The arrow length $a_k$ at breadth-first level $k$ is the hypotenuse of a rectangular triangle, with sines $\alpha_{q\in \{2,4\}}$ (Eq.\ref{eq:angleq}) and $\alpha_{p\in \{1,2,3,4\}}$ (Eq.\ref{eq:anglep}) of the angle of the hypotenuse with the horizontal axis $x$. 
\begin{align}
&\sm{\textrm{Breadth-first level } k \\ \textrm{Length hypotenuse } a_k \textrm{ given } a_1=1/9 : }\quad\sm{1\\0.111},\sm{2\\0.998},\sm{3\\0.088},\sm{4\\0.078},\sm{5\\0.069}, \cdots,\sm{18\\0.015},\cdots,\sm{100\\3.736\cdot10^{-7}},\cdots
\label{eq:fracmag} 
\end{align} 
The sines $\alpha_{q}$ (Eq.\ref{eq:angleq}) and $\alpha_{p}$ (Eq.\ref{eq:anglep}) depend on the opposite sides of angles in their nominators  which are, relative to each parent number $n$ on a logarithmic scale with base $2$, the increase or decrease of $n$'s (blue) upward child $U(n)=g^q(n)=n\cdot2^q$ and of $n's$ (red) leftward child $L(n)=fg^p(n)=(n-1)/3\cdot2^p$, divided by the maximum increase on this logarithic scale. The maximum increase, obtained for any (blue) upward child of a parent $n\in[16]_{18}$ (Defs.\ref{def:Uw},\ref{def:Lw}), gives a division of the increase or decrease by $\log_2(2^4)=4$, or a multiplication by $\frac{1}{4}$ (Eqs.\ref{eq:angleq},\ref{eq:anglep}). The possible periodicities, corresponding to all possible powers $p\in \{1,2,3,4\}$ and $q \in \{2,4\}$ in the sine equations $\alpha_{q\in \{2,4\}}$ (Eq.\ref{eq:angleq}) and $\alpha_{p\in \{1,2,3,4\}}$ (Eq.\ref{eq:anglep}) define the upward function $U=g^p$ and leftward function $L=fg^p$ (Defs.\ref{def:Uw},\ref{def:Lw}) in Section \ref{sec:cc}. 
\begin{align} 
\label{eq:angleq}
 \sin{\alpha_{q\in \{2,4\}}} &:&\quad\{\frac{1}{4}\log_2(2^2)\;=\;\frac{1}{2},&\quad\frac{1}{4}\log_2(2^4)\;=\;1 \} \notag \\
\;\;\textrm{Degree angle } \alpha_{q\in \{2,4\}}&:&\quad\{\;30^{\circ},&\quad90^{\circ}\;\}
\end{align}
The triangle sides opposite to the angles of the triangle hypotenuses with the horizontal axis, measure the logarithmic length increases of (red) leftward children $L:n\to L(n)$ relative to $n$, converging to $4$ angles. 
\begin{align} 
\label{eq:anglep}
 \sin{\alpha_{p\in \{1,2,3,4\}}} &:\{\;\;\approx\frac{1}{4}\log_2 2^1/3,\;
 \;\approx\frac{1}{4}\log_2 2^2/3,
 \;\approx\frac{1}{4}\log_2 2^3/3,
 \;\approx\frac{1}{4}\log_2 2^4/3\;\} \notag \\ 
\;\textrm{Angle } \alpha_{p\in \{1,2,3,4\}} \textrm{ converges to}&:\{\;
\;\approx193^{\circ}  (\approx -13^{\circ} ),
\;\approx171^{\circ},
\,\quad\quad\approx158^{\circ},
\;\;\;\;\quad\quad\approx148^{\circ}\;\}
\end{align}
For the Collatz $3n+1$ function the first wondrously ordered root paths $y$ paired with well-ordered numbers $n=3{\textup{n}_\textup{odd}}+1$ numbers (Eq.\ref{eq:lowb}) contracted to the open interval $(0,1)$\; (Eq.\ref{eq:frac0}) are: 
\begin{align}
&\textrm{Wondrously ordered root paths:}\quad &\sm{\textcolor{red}{0}\\4}\quad\;,&\quad\sm{\textcolor{blue}{1}\\16},&\sm{\textcolor{blue}{1}\textcolor{red}{0000}\\22}\;,&\;\;\;\sm{\textcolor{blue}{1}\textcolor{red}{000}\\34}\quad\quad&\quad\sm{\textcolor{blue}{1}\textcolor{red}{0}\\40},\quad&\quad\sm{\textcolor{blue}{1}\textcolor{red}{00}\\52},&\cdots& \label{eq:lowb} \\
&\textrm{Well-ordered $y$-projections (0,1):}\quad
&\sm{4\\0.11111},\;&\sm{16\\0.14053},&\sm{22\\0.14701},&\;\sm{34\\0.16117},&\sm{40\\0.16627},&\;\sm{52\\0.17474},&\cdots&
\label{eq:frac0}\qed
\end{align}
\subsubsection{Cayley tree Conway 1$n$+1 conjecture: division by 1, no Conway fraction game, self-evident}
\begin{figure}[h!]
\caption{Cayley graph $T_{\ge0}$ of the Conway $1n+1$ function}
\label{fig:treeP1}
\begin{center}
\includegraphics[width=0.87\textwidth]{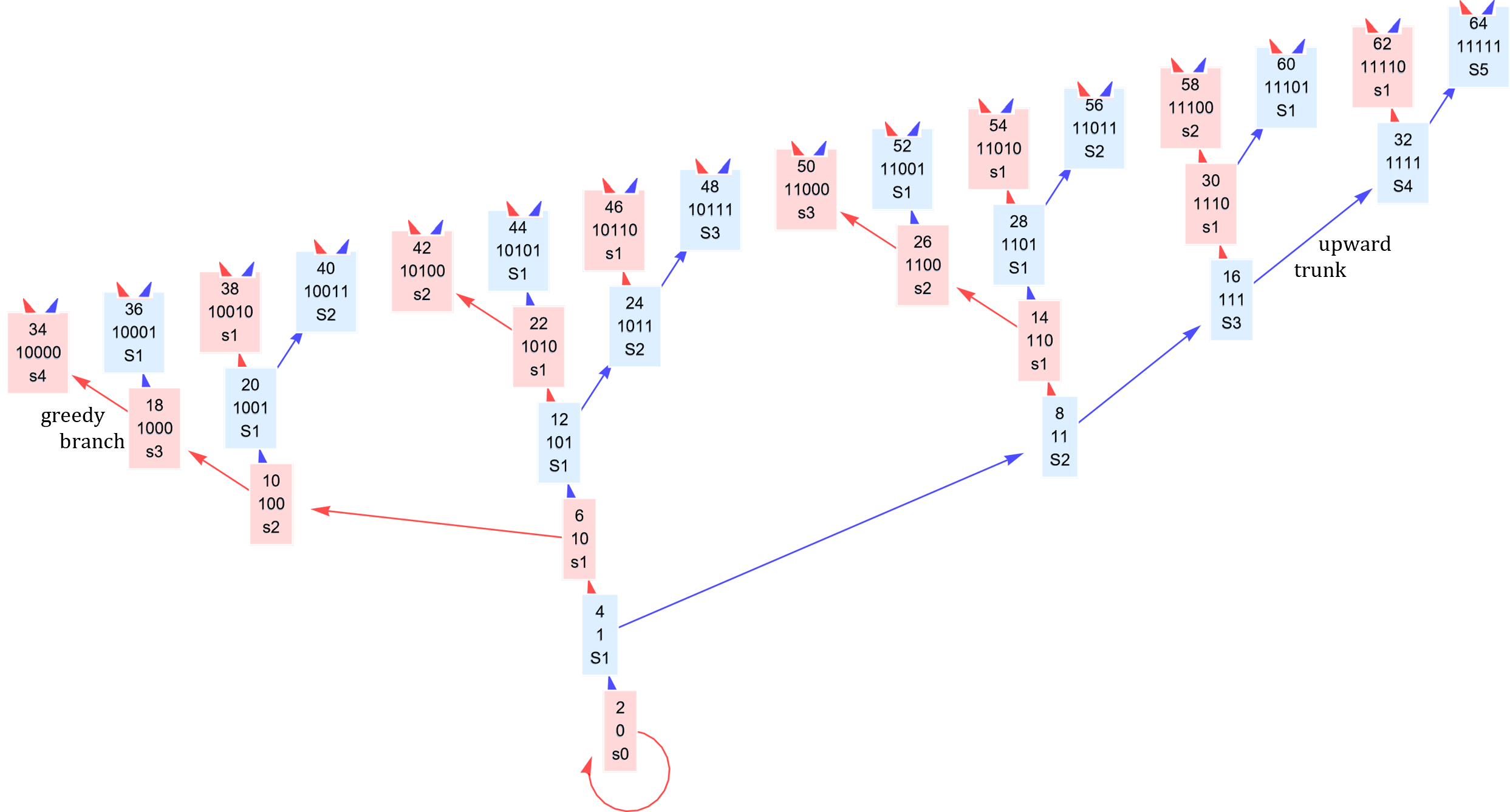}
\end{center}
\small{\textbf{Legend:} The Conway $1n+1$ function generates from its cyclic root $\rotatebox[origin=c]{270}{\(\textcolor{red}{\circlearrowright}\)}\,\sm{\textcolor{red}{0}\\c=1\cdot1+1}$ the breadth-first even number pairs line $\sm{y=n/2-1\\n=2(y+1)}$ with alternatingly a red $L$eftward child  $L:\sm{y\\n}\textcolor{red}{\to}\sm{y\textcolor{red}{0}\\2n-2=y}$, e.g. $L:\sm{100\\10}\textcolor{red}{\to}\sm{100\textcolor{red}{0}\\18}$, and a blue $U$pward child $U:\sm{y\\n}\textcolor{blue}{\to}\sm{y\textcolor{blue}{1}\\2n=y+2}$, e.g. $U:\sm{100\\10}\textcolor{blue}{\to}\sm{100\textcolor{blue}{1}\\20}$}.
\end{figure}

The Conway $1n+1$ function (Def.\ref{def:treeP1}, Fig.\ref{fig:tree5p1}) giving even numbers $n=1\cdot n_{\textup{odd}}+1$ does not give rise to a `fraction game' by the function $f:n\to (n-b)/a=(n-1)/1=n-1$ \cite{fractran,RN3209}, as $n{_\textup{even}}-1$ is invariably an odd integer number, implying non-crossing root paths. The $L$eftward function $L=fg^1:n\tcr{\to}(n-1)/1\cdot2^1=2n-2\in[2]_4$ yields red, Leftward successors $[2]_4$. The $U$pward function yields blue $U$pward successors $U=g^p=g^1=g:n\in[0]_2\tcb{\to}2n\in[0]_4$. The two functions provide a breadth-first traversal of its Cayley tree passing all well-ordered even 1$n$+1 numbers $n_{>0}\in[0,2]_4$ (Fig.\ref{fig:treeP1}).

\begin{definition} The proven $1n+1$ conjecture 
\label{def:treeP1}
\begin{align*}
f^{-1}:n\to 1n+1 \textrm{ for $n\ge1$ from odd class }n\in[1]_2; &&f:\;n&&\in[0]_2&&\to\;(n-1)/1\in[1]_2 &&\\
&&L=fg^{p=1}:\;n&&\in[0,2]_4&&\tcr{\to}\;\;(n-1)2^1\in[2]_4 &&\\
g^{-1}:n\to n/2 \;\;\;\textrm{ for $n\ge1$ from even class }n\in[0]_2;&& U=g^{q=1}=g:\;n&&\in[0,2]_4&&\tcb{\to}\quad\quad\quad\;2n\in [0]_4 &&
\end{align*}
\begin{align*}
\quad\textrm{Even branching numbers } S_{\ge0}=[0,2]_4\stackrel{\textup{bin}}{=}[0,10]_{100};\quad [L,U]\textrm{-output classes}=\{[2]_4,[0]_4\}\stackrel{\textup{bin}}{=}[10,0]_{100} 
\end{align*}
\end{definition}

\subsubsection{Cayley tree Bařina $7n\pm1$: second Mersenne prime conjecture, to be proven}

\begin{figure}[hbt!]
\caption{Cayley graph $T_{\ge0}$ of the Bařina $7n\pm1$ function}
\label{fig:tree7pm1}
\begin{center}
\includegraphics[width=0.87\textwidth]{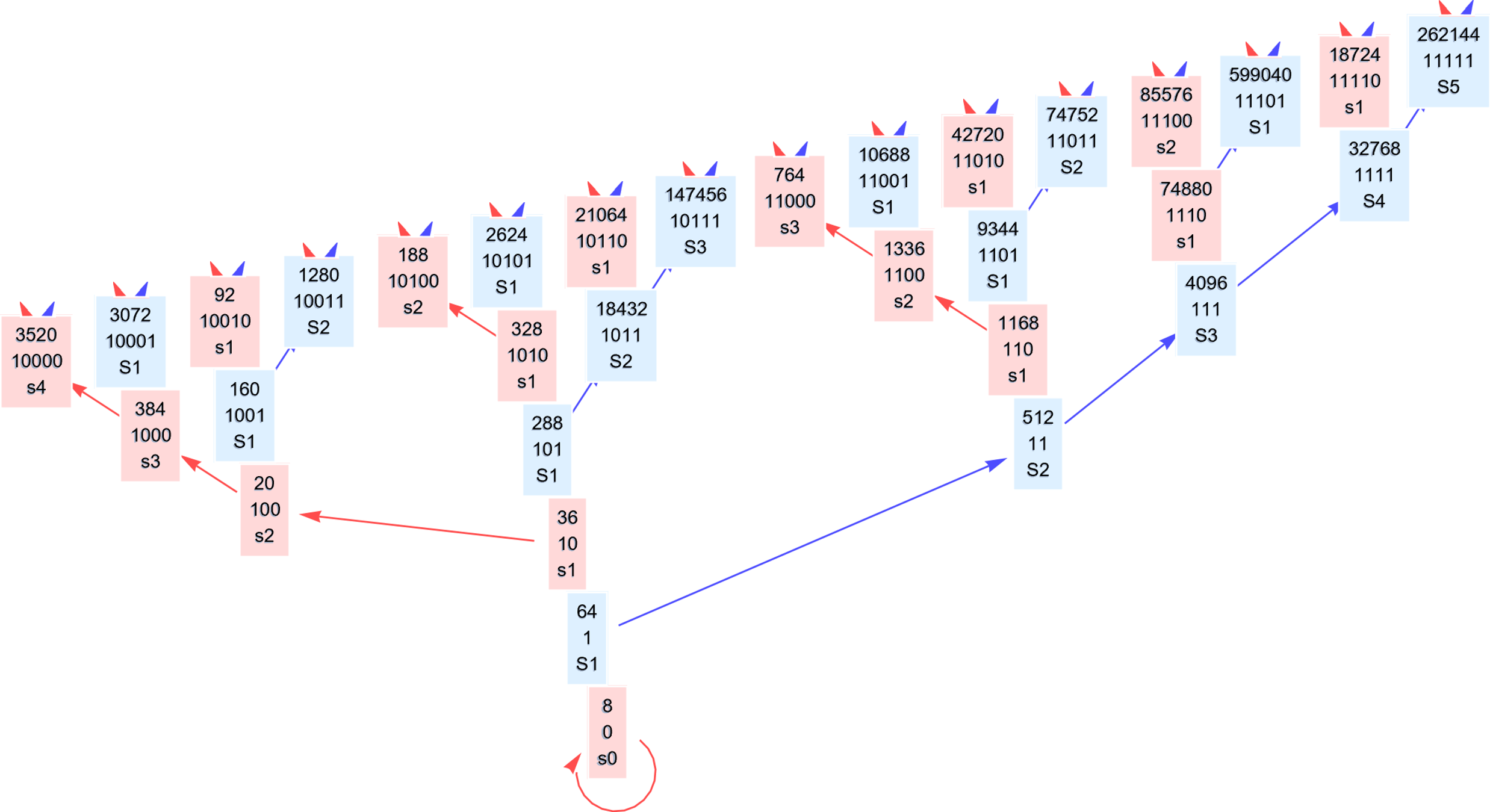}
\end{center}
\end{figure}

The $7n\pm1$ conjecture (Def.\ref{def:f7pm1},Fig.\ref{fig:tree7pm1}), proposed by David Bařina \cite{RN102} as a \textit{close relative} of the Collatz $3n+1$ function, has $v=2$ different subfunctions (Eq.\ref{eq:vnum}) for the odd subclasses $[1]_4,[3]_4\stackrel{\textup{bin}}{=}[1,11]_{100}$, giving $(a-1)v=6\cdot 2=12$ purple classes $[8,36,\cdots\ 160,188]_{2^v7^2=196}$ for which $a\nmid(n\pm1)/7$ (Eq.\ref{eq:dSge0}, Def.\ref{def:f7pm1}).
\noindent 
\begin{definition} The Bařina $7n\pm1$ function
\label{def:f7pm1}
\begin{align*}
f^{-1}:n&\to7n+1 \textrm{ for the odd subclass }[1]_4;&\quad\quad f:n\in[8]_{28} &\to (n-1)/7\in [1]_4 \\
f^{-1}:n&\to7n-1 \textrm{ for the odd subclass }[3]_4;&\quad\quad f:n\in[20]_{28} &\to (n+1)/7\in [3]_4 \\
g^{-1}:n&\to n/2 \textrm{ for the even class }[0]_2;&\quad\quad g:n\in[0]_2 &\to 2n\in [0]_2=[0,2]_4 
\end{align*}\vspace{-0.8cm}
\begin{align*}
\quad\quad(a-1)v=12 \textrm{ branching }&\textrm{classes }\textrm{ for which }7\nmid (n\pm 1)/7:\;\; &S_{\ge0}= &{\left[ \begin{array}{c}
\mathrm{8,36,64,92,} \\ 
\mathrm{120,176,20,76,} \\ 
104,132,160,188 \end{array}
\right]}_{196}
\end{align*}
\end{definition}
\begin{align}
&\sm{\textrm{Bařina 7$n\pm$1 breadth-first line graph }\\\sm{\textrm{breadth-first }k=\croot1\tok2\tok3\tok\;\;\; \\\textrm{successor generation pairs }s_i\textcolor{blue}{\to}S_j}}
&& \rotatebox[origin=c]{270}{\(\textcolor{red}{\circlearrowright}\)}\,
\sm{\textcolor{red}{0}^0\\\sm{8\\\sm{L\\s_0}}}
\textcolor{blue}{\to}\sm{\textcolor{blue}{1}^1\\\sm{64\\\sm{LU\\S_1}}}
\tok
\sm{1\textcolor{red}{0}^1\\\sm{36\\\sm{UL^1\\s_1}}}
\textcolor{blue}{\to}
\sm{1\textcolor{blue}{1}^1\\\sm{512\\\sm{LU^2\\S_2}}}
\tok\
\sm{1\textcolor{red}{0}^2\\\sm{20\\\sm{UL^2\\s_2}}}
\textcolor{blue}{\to}
    \sm{10\textcolor{blue}{1}^1\\\sm{288\\\sm{ULU^1\\S_1}}}
\textcolor{red}{\to}
\sm{11\textcolor{red}{0}^1\\\sm{1168\\\sm{U^2L^1\\s_1}}}
\textcolor{blue}{\to}
\sm{\textcolor{blue}{1}^3\\\sm{4096\\\sm{LU^3\\S_3}}}
\tok
\sm{1\textcolor{red}{0}^3\\\sm{384\\\sm{UL^3\\s_3}}}
\textcolor{blue}{\to}
\tag{\ref{eq:bf7pm1}}
\end{align}

\subsubsection{Cayley trees $3n-1$ and Kontorovich \& Lagarias $5n+1$ conjectures: evidently refuted}
\begin{figure}[h]
\caption{The $3n-1$ function, for $k=5$ breadth-first levels rooted in $c=3\cdot1-1=2$}
\vspace{-0.4cm}
\noindent \center{\includegraphics[width=0.75\textwidth]{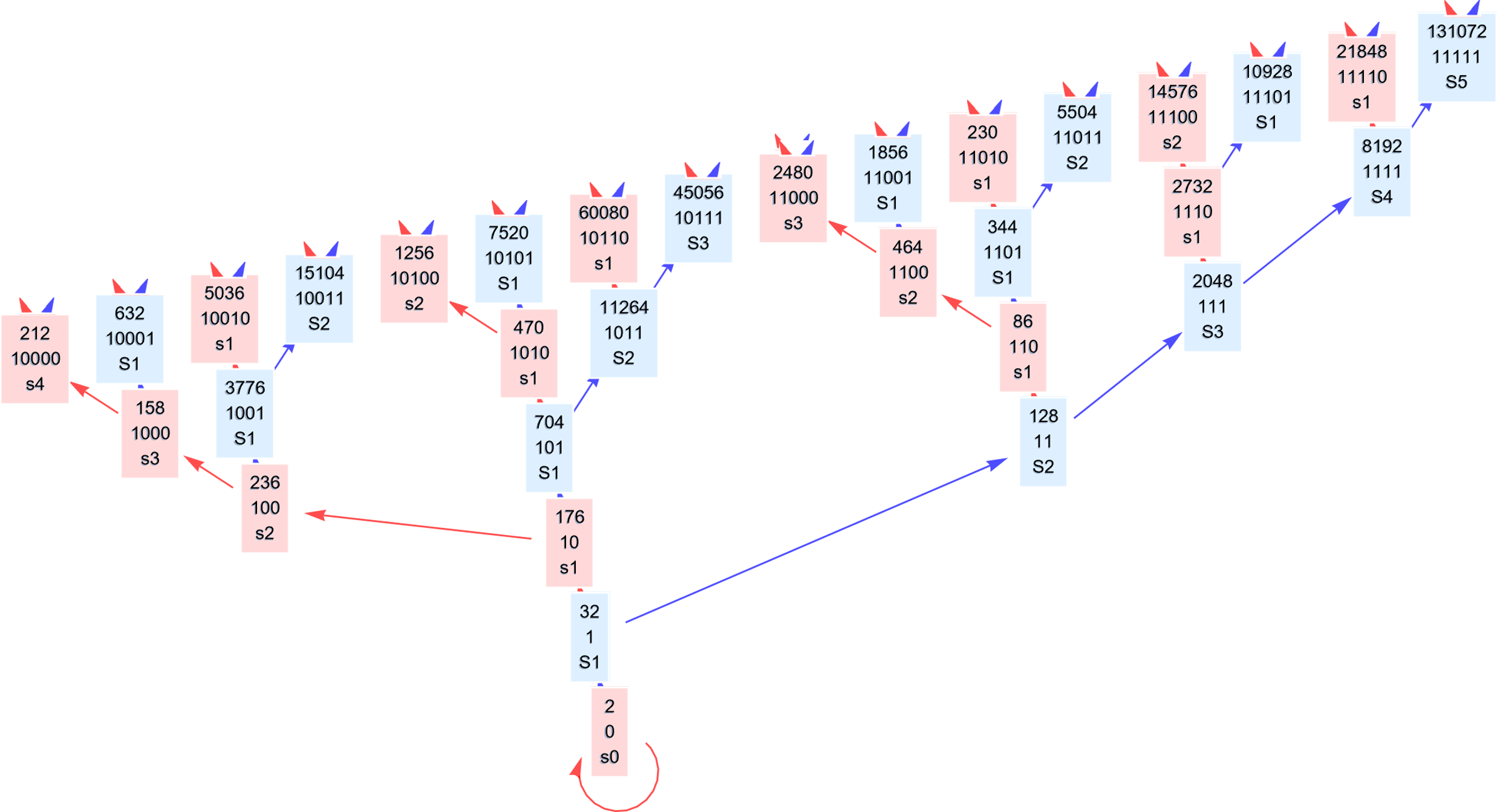}}
\label{fig:tree3m1}
\end{figure}

\begin{figure}[h!]
\caption{The $5n+1$ function, for $k=4$ breadth-first levels rooted either in $c=5\cdot1+1=6$ or $L(c)=16$}
\noindent \center{\includegraphics[width=0.75\textwidth]{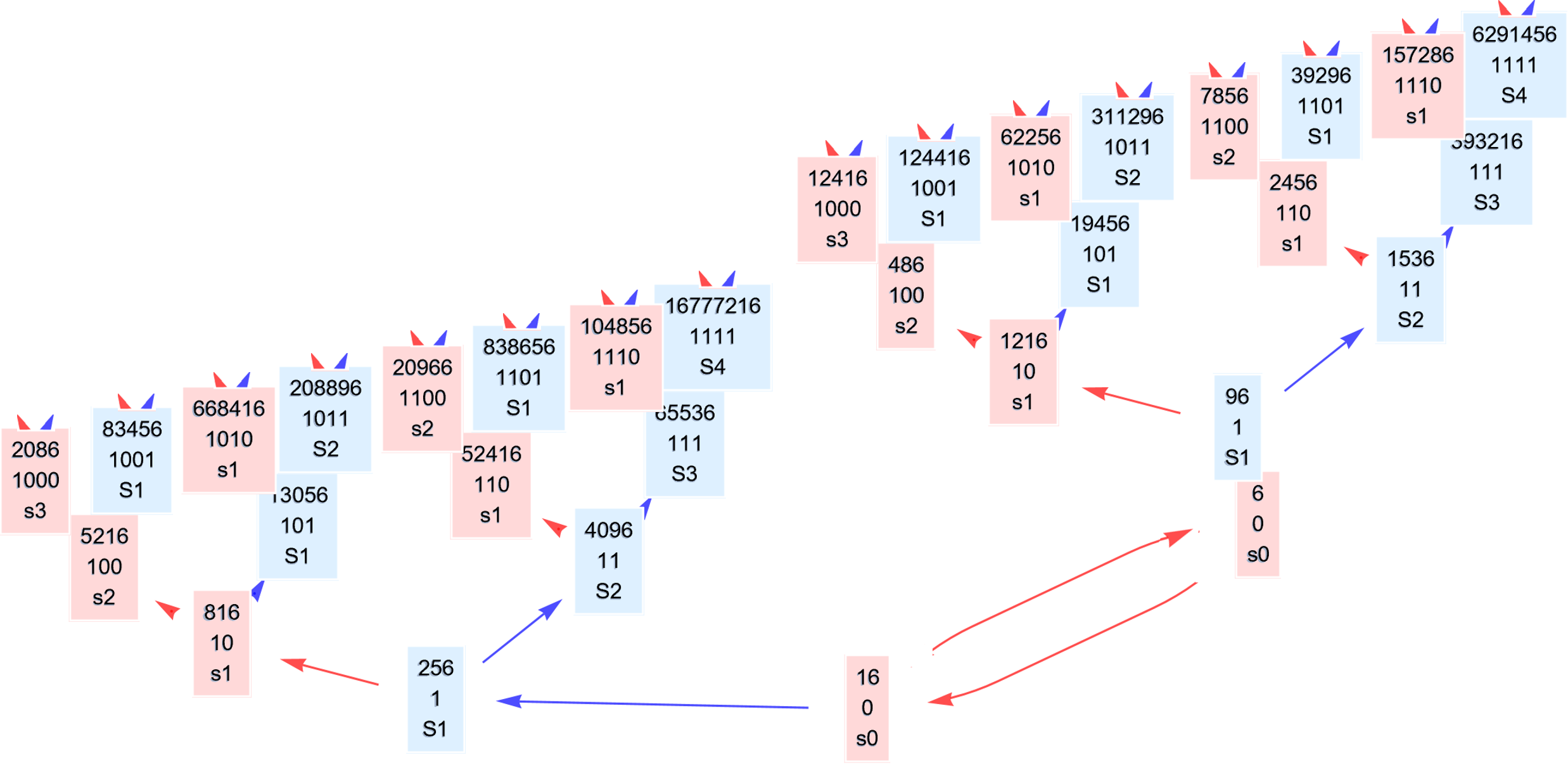}}
\label{fig:tree5p1}
\end{figure}

The $3n-1$ Cayley tree (Fig.\ref{fig:tree3m1}, Eq.\ref{eq:bf3m1}) with trivial root number pair $\sm{0\\3\cdot1-1=2}$ has a lowest trunk number $U(c)=32$ higher than its periodicity $(U(c=2)=2^{1+4}=32)>(2^va^2=2^13^2=18)$. The  $5n+1$ Cayley graph (Fig.\ref{fig:tree5p1}, Eq.\ref{eq:bf5p1}), is a shrub instead of a one-trunk-tree, with a trivial root $c=5\cdot1+1=6$ that is not a power-of-2, and a lowest trunk number higher than its periodicity $U(c=6)=6\cdot2^{4}=96 > 2^va^2=2^15^2=50$. 
\begin{align}
\label{eq:bf3m1}
&\sm{y,\;y\textrm{-axis, $3n-1$ breadth-first line graph }\;\\ \textcolor{white}{ b}\;\textrm{ with cyclic leading zero root }
\quad\quad\;\;\;\;}
&& \rotatebox[origin=c]{270}{\(\textcolor{red}{\circlearrowright}\)}\,
\sm{\textcolor{red}{0}\\\sm{2\\s_0}}
\textcolor{blue}{\to}\sm{\textcolor{blue}{1}\\\sm{32\\S_1}}
\tok
\sm{1\textcolor{red}{0}\\\sm{176\\s_1}}
\textcolor{blue}{\to}
\sm{1\textcolor{blue}{1}\\\sm{128\\S_1}}
\;\tok\;
\sm{10\textcolor{red}{0}\\\sm{236\\s_2}}
\textcolor{blue}{\to}
\sm{10\textcolor{blue}{1}\\\sm{704\\S_1}}
\textcolor{red}{\to}
\sm{11\textcolor{red}{0}\\\sm{86\\s_2}}
\textcolor{blue}{\to}
\sm{11\textcolor{blue}{1}\\\sm{2048\\S_3}}
\;\tok\;
\sm{100\textcolor{red}{0}\\\sm{158\\s_3}}
\textcolor{blue}{\to\cdots}
\end{align}

\vspace{-1cm}
\begin{align}
&\sm{y,\;y\textrm{-axis, $5n+1$ breadth-first line graph }\;\\ \textcolor{white}{ b}\;\textrm{ with cyclic leading zero root }
\quad\quad\;\;\;\;}
&& \rotatebox[origin=c]{270}{\(\textcolor{red}{\circlearrowright}\)}\,
\sm{\textcolor{red}{0}\\\sm{6\\s_0}}
\textcolor{blue}{\to}\sm{\textcolor{blue}{1}\\\sm{96\\S_1}}
\tok
\sm{1\textcolor{red}{0}\\\sm{96\\s_1}}
\textcolor{blue}{\to}
\sm{1\textcolor{blue}{1}\\\sm{128\\S_1}}
\;\tok\;
\sm{10\textcolor{red}{0}\\\sm{236\\s_2}}
\textcolor{blue}{\to}
\sm{10\textcolor{blue}{1}\\\sm{704\\S_1}}
\textcolor{red}{\to}
\sm{11\textcolor{red}{0}\\\sm{86\\s_2}}
\textcolor{blue}{\to}
\sm{11\textcolor{blue}{1}\\\sm{2048\\S_3}}
\;\tok\;
\sm{100\textcolor{red}{0}\\\sm{486\\s_3}}
\textcolor{blue}{\to\cdots}
\label{eq:bf5p1}
\end{align} 
\vspace{-0.5cm}

\subsection {Eq. series 2: depth-first / breadth-first traversal, serial paths / parallel arrows}
\label{sec:32eqs}
\begin{align} 
&\hspace{-0.25cm}\textrm{Each }T_{\ge0}\textrm { pair }\sm{y\\n} \textrm{ shows: serial \textbf{depth-first} paths to} \textrm{ successor \textit{nodes}}&\nonumber\\
&\quad\textrm{a red \textit{L}(eftward)-arrow to its red \textit{L}-child in red subtree }t_{0;\ge1}&\;\;\sm{\textcolor{red}{0}\\L}: \sm{y\\n}\in T_{\ge0}\tcr{\to}\sm{y\tcr{0}\\L(n)}\in t_{0;\ge1} \label{eq:redChild}
\\
&\quad\quad\textrm{and $L$'s blue serial $U$-path to generation } {S}_{j=1,2,3,\cdots}\textrm{ successors}
&\;\;
\sm{\textcolor{red}{0}\;\textcolor{blue}{1}^j\\LU^j}:
\sm{y\\n}\textcolor{blue}{\to}
\sm{y\textcolor{red}{0}\textcolor{blue}{1}^j_{\textcolor{blue}{S}j}\\LU^j(n)}
\label{eq:dfSj}
\\
&\quad\quad\textrm{(or, by Eq.\ref{eq:ucoar}, blue \textit{U}-path + conjugate red \textit{L}-arrow})&
\sm{\textcolor{blue}{1}^j\textcolor{red}{0}_j\\U^jL_j}:
\sm{y\\n}\textcolor{blue}{\to}
\sm{y\textcolor{red}{0}\textcolor{blue}{1}^j_{\textcolor{blue}{S}j}\\LU^j(n)}
\label{eq:dfSjconjug}
\\
&\quad\quad\quad\Rightarrow\textrm{Theorem \ref{the:cartesian}: $L$'s Blue $U$-successors in $QI$ Cartesian $C_H$} \notag
\\
&\quad\quad\quad\textrm{of breadth-first $\sm{y\\n}$ pairs (Eqs.\ref{eq:sorS},\ref{eq:bfpar}) on y-axis $(0,y)\subset C_H$}  
\label{eq:QI}
\\
&\quad\quad\textrm{random access }\sm{y\\n}\textrm{ to $L$'s blue \textit{U}- successors } 
&\;\;\sm{y\\n}\textcolor{blue}{\to}
\sm{\textcolor{red}{0}\textcolor{blue}{1}^1_{\textcolor{blue}{S}_{\,1}}\\LU^1(n)}; \sm{y\\n}\textcolor{red}{\to}
\sm{\textcolor{red}{0}\textcolor{blue}{1}^2_{\textcolor{blue}{S}_{\,2}}\\LU^2(n)};\textcolor{blue}{\cdots}
\label{eq:dfSjpar}
\\
&\quad\quad\textrm{random access }\sm{\textcolor{red}{0}\\c}\textrm{ to upward trunk cotree roots } \textcolor{blue}{T}_{j=1,2,3,\cdots} 
&\;\;\sm{\textcolor{red}{0}\;\textcolor{blue}{1}^j\\LU^j}:
\sm{\textcolor{red}{0}\\c}
\textcolor{blue}{\to}
\sm{\textcolor{red}{0}\textcolor{blue}{1}^j_{\textcolor{blue}{S}j}\\LU^j(c)}=
\sm{\textcolor{blue}{1}^j_{\textcolor{blue}{S}j}\\U^j(c)}
\label{eq:dfSjcon}
\\
&\quad\quad\quad\textrm{inverse subpath } {S}_{j=\cdots,3,2,1} \textrm{ successors to } \sm{y\\n} \textrm{ and } \cdots\;\cabroot 
&\;\;
\sm{\textcolor{blue}{1^{-j}}\;\textcolor{red}{0}^{-1}\\U^{-j}L^{-1}}:
\sm{y\textcolor{red}{0}\textcolor{blue}{1}^j_{\textcolor{blue}{S}j}\\LU^j(n)}
\textcolor{blue}{\to}\sm{y\\n}
\label{eq:IdfSj}
\\
&\quad\quad\textit{d}\textrm{ensity (red \textit{L}eftward subtree } t_{0;\ge1}) = \textit{d }\textrm {(red \textit{L}eftward cotrees)} &\;\;d(t_{0;\ge1})=\sum_{i=1}^{\infty}d(t_{i}) \label{eq:redsubco}
\\
&\quad\textrm{and a blue \textit{U}(pward)-arrow to its blue \textit{U}-child in blue subtree }T_{\ge1}&\;\;\sm{\textcolor{blue}{1}\\U}: \sm{y\\n}\in T_{\ge0}\tcb{\to}\sm{y\tcb{1}\\U(n)}\in T_{\ge1} \label{eq:blueChild}\\
&\quad\quad\textrm{and $U$'s red serial $L$-path to generation } {s}_{i=1,2,3,\cdots}\textrm{ successors}
&\;\;\sm{\textcolor{blue}{1}\;\textcolor{red}{0}^i\\UL^i}:\sm{y\\n}\textcolor{red}{\to}
\sm{y\textcolor{blue}{1}\textcolor{red}{0}^i_{\textcolor{red}{s}_{\,i}}\\UL^i(n)}
\label{eq:dfsi}
\\
&\quad\quad\textrm{(or, by Eq.\ref{eq:lcoar}, red $L$-path + conjugate blue $U$-arrow)} &
\sm{\textcolor{red}{0}^i\textcolor{blue}{1}_i\\L^iU_i}:
\sm{y\\n}\textcolor{blue}{\to}
\sm{y\textcolor{blue}{1}\textcolor{red}{0}^i_{\textcolor{red}{s}_{\,i}}\\UL^i(n)}
\label{eq:dfsiconjug}
\\
&\quad\quad\quad\Rightarrow\textrm{Theorem \ref{the:cartesian}: $U$'s red $L$-successors in $QII$ Cartesian $C_H$} \notag
\\
&\quad\quad\quad\textrm{of breadth-first $\sm{y\\n}$ pairs (Eqs.\ref{eq:sorS},\ref{eq:bfpar}) on y-axis $(0,y)\subset C_H$} 
\label{eq:QII}
\\
&\quad\quad\textrm{random access }\sm{y\\n}\textrm{ to $U$'s red \textit{L}-successors }
&\;\;\sm{y\\n}\textcolor{red}{\to}
\sm{\textcolor{blue}{1}\textcolor{red}{0}^1_{\textcolor{red}{s}_{\,1}}\\UL^1(n)}; \sm{y\\n}\textcolor{red}{\to}
\sm{\textcolor{blue}{1}\textcolor{red}{0}^2_{\textcolor{red}{s}_{\,2}}\\UL^2(n)};\textcolor{red}{\cdots}
\label{eq:dfsipar}
\\
&\quad\quad\textrm{random access }\sm{\textcolor{red}{0}\\c}\textrm{ to greedy branch cotree roots } \textcolor{red}{t}_{i=1,2,3,\cdots} 
&\;\;\sm{\textcolor{blue}{1}\;\textcolor{red}{0}^i\\UL^i}:\sm{\textcolor{red}{0}\\c}\textcolor{red}{\to}
\sm{\textcolor{blue}{1}\textcolor{red}{0}^i_{\textcolor{red}{s}_{\,i}}\\UL^1(c)}
\label{eq:dfsicon}
\\
&\quad\quad\quad\textrm{inverse subpath } {s}_{i=\cdots,3,2,1} \textrm{ successors to } \sm{y\\n} \textrm{ and } \cdots\;\cabroot 
&\;\;
\sm{\textcolor{red}{0^{-i}}\;\textcolor{blue}{1}^{-1}\\L^{-i}U^{-1}}:
\sm{y\textcolor{red}{0}\textcolor{blue}{1}^j_{\textcolor{blue}{S}j}\\LU^j(n)}
\textcolor{blue}{\to}\sm{y\\n}
\label{eq:Idfsi}
\\
&\quad\quad\textit{d}\textrm{ensity (blue \textit{U}pward subtree } T_{\ge1}) = \textit{d }\textrm {(blue \textit{U}pward cotrees)} &\;\;d(T_{\ge1})=\sum_{j=1}^{\infty}d(T_{j}) \label{eq:bluesubco}
\\
&\textrm
{and shows parallel arrows to \textit{sites} above nodes, bordered on the}&\;\;\nonumber\\ 
&\quad\quad\textrm{left by random accessible sites of $L$'s } \textcolor{blue}{s}_{j=1,2,3,\cdots} \textrm{ \textit{U}-successors}&\;\;
\sm{y\\n}\textcolor{blue}{\to}
\sm{y\textcolor{red}{0}\textcolor{blue}{1}^{1}_{\textcolor{blue}{s}_{j=1}}\\LU^{1}(n)},\;
\sm{y\\n}\textcolor{blue}{\to}
\sm{y\textcolor{red}{0}\textcolor{blue}{1}^{j=2}_{\textcolor{blue}{s}_{j=2}}\\LU^{2}(n)},\;\cdots
\label{eq:sitesleft}\\ 
&\quad\quad\textrm{right by random accessible sites of $U$'s }\textcolor{red}{s}_{i=1,2,3,\cdots} \textrm{ \textit{L}-successors}&\;\;
\sm{y\\n}\textcolor{red}{\to}
\sm{y\textcolor{blue}{1}\textcolor{red}{0}^{1}_{\textcolor{red}{s}_{i=1}}\\UL^{1}(n)},\;
\sm{y\\n}\textcolor{red}{\to}
\sm{y\textcolor{blue}{1}\textcolor{red}{0}^{2}_{\textcolor{red}{s}_{i=2}}\\UL^{2}(n)},\;\cdots
\label{eq:sitesright}
\end{align}
\newpage
\begin{flalign}
\noindent\textrm{and their root }\cabroot\textrm{ shows \textit{well-ordered} \textbf{breadth-first} root paths to \textit{wondrously ordered} a$n$+b numbers:} \quad\quad\quad\quad\quad\quad\quad\quad\quad\quad\quad\quad && \nonumber
\end{flalign}
\vspace{-0.7cm}
\begin{align}
&\sm{\textrm{a$n$+b breadth-first tape}\quad\quad\quad\quad\\\sm{\textrm{level }k=\croot1\tok2\tok3\tok\cdots\quad\; \\\quad\textrm{a-symmetric generations }\tcr{\to}s_i\textcolor{blue}{\to}S_j}}
&& \rotatebox[origin=c]{270}{\(\textcolor{red}{\circlearrowright}\)}\,
\sm{\textcolor{red}{0}^0\\\sm{L\\s_0}}
\textcolor{blue}{\to}\sm{\textcolor{blue}{1}^1\\\sm{LU\\S_1}}
\;\tok\;
\sm{1\textcolor{red}{0}^1\\\sm{UL^1\\s_1}}
\textcolor{blue}{\to}
\sm{\textcolor{blue}{1}^2\\\sm{LU^2\\S_2}}
\;\tok\;
\sm{1\textcolor{red}{0}^2\\\sm{UL^2\\s_2}}
\textcolor{blue}{\to}
\sm{10\textcolor{blue}{1}^1\\\sm{ULU^1\\S_1}}
\textcolor{red}{\to}
\sm{11\textcolor{red}{0}^1\\\sm{U^2L^1\\s_1}}
\textcolor{blue}{\to}
\sm{\textcolor{blue}{1}^3\\\sm{LU^3\\S_3}}
\;\tok\;
\sm{1\textcolor{red}{0}^3\\\sm{UL^3\\s_3}}
\textcolor{blue}{\to\cdots} \notag
% \label{eq:sorSn}
\\&\sm{k\;\textrm{ digits binary root path number }\\\sm{n_b: \;\#\textrm{ (root paths to) nodes level }k }\\\sm{\;\;\;\textrm{ yielding } 1\,+\,\sum_{k=1}^{2,3,\dots,K}\,2^{k-1}\,=\;1+ \\\quad\; (2^K-1)=2^K\textrm{ (root paths to) nodes}}} 
&&\sm{\textcolor{red}{1}\\\sm{1\\1+2^{k-1=0}}},
\sm{\textcolor{blue}{1}\\\sm{1\\=2^{\textcolor{white}{1}}   }};
\sm{\textcolor{red}{2}\\\sm{2\\2^{2-1=1}}},
\sm{\textcolor{blue}{2}\\\sm{2\\=2^{\textcolor{white}{1}} }};
\quad\quad\;\;\sm{\textcolor{red}{3}\\\sm{3\\{\textcolor{white}{1^1}}   }},
% \sm{\textcolor{blue}{3}\\3},
\sm{\textcolor{blue}{3}\\\sm{\textcolor{blue}{3}\\2^{3-1=2}}}
\sm{\textcolor{red}{3}\\\sm{\textcolor{blue}{3}\\=4^{\textcolor{white}{1}} }},
\sm{\textcolor{blue}{3}\\ \sm{3\\{\textcolor{white}{1^1}} }};
\quad\quad\quad\sm{\textcolor{red}{4}\\\sm{4\\2^{k-1}}}, \cdots
\label{eq:sorS} 
\\&\sm{\textrm{Conway 1$n$+1 breadth-first tape}\quad\quad
\\\sm{\textrm{level }k=\croot1\tok2\tok3\tok\dots\textcolor{white}{\to}\quad\quad\\\textrm{a-symmetric generations  }\tcr{\to}s_i\textcolor{blue}{\to}S_j}}
&& \rotatebox[origin=c]{270}{\(\textcolor{red}{\circlearrowright}\)}\,
\sm{\textcolor{red}{0}^0\\\sm{2\\\sm{L\\s_0}}}
\textcolor{blue}{\to}\sm{\textcolor{blue}{1}^1\\\sm{4\\\sm{LU\\S_1}}}
\tok
\sm{1\textcolor{red}{0}^1\\
\sm{6\\\sm{UL^1\\s_1}}}
\textcolor{blue}{\to}
\sm{\textcolor{blue}{1}^2\\\sm{8\\\sm{LU^2\\S_2}}}
\tok\
\sm{1\textcolor{red}{0}^2\\
\sm{10\\\sm{UL^2\\s_2}}}
\textcolor{blue}{\to}
    \sm{10\textcolor{blue}{1}^1\\\sm{12\\\sm{ULU^1\\S_1}}}
\textcolor{red}{\to}
\sm{11\textcolor{red}{0}^1\\\sm{14\\\sm{U^2L^1\\s_1}}}
\textcolor{blue}{\to}
\sm{\textcolor{blue}{1}^3\\\sm{16\\\sm{LU^3\\S_3}}}
\tok
\sm{1\textcolor{red}{0}^3\\
\sm{18\\\sm{UL^3\\s_3}}}
\textcolor{blue}{\to}
\label{eq:bf1p1}
\\&\sm{\textrm{Collatz 3$n$+1 breadth-first tape }\quad\quad\\
\sm{\textrm{level }k=\croot1\tok2\tok3\tok\cdots\quad\quad 
\\\quad\textrm{a-symmetric generations  }\textcolor{red}{\to} s_i\textcolor{blue}{\to}S_j\;}}
&& \rotatebox[origin=c]{270}{\(\textcolor{red}{\circlearrowright}\)}
\sm{\textcolor{red}{0}^0\\\sm{4\\\sm{L\\s_0}}}
\textcolor{blue}{\to}\sm{\textcolor{blue}{1}^1\\\sm{16\\\sm{LU\\S_1}}}
\tok
\sm{1\textcolor{red}{0}^1\\
\sm{40\\\sm{UL^1\\s_1}}}
\textcolor{blue}{\to}
\sm{\textcolor{blue}{1}^2\\\sm{256\\\sm{LU^2\\S_2}}}
\tok\
\sm{1\textcolor{red}{0}^2\\\sm{52\\\sm{UL^2\\s_2}}}
\textcolor{blue}{\to}
    \sm{10\textcolor{blue}{1}^1\\\sm{160\\\sm{ULU^1\\S_1}}}
\textcolor{red}{\to}
\sm{11\textcolor{red}{0}^1\\\sm{340\\\sm{U^2L^1\\s_1}}}
\textcolor{blue}{\to}
\sm{\textcolor{blue}{1}^3\\\sm{1024\\\sm{LU^3\\S_3}}}
\tok
\sm{1\textcolor{red}{0}^3\\
\sm{34\\\sm{UL^3\\s_3}}}
\textcolor{blue}{\to}
\label{eq:cywell}
\\&\sm{\textrm{Bařina 7$n\pm$1 breadth-first tape }\;\;\;\quad\\\sm{\textrm{level }k=\croot1\tok2\tok3\tok\dots\quad\;\; \\\quad\textrm{a-symmetric generations }s_i\textcolor{blue}{\to}S_j\;\;}}
&& \rotatebox[origin=c]{270}{\(\textcolor{red}{\circlearrowright}\)}\,
\sm{\textcolor{red}{0}^0\\\sm{8\\\sm{L\\s_0}}}
\textcolor{blue}{\to}\sm{\textcolor{blue}{1}^1\\\sm{64\\\sm{LU\\S_1}}}
\tok
\sm{1\textcolor{red}{0}^1\\
\sm{36\\\sm{UL^1\\s_1}}}
\textcolor{blue}{\to}
\sm{\textcolor{blue}{1}^2\\\sm{512\\\sm{LU^2\\S_2}}}
\tok\
\sm{1\textcolor{red}{0}^2\\
\sm{20\\\sm{UL^2\\s_2}}}
\textcolor{blue}{\to}
    \sm{10\textcolor{blue}{1}^1\\\sm{288\\\sm{ULU^1\\S_1}}}
\textcolor{red}{\to}
\sm{11\textcolor{red}{0}^1\\\sm{1168\\\sm{U^2L^1\\s_1}}}
\textcolor{blue}{\to}
\sm{\textcolor{blue}{1}^3\\\sm{4096\\\sm{LU^3\\S_3}}}
\tok
\sm{1\textcolor{red}{0}^3\\
\sm{384\\\sm{UL^3\\s_3}}}
\textcolor{blue}{\to}
\label{eq:bf7pm1}
\\&\sm{\textrm{random accessible }\\ \textrm{root successors }\;\;\;\;\,}
&& \rotatebox[origin=c]{270}{\(\textcolor{red}{\circlearrowright}\)}\,
\sm{\sm{\textcolor{red}{0}^0\\{{s}_0}}\\c}
\textcolor{blue}{\to}
    \sm{\sm{\textcolor{blue}{1}^1\\S_1}\\U(c)}
,\;
\rotatebox[origin=c]{270}{\(\textcolor{red}{\circlearrowright}\)}\,
\sm{\sm{\textcolor{red}{0}^0\\{{s}_0}}\\c}
\textcolor{red}{\to}
\sm{\sm{1\textcolor{red}{0}^1\\{{s}_1}\\UL(c)}}
,\;
\rotatebox[origin=c]{270}{\(\textcolor{red}{\circlearrowright}\)}\,
\sm{\sm{\textcolor{red}{0}^0\\{{s}_0}}\\c}
\textcolor{blue}{\to}
\sm{\sm{\textcolor{blue}{1}\\{{S}_2}}\\U^2(c)}
,\;
\rotatebox[origin=c]{270}{\(\textcolor{red}{\circlearrowright}\)}\,
\sm{\sm{\textcolor{red}{0}^0\\{{s}_0}}\\c}
\textcolor{red}{\to}
\sm{\sm{1\textcolor{red}{0}^2\\{{s}_2}\\UL^2(c)}}
\textcolor{blue},\;\;\cdots
&\label{eq:bfpar}
\\&\textrm{with as $y$-axis projections:}&&\textrm{ 
\textit{wondrously ordered} root paths to  
\textit{well-ordered} a$n$+b numbers}
\\&\sm{\textrm{wondrously ordered root paths }y \textrm{ to }\\\sm{ \tcr{\to}\textcolor{red}{s}-\textrm{ and }\tcb{\to}\textcolor{blue}{S}-\textrm{generations}\\ \textrm{well-ordered a\textit{n}+b \textit{y}-axis, 3\textit{n}+1:}}}
&& 
\sm{\sm{0\\{\textcolor{red}{s}_0}\\4}}
\to\sm{\sm{1\\{\textcolor{blue}{S}_1}\\16}}
\to\sm{\sm{10^4\\{\textcolor{red}{s}_4}\\22}}
\to\sm{\sm{10^3\\{\textcolor{red}{s}_3}\\34}}
\to\sm{\sm{10\\{\textcolor{red}{s}_1}\\40}}
\to\sm{\sm{10^2\\{\textcolor{red}{s}_2}\\52}}
\to\sm{\sm{10^410\\{\textcolor{red}{s}_1}\\58}}
\to\sm{\sm{1010^2\\{\textcolor{red}{s}_2}\\70}}
\to\cdots
\label{eq:xpar}
\\&\Rightarrow \textrm{Lemma \ref{lem:2comb} (\textit{S}-challenge)}:&&\textrm{left-right asymmetric $\tcr{\to} s_i\textcolor{blue}{\to S_j}$ combinators } \tcr{\to}UL^{i=1,2,3,\dots}\tcb{\to}LU^{j=1,2,3,\dots}\notag
\\&\quad\;\Rightarrow\textrm{Lemmas \ref{lem:cantorno}-\ref{lem:wolframno}:}&&\textrm{no a\textit{n}+b-conjecture undecidability by Cantor diagonals,}\notag
\\& &&\textrm{Gödel numbers, Conway unpredictability, Wolfram irreducibility}\notag
\\&\Rightarrow\textrm{Theorem \ref{the:fractal}:}
&&
\textrm{$y$-axis $T_{\ge0}$ (Fig.\ref{fig:tge0}b): well-ordered a$n$+b, wondrously ordered root paths $y$} \notag
\\&\quad\quad\quad&&
\textrm{$y$-axis $H_C$ (Figs.\ref{fig:hamiltonCo},\ref{fig:hamiltonBa}): well-ordered root paths $y$, wondrously ordered a$n$+b}\notag
\end{align}

\noindent Leibniz' study of successors in dynastic and biblical successor generations  \citep[Nachfolger, Generation]{leibnizKatalog} anticipated the computer science distinction between the \textit{depth-first, children-first} traversal, of a rooted graph (Eqs.\ref{eq:redChild}-\ref{eq:sitesright}) and its \textit{breadth-first, brothers/sisters-first} traversal (Eqs.\ref{eq:sorS}-\ref{eq:xpar})  \cite{whlisp,RN23,RN32}.

The red cyclic root $\cabroot$ (Eq.\ref{eq:Croot}) implies  \textit{biperiodic} depth-first traversals by red $L$(eftward) and blue $U$(pward) depth-first arrow function digits (Eqs.\ref{eq:rootPairs}, \ref{eq:redChild}, \ref{eq:blueChild}) generating unique, non-commutative binary numbers of root paths that are non-confluent (Eq.\ref{eq:tleg}). 
Biperiodic functions are renowned due to Wiles's proof of Fermat's Last Theorem, and the Clay Millennium Prize for proving the Birch and Swinnerton-Dyer Conjecture, based on the independent growth of elliptic curves (and of ellipses known from high school) along the Cartesian $x$- and $y$-axes.

Applied to \textit{even positive} Collatz-Conway numbers $n$=a$n_{odd}$+b, the depth-first red and blue function steps, $L:n\textcolor{red}{\to}L(n)=fg^p(n)=(n-b)/a\cdot2^p$ and $U:n\textcolor{blue}{\to}U(n)=g^q(n)=n\cdot2^q$ return respectively $(n-b)/a\cdot2^p$ (Eq.\ref{eq:redChild}) and $n\cdot2^q$ (Eq.\ref{eq:blueChild}). Every binary root path number $y$ has a leading $1$, except for the even, leftward, self-regenerating, root path number $\croot y_{=\textcolor{red}{0}}$ in generation $\textcolor{red}{s}_0$ of numbers without an upward predecessor, including the cyclic root number pair $\rotatebox[origin=c]{270}{\(\textcolor{red}{\circlearrowright}\)}\,\sm{\textcolor{red}{0}\\c=a\cdot1+b}$  (Eq.\ref{eq:Croot}). 
A root path $y$ concatenates red and blue binary arrow digits $\textcolor{red}{0}$ and $\textcolor{blue}{1}$ representing the depth-first, children-first binary function digits $L$ and $U$ (Eqs.\ref{eq:redChild},\ref{eq:blueChild}). Each number pair $\sm{y\\n}$ is also depth-first paired to either its red successor generation $\textcolor{red}{s}_i$ if $y$ ends in $1\textcolor{red}{0}^i$ generated by trailing $UL^i$-iterations (Eq.\ref{eq:dfsi}), or its blue successor generation $\textcolor{blue}{S}_j$ if $y$ ends in $0\textcolor{blue}{1}^j$ generated by trailing $LU^j$-iterations (Eq.\ref{eq:dfSj}). Generation $s_0$ comprises root path numbers $\croot\textcolor{red}{0}$ with a trailing $\textcolor{red}{0}$ without a leading $\textcolor{blue}{1}$. The \textit{d}ensity of red and blue successors of $T_{\ge0}$ number pairs in the red subtree $t_{0;\ge0}$ and the blue cotree $T_{\ge1}$ is equal to the density of all red successor generations $s_i$ in red cotrees $t_i$, respectively of all blue successor generations $S_j$ in blue cotrees $T_j$ (Eqs.\ref{eq:redsubco},\ref{eq:bluesubco}). 

The breadth-first, brothers-sisters-first traversal of a Cayley a\textit{n}+b graph encounters in each successive depth-first generation \textit{i} alternatingly red, L(eftward) to blue U(pward) number pairs $s_i\textcolor{blue}{\to}S_i$ (Eq.\ref{eq:sorS}; specified for Collatz 3$n$+1 function in Eq.\ref{eq:cywell}). The regenerating $\textcolor{red}{0}$ pairs with the lowest a$n$+b branching number $c=a\cdot1+b$ (Eq.\ref{eq:sorS}), which is for the Collatz $3n+1$ function $c=a\cdot1+b=4$  (Eq.\ref{eq:sorS}), regenerating by $L:c\,\textcolor{red}{\to}\, L(c)=fg^2(c)=(c-1)/3\cdot 2^2=1\cdot2^2=4=c$. Climbing red $\tok$ arrows show staircase steps to the next higher breadth-first-generation level. The breadth-first Eqs.\ref{eq:sitesleft}-\ref{eq:xpar}
cover the generalization to graph transformations by Hartmut Ehrig (1944-2016) \cite{ehrig1979}, inspired by the serial-paths-to-parallel-arrows conversion in electric circuits by Kirchhoff \citep[6A]{RN3245}. A finite tree has only $1$ site surrounding each of its nodes. In the infinite Cayley tree $T_{\ge0}$ each node number pair is also the site number pair of its upwardly narrowing but unclosed infinite site above it, bordered to the left by a path of blue serial arrows (Eq.\ref{eq:dfSj}) and to the right by a path of red serial arrows (Eq.\ref{eq:dfsi}). This holds also for the cyclic root $\cabroot$, which is the only node number pair with its own site also below it, and is therefore bordered to the left by the upward trunk of blue arrows (Eq.\ref{eq:dfSjcon}), and bordered to the right by the greedy branch of red arrows (Eq.\ref{eq:dfsicon}).  

A question in the post WW-II Turing/Von Neumann age of serial computing was whether “massively parallel” computations  \citep[p999, 1035]{RN23} could bypass the computational irreducibility of iterations of serial `rules' / functions. In the $21^{th}$ century, the serial-to-parallel question was answered by ``massively parallel" arrows to, ``random access" arrows to, (colour) pixels in aver larger rectangles within an infinite Cartesian Half graph $C_H$. The serial-paths-to-nodes to parallel-arrows-to-sites transformation enables the metamorphosis of the 3-regular, left-right asymmetric, Cayley tree $T_{\ge0}$ of serial paths of arrows, as the 4-regular integer Cartesian upper half graph $C_H$ of successive parallel, random access, left-right-symmetric, arrows to blue site number pairs in Quadrant I and red site number pairs in Quadrant II to the right, respectively to the left of the $y$-axix (Figs.\ref{fig:hamiltonCo},\ref{fig:hamiltonBa}). Graph $C_H$ enables either an a$n$+b conjecture proof or refutation by a Hamiltonian tour of all $a\nmid$ a$n_{\textup{odd}}$+b numbers from and back to the trivial root $\cabroot$ (Theorem \ref{the:cartesian} below). 

The depth-first and breadth-first equations in this subsection suffice to present the Cayley graph $T_{\ge0}$ as a specification of the infinite propositional logic graph proposed by Donald Knuth \cite{knuthAlltrees,knuthBDD} (Lemma \ref{lem:knuth}), and to solve the (one?) $S$-combinator challenge proposed by Stephen Wolfram \cite{combinators} (Lemma \ref{lem:2comb}).

\subsubsection{The infinite Cayley-Cantor tree specifies Knuth’s infinite propositional logic tree}

Finite propositional logic truth tables combinations of either $\sm{\text{True}\\\top}$ or $\sm{\text{False}\\\bot}$ propositions are presented by the binary number digits $\sm{\textcolor{blue}{1}\\\top,\text{True}}$ and $\sm{\textcolor{red}{0}\\\bot,\text{False}}$. However, `at least at a formal level', the numbers `0' or `1' for truth ($\sm{\textrm{True}\\\top}$) or falsehood ($\sm{\textrm{False}\\\bot}$) are arbitrary \citep[807]{RN32}. Donald Knuth casts
the infinite propositional logic tree  as the tree of all possible combinations of $\top$-propositions and $\bot$-propositions without a cyclic root \cite{knuthAlltrees,knuthBDD}.  

\begin{lemma}
\label{lem:knuth}
The Cayley tree $T_{\ge0}$ of a Collatz-Conway a\textit{n}+b function paired with the Cantor tree of binary root path numbers (Theorem \ref{the:caytree}) transforms into Donald Knuth's infinite propositional logic tree \cite{knuthAlltrees,knuthBDD} with cyclic root $\croot\sm{\tcr{0}\\\top,\text{\footnotesize{\textup{true}}}}$ and non-arbitrary truth values $\sm{\tcr{0}\\\top\text{\footnotesize{\textup{,True}}}}$  and $\sm{\tcb{1}\\\bot\text{\footnotesize{\textup{,False}}}}$.

\noindent\textit{Proof: }
\textup{The Cayley tree $T_{\ge0}$ of a Collatz-Conway a\textit{n}+b function paired with the Cantor tree of binary root path numbers (Theorem \ref{the:caytree}) is the tree of all possible binary root path numbers of combinations of the binary digits $\sm{\tcr{0}\\\text{\footnotesize{s-\textup{generation, }}}L\textup{eftward }, }$  and $\sm{\tcb{1}\\\text{\footnotesize{S-\textup{generation, }}}\text{\footnotesize{U\textup{pward}}}}$ with cyclic root $\sm{\tcr{0}\\c}$. The cyclic root is the tautological proposition $\sm{\tcr{0}\\c:\;\top,\text{\footnotesize{\textup{ 
True by definition}}}}$ that a definition, an axiom, a constant, or any proposition $c$ of which the truth is non-falsifiable and non-deniable, is $\sm{\tcr{0}\\c:\top}$. A blue cyclic, tautological, proposition is non-existent as the denial, refutation, or falsification of a False proposition is a non-cyclic, non-tautological, $\sm{\tcr{0}\\\top\textup{rue}}$ proposition.\qed}
\end{lemma}

\subsubsection{Solution $S$-combinator challenge}

Stephen Wolfram launched the (one?) $S$-combinator challenge whether mathematical functions can be rewritten as iterators of one $S$-combinator \cite{combinators}, as pictured in \textit{A new kind of science} \cite{RN23}, given that combinations of non-trivial mathematical logic functions by Mozes Schönfinkel (1888-1942), implemented by Wolfram in \textit{Mathematica} \cite{RN32},  require two combinators.

\begin{lemma} 
\label{lem:2comb} 
\textbf{Solution \textit{S}-combinator challenge}. Generating the infinite propositional logic Cayley tree $T_{\ge0}$ (Lemma \ref{lem:knuth}) requires a depth-first $S$-generations combinator and a depth-first $s$-generations combinator. 

\noindent{Proof}:
\textup{Breadth-first iterations of a biperiodic a$n$+b function Eq.\ref{eq:sorS} show \textit{two alternating depth-first combinators}, $\tcr{\to}UL^{i=1,2,3,\cdots}$ (Eq.\ref{eq:dfsi}) and $\tcb{\to}LU^{j=1,2,3,\cdots}$ (Eq. \ref{eq:dfSj}), generating from the cyclic root $\cabroot$ (Eq.\ref{eq:Croot}), at each breadth-first level $k=\croot1\tok2\tok3\tok\cdots$, alternately a$n$+b nodes/propositions in \textit{left-right a-symmetric}  red leftward $\tcb{\to}S_j\tcr{\to}s_i$ and blue upward $\tcr{\to}s_i\tcb{\to}S_j$-generations $i$ and $j$ (Eqs.\ref{eq:sorS}-\ref{eq:bf7pm1}). The a-symmetry of the red $s$-combinator and the blue $S$-combinator is easily observed as each breadth-first proposition level $k=1,2,3,\dots$ shows a leftmost proposition $\sm{1\tcr{0}^{k-1}\\\top}$, and a rightmost proposition $\sm{\tcb{1}^k\\\bot}$ (Figs.\ref{fig:treeP1},\ref{fig:tge0},\ref{fig:tree7pm1},\ref{fig:tree3m1}). }
\end{lemma}

\subsubsection{No a\textit{n}+b conjecture is undecidable by Cantor diagonals}

Cantor proved from 1874 onward \cite{canzer} in various ways the no longer contested proposition that `the continuum' $L$ of real numbers $0\le z\le1$ in the open interval (0,1) cannot be counted by, and cannot be paired one-by-one to, the countably infinite set $M$ of integer even and odd binary numbers with a terminal digit  $\tcr{0}$ or $\tcr{1}$ \citep[77]{cantor1890}, whose count/cardinality is denoted as $\aleph_0$. Cantor proved that the $2^{\aleph_0}=\aleph_1$ numbers in subsets of the infinite number set cannot be counted by, and not be paired one-by-one to the  $\aleph_0$ numbers in the countably infinite set $\mathbb{N}_1=1,2,3,\dots$, implying an infinite hierarchy of $\aleph_{0,1,2,\dots}$ infinities. His famous 1890/1891 `diagonal proof' \cite{cantor1890} of $\aleph_{\ge1}>\aleph_0$ cardinalities assumes an infinite square number matrix with one finite corner and  successor diagonals of exceeding length $2,3,4,\dots$. 

\begin{lemma} 
\label{lem:cantorno}
No Collatz-Conway \textup{a}n+\textup{b} conjecture with Cayley tree $T_{\ge0}$ is undecidable by Cantor diagonals.

\noindent Proof: \textup{
The Cayley tree $T_{\ge0}\subseteq G_{\ge0}$ (Figs.\ref{fig:tge0},\ref{fig:treeP1},\ref{fig:tree7pm1},\ref{fig:tree3m1},\ref{fig:tree5p1}) of well-ordered, breadth-first-ordered non-negative even and odd integer numbers $y_{\ge0}\in([0,1]_{2}\stackrel{\textup{bin}}{=}[0,1]_{10})\,$ in the binary Cantor-tree \citep[S6]{cantor1883bin} of root paths, leading to Collatz-Conway $a\nmid n=$\textup{a} n$_{\textup{odd}}$+\textup{b} numbers, shows up also as an infinite square matrix with:
\begin{itemize}
\item the root-trunk pair of its cyclic root and its upward trunk arrowhead node $\cabroot\tcb\to\sm{\tcb{1}\\U(c)=c\cdot2^q},\;$ of $1+2^{k-1}=1+2^{1-1}=2$ root paths of $k=1$ binary digit up to tree level $k=1$, 
\item exclusively power-of-2 diagonal tree successor levels $k=2,3,4,\dots$, comprising $2^{k-1}$ alternately even, red, $L$eftward binary numbers and odd, blue, $U$pward binary numbers of root paths with terminal digits $\tcr{\to 0}\tcb{\to 1}$, each comprising $k$ binary digits (Eq.\ref{eq:sorS}), amount to a Cantor totality (`aktualer Gesamtheit' \citep[75]{cantor1890},\citep[p374]{canzer}) of $2^K$  binary root paths numbers led to $2^K$ alternately even/red and odd/blue terminal arrows $\tcr{\to0}\tcb{\to1}$ of $2^K$ root paths:
\begin{align}
&\quad\quad\quad\quad1+\sum_{k=1}^{K}2^{k-1}\quad=\quad1+(2^{K+1}-1)\quad=\quad1+(2^K-1)\quad=\quad2^K \label{eq:2K}\\
&\textup{in which}\nonumber\\ 
&k=1\tok2\tok3\tok\dots:\textup{$ $breadth-first tree level $k$ with binary root path numbers of $k$ binary digits (Eq.\ref{eq:sorS})} \nonumber\\
&\textup{$K $:\; Cantor totality of breadth-first tree levels $k$ for which $K+1=K$}\nonumber
\end{align}
\end{itemize}
}

\noindent\textup{Eq.\ref{eq:2K} falsifies the $\aleph$-hierarchy as Cayley tree $T_{\ge0}$ of non-negative integer numbers shows equal cardinalities $(2^K=2^{\aleph_0}=\aleph_1)$
$=(K=\aleph_0)$. Pairwise flipping successive terminal digits in the `diagonal' breadth-first levels of Cayley tree $T_{\ge0}$
merely reverses their left-to-right-reading order $\tcr{\to0}\tcb{\;\to1}$ (e.g. Sumerian, Roman, Western) into a right-to-left-reading order $\tcb{1\leftarrow}\;\tcr{0\leftarrow}$ (e.g. Hebrew, Arabic, Persian)
}\qed.
\end{lemma}

\subsubsection{No a\textit{n}+b conjecture is undecidable by Gödel numbering}

\noindent Gödel numbering, developed by Kurt Gödel (1906-1978), pairs logical propositions with numbers, which enables a proof of the undecidability, of neither a proof nor a refutation, of a proposition stating its own undecidability: "Wir haben also einen Satz vor uns, der seine eigene Unbeweisbarkeit behauptet" \citep[p.175]{godel}. Prototypical examples are the liar paradox "I am lying" \cite{godel,nagnew} and the unique paradox "you are not unique, as everybody is unique". Gödel numbers do not recognize the, often implicit, symbol : representing citations, attributions, embedded worlds, and (iterated) mathematical functions. A proposition asserting: its unprovability is transformed to, and recognized as, the non-contradictory, non-paradoxal, proposition "$x$ (always lying, always unique) is neither proven  
$\sm{\tcr{0}\\\top,\text{True}}$, nor refuted  $\sm{\tcb{1}\\\bot,\text{False}}$. 

\begin{lemma}
\label{lem:godelno} 
No Collatz-Conway \textup{a}n+\textup{b} conjecture is undecidable by Gödel numbering. Proof:
\begin{enumerate}
\item \textup{For the sequence of binary numbers $\mathbb{N}_0=0,1,\,10,11,\,100,\dots$ the, not-biperiodic, Gödel successor function $f$ of the number $0$ and of arbitrary `natural' numbers $a\in1,\,10,11,\,100,\dots$ gives in positional/diagrammatic/leftmost-innermost number notation the sequence $a,af^1,af^2,af^3,\dots$ \citep[177]{godel}, not distinguishing the biperiodic, left-right-asymmetric, sequences given by the two binary Cantor-functions  \citep[\S6]{cantor1883bin}(Eq.\ref{eq:sorS})\; $\tcr{0}:(=g:)\,a\tcr{\to}10a+0=a\tcr{0}$ and $\tcb{1}:(=f:)\,a\tcb{\to}10a+1\stackrel{\textup{bin}}{=}a\tcb{1}$, when applied to the lowest trunk number $a=\tcb{1}$ giving $1\tcr{0}^1,1\tcr{0}^2,1\tcr{0}^3,\dots$ and $\tcb{1}^2,\tcb{1}^3,\tcb{1}^4,\dots$. Since the not-biperiodic Gödel successor function $f$ gives either only even or only odd binary (root path) numbers paired with a$n$+b numbers, Gödel numbers do not represent (root paths to) all a$n$+b numbers in the Cayley graph $T_{\ge0}$ of a Collatz-Conway a$n$+b function.}
\item\textup{Reversely, not a single Cayley graph $T_{\ge0}$ of a Collatz-Conway a$n$+b function includes a proposition that is neither proven nor refuted, as the terminal digits of all its binary root path numbers denote by Lemma \ref{lem:knuth} either red Leftward proven $\sm{\tcr{0}\\\top,\text{True}}$ or blue Upward refuted $\sm{\tcb{1}\\\bot,\text{False}}$ propositions.}\qed
\end{enumerate}
\end{lemma}

\noindent The implicit citation/attribution/function symbol : in the liar paradox and the unique paradox is explicit in politics, religion, economics, literature, and in communication content analysis language models \cite{osgood,Cuil86graph,Cuil88AI,att2008semantic,att2021validity,attlaurer2024less}.

\subsubsection{No a\textit{n}+b conjecture is undecidable by Conway unpredictability}
In the 1971 article in which John H. Conway generalized the Collatz 3$n$+1 conjecture to a$n$+b conjectures \cite{RN3209}, he states their undecidability because of unpredictability by a \textit{Minsky program} whether its iterations reach for all numbers integer odd numbers, and finally 1, instead of `undefined' fractions \cite{fractran,endrullisGrabmayer}. Minsky iterations both enable (repeated) addition and subtraction, and (repeated) mulriplication and division \cite{RN3209}.
\begin{lemma}
\label{lem:conwayno} 
No \textup{a}n+\textup{b} conjecture is undecidable by Conway unpredictability. 

\noindent Proof:
\textup{An a$n$+b function is not undecidable by Conway unpredictability resting on the non-applicability of a Minsky program, since in a breadth-first traversal $n$=a$n_\textrm{odd}$+b numbers on the nodes and sites of the a$n$+b Cayley tree $T_{\ge0}$  (Eqs.\ref{eq:sorS}-\ref{eq:xpar}) are paired with consecutive non-negative binary root path numbers in the Cantor-tree  $y\in\mathbb{N}_0=0,1\tok10,11\tok100,101,110,111\tok1000,\dots$ (Eq.\ref{eq:sorS}), enabling a Minsky program allowing for (repeated) addition and subtraction, and (repeated) multiplication and division.}
\qed 
\end{lemma}

\subsubsection{No a\textit{n}+b conjecture is undecidable by Wolfram irreducibility}
Wolfram pictured iterations of thousands of `simple', one-periodic, functions (“rules”), generating “trivariate”, 2-colour (black-white) trees \cite{RN23,combinators}, starting from “initial conditions” \citep[index]{RN23}---`origins', `roots', `poles'  \cite{vhcombi,nieuw1976}. The apparent a-periodicity of black cells generated by successive serial iterations of some of these simple `rules' (or functions, e.g. Wolfram's rule 30) underlies Wolfram's hypotheses of the "computational irreducibility" of simple serial iteration rules \citep[ch.11]{RN23}  and their "computational equivalence" \citep[ch.12]{RN23} with arbitrarily complex serial iteration rules. "Whenever computational irreducibility exists in a system it means that in effect there can be no way to predict how the system will behave except by going through almost as many steps of computation as the evolution of the system itself" \citep[p739]{RN23}.  

Wolfram defines computational irreducibility by a “practical definition of randomness”, “which is more or less that something is random if the picture” of function/rule iterations “looks random and passes the standard statistical tests for randomness” \citep[209]{nksWolframGray}. \textit{Repetition} and \textit{nesting} do not look random \citep[58,106]{RN23}. Wolfram notes that the Collatz 3$n$+1 function \citep[index]{RN23} is a “simple rule for which unpredictable behavior is found, but the main focus is on proving a simple result about it” \citep[p879, index]{RN23}. “Massively parallel” computations \citep[p999, 1035]{RN23} could bypass the computational irreducibility of serial rules.

\begin{lemma}
\label{lem:wolframno} 
No \textup{a}n+\textup{b} conjecture is undecidable by Wolfram serial path irreducibility. Proof:
\end{lemma}
\noindent All infinite breadth-first and depth-first serial paths of arrows to number pair nodes in the Cayley graph are also infinite sequences of parallel, random-access arrows to same-numbered sites.   
\begin{itemize}
\item the serial, breadth-first, successor path $\tcb{\to}\tcr{\to}\tcb{\to}\tcr{\to}\tcb{\to}\tcr{\to}\cdots$ of the trivial root $\cabroot$  to successive alternating blue, $U$pward, $S_j$ successor generation number pairs, and red, $L$eftwward, $s_i$ successor generation number node pairs (Eq.\ref{eq:sorS}, Lemma \ref{lem:2comb}) are also parallel, random access, arrows 
$\tcb{\uparrow}\tcr{\downarrow}\tcb{\uparrow}\tcr{\downarrow}\tcb{\uparrow}\tcr{\downarrow}\cdots$ to same-numbered sites (Eq.\ref{eq:bfpar}) 
\begin{itemize}
\item blue $U$pward generations $S_j$ whose breadth-first serial paths ($\tcb{\to}\tcb{\to}\tcb{\to}\cdots$, Eq. \ref{eq:dfSj}) to cotree $T_j$ nodes are sequential parallel, random access, arrows to same-numbered sites ($\tcb{\uparrow}\tcb{\uparrow}\tcb{\uparrow}\cdots$, Eq.\ref{eq:dfSjpar}) 
\item red $L$eftward generations $s_i$ whose breadth-first serial paths   ($\tcr{\to}\tcr{\to}\tcr{\to}\cdots$, Eq. \ref{eq:dfsi}) to cotree $t_i$ nodes are sequential parallel, random access, arrows to same-numbered sites  ($\tcr{\downarrow}\tcr{\downarrow}\tcr{\downarrow}\cdots$ (Eq.\ref{eq:dfsipar}). 
\end{itemize}
\end{itemize}
In a mathematical tree traversal, parallel, random-access, arrows are traversed instantaneously and simultaneously, thereby disproving the undecidability of a$n$+b conjectures by serial path irreducibility.\qed

\subsubsection{Cantor 1$n$+0 binary tree, Gödel colons, Conway 1$n$+1 tree, decidable Wolfram $T_L$ ruliad}
\label{sec:cartesian}

The disproofs of undecidability by Cantor diagonals, Gödel numbers, Conway unpredictability and Wolfram irreducibility (Lemmas \ref{lem:cantorno}-\ref{lem:wolframno}) justify:
\begin{enumerate}
\item pairing the arrows ($E$dges), nodes ($V$ertices) and sites (sur$F$aces) within a single Cayley $n$=an$_{\text{odd}}$+b tree $T_{\ge0}$ (Figs.\ref{fig:tge0},\ref{fig:treeP1},\ref{fig:tree7pm1},\ref{fig:tree3m1},\ref{fig:tree5p1}), occupying one column of the infinite Turing–Leibniz–Wolfram-ruliad tape, with those in the binary Cantor 1$n$+0 root-path number tree \citep[S6]{cantor1883bin};
\item pairing the arrows, nodes and sites of the $0$-root Conway 1$n$+1 even-numbers tree with those of Cayley a$n$+b successor trees across columns of the decidable infinite Turing–Leibniz tape $T_L$ of proven a$n$+b conjectures;
\item that same tape $T_L$ being a decidable, bi-periodic, $3$-colour, Wolfram ruliad tape \cite{wolfram21rulial};
\item pairing, across rows of the tape, the successive binary Cantor 1$n$+0 root-path numbers, with Turing-computable cyclic roots $c$=1$a$+1=$2^{2^{\eta-1}+1}$ for $\eta=1,2,3,\dots$,  and other $c$-expressible features of proven a$n$+b conjectures, including a lowest trunk node $c^2$, Lean-certified cotree-density Hamiltonian-tour proofs, and root-trunk proofs $c^2<\frac{c}{2}(c-1)^2$, by means of attributions, citations, number and colour functions not signalled by an implicit colon $:$ as in Gödel  numbers of the liar paradox and the unique paradox yielding undecidability by Gödel numbers (Lemma \ref{lem:godelno}), but by the explicit colon $:$.
\end{enumerate}
See Tab.\ref{fig:TLEE} of the decidable Turing-Leibniz tape $T_L$, showing in its rows an infinite number $\mathbb{N}_0$ of cyclic roots $c=1a+1$ of Cayley trees $T_{\ge0}$ for proven rather than refuted Collatz/Conway a$n$+b conjectures, and showing in its columns paired $c$-expressible parameters, including the column numbers $\eta_{0,1,2,\dots}$ and including 'cotree density' proofs of Hamiltonian tours from and back to $c$, and 'root-trunk' proofs of lowest trunk nodes $c^2$. 

\subsubsection{The 3-regular Cayley tree $T_{\ge0}$ as the 4-regular integer Cartesian upper half graph $C_H$}
\label{sec:cartesian}

\begin{figure}[htb!]
\caption{The Hamiltonian tour of Collatz $3n+1$ branching numbers in the 4-regular Cartesian graph $C_H$} 
\label{fig:hamiltonCo}
\noindent \includegraphics[width=1\textwidth]{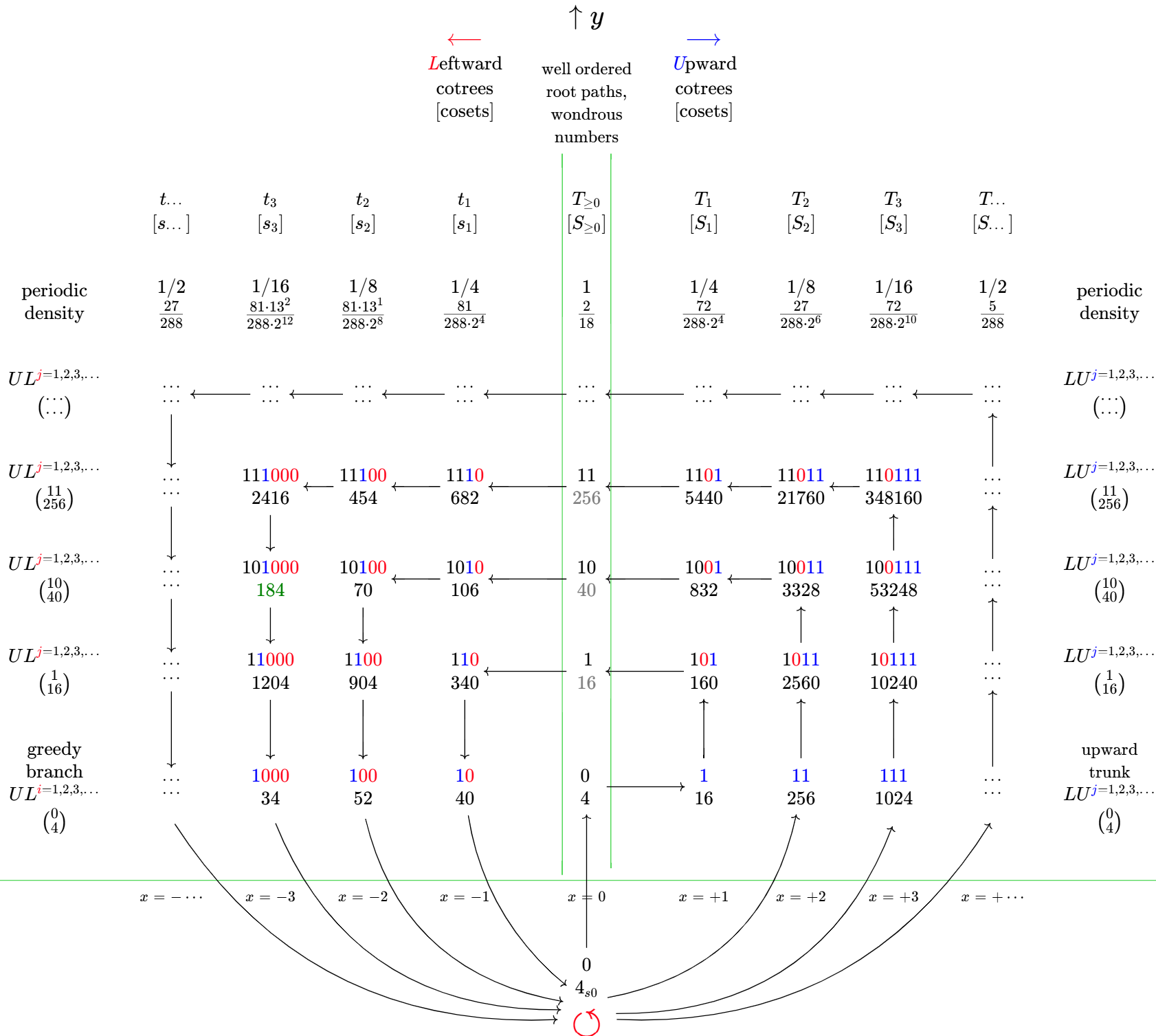}
\end{figure}

\begin{figure}[h!]
\caption{The Hamiltonian tour of Bařina $7n\pm1$ branching numbers in the 4-regular Cartesian graph $C_H$} 
\label{fig:hamiltonBa}
\noindent \includegraphics[width=1\textwidth]{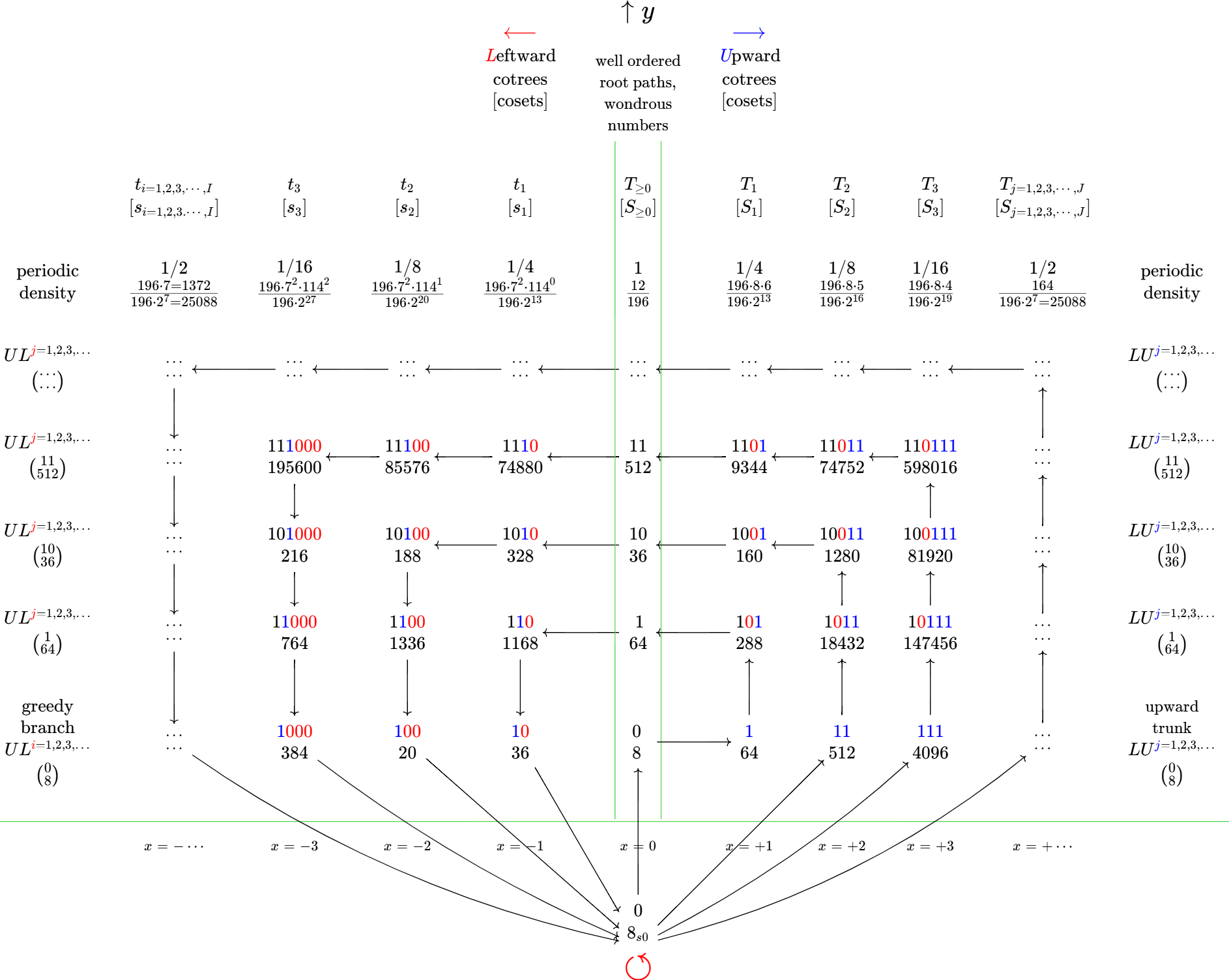}
\end{figure}

\begin{theorem}
\label{the:cartesian}
The 3-regular, left-right-asymmetric, non-commutative/non-confluent Cayley a$n$+b tree  $T_{\ge0}\subseteq G_{\ge0}$ (Figs.\ref{fig:tge0},\ref{fig:treeP1},\ref{fig:tree7pm1},\ref{fig:tree3m1},\ref{fig:tree5p1}, Theorem \ref{the:caytree}) rooted in the cyclic trivial root number pair $\croot\sm{\textcolor{red}{0}\\c=a\cdot1+b}$ (Eq.\ref{eq:Croot}) shows up as the 4-regular, left-right-symmetric, commutative/confluent integer coordinates filling \textit{Cartesian half graph} $C_H$ \textit{with }$\cabroot$ \textit{as its cyclic root and origin }\textup{(0,0) }(Figs.\ref{fig:hamiltonCo},\ref{fig:hamiltonBa}), providing a :
\textup{
\begin{enumerate}
\item a \textbf{Hamiltonian tour} traversing in increasingly wider quadrangles all its coordinates $(x,y)$ from and back to $\sm{\textcolor{red}{0}\\c=a\cdot1+b}$ at its origin $(0,0)$, collecting blue $U$pward cotree number pairs from Quadrant I, and red $L$eftward cotree number pairs from Quadrant II, while skipping $y$-axis coordinates $(0,y)$ (Figs.\ref{fig:hamiltonCo},\ref{fig:hamiltonBa}).
\item a \textbf{cotree density test} whether the cumulative periodic density $d(T_{\ge0})$ of blue $U$pward successor number pairs from Quadrant I columns and of red $L$eftward successor number pairs from Quadrant II columns is equal to, or lower than, the periodic density
$d(S_{\ge0})$ of all $n$=a$n_{\textrm{odd}}$+b numbers indivisible by $a$ (top row in Figs.\ref{fig:hamiltonCo},\ref{fig:hamiltonBa}). Upcoming subsections (Eqs.\ref{eq:dSge0},\ref{eq:cdp}) specify this proof or refutation.
\end{enumerate}
}
\end{theorem}

\noindent\textit{Proof / Explanation: }
Each number pair $\sm{y\\n}$ in Cayley tree $T_{\ge0}$ shows up at two distinct integer coordinates $(x,y)$ in the Cartesian half graph $C_H$:
\begin{itemize}
    \item at its cyclic root and at is origin (0,0) the cyclic root $L:\croot\sm{\textcolor{red}{0}\\c=a\cdot1+b}\tcr{\to}\sm{\textcolor{red}{0}\\c=a\cdot1+b}$ (Eq.\ref{eq:Croot}).
\item At $(0,y)$ coordinates along its $y$-axis the well-ordered, breadth-first ordered, number pairs $\sm{y\\n}$ in $T_{\ge0}$ (Eq.\ref{eq:sorS}). The serial, breadth-first, successor path $\tcb{\to}\tcr{\to}\tcb{\to}\tcr{\to}\tcb{\to}\tcr{\to}\cdots$ of the trivial root $\cabroot$  to successive alternating blue, $U$pward, $S_j$ successor generation number pairs, and red, $L$eftwward, $s_i$ successor generation number pairs (Eq.\ref{eq:sorS}, Lemma \ref{lem:2comb})  matches parallel, random access, arrows 
$\tcb{\uparrow}\tcr{\downarrow}\tcb{\uparrow}\tcr{\downarrow}\tcb{\uparrow}\tcr{\downarrow}\cdots$ (Eq.\ref{eq:bfpar}) to these successive (0,y)-coordinates along its $y$-axis.
\item at ($x=j,y$) coordinates in successive Quadrant I columns (Eq.\ref{eq:j+x}) the blue $U$pward successor number pairs in successive generations $S_{j=1,2,3,\dots}$ / cotrees $T_{j=1,2,3,\dots}$ (Eqs.\ref{eq:dfSj},\ref{eq:dfSjpar}),
\item at ($x=-i,y$) coordinates in successive Quadrant II columns (Eq.\ref{eq:i-x}) the red $L$eftward successor number pairs in successive generations $s_{i=1,2,3,\dots}$ / cotrees $t_{i=1,2,3,\dots}$ (Eqs.\ref{eq:dfsi},\ref{eq:dfsipar}). 
\item at successive ($x=j,y=0$) coordinates on the $x$-axis in Quadrant I successive upward trunk root pairs being cyclic roots of blue $U$pward cotrees $T_j$. At successive $(x=j,y)$ coordinates in Quadrant I breadth-first ordered number pairs $\sm{x=j\\n\in T_j}$ in blue $U$pward cotrees $T_{j}$ (Fig.\ref{fig:utree}bdf). Breadth-first traversals of blue $U$pward generations $S_j$ / cotrees $T_j$ ($\tcb{\to}\tcb{\to}\tcb{\to}\cdots$, Eq. \ref{eq:dfSj}) match parallel, random access, arrows ($\tcb{\uparrow}\tcb{\uparrow}\tcb{\uparrow}\cdots$ to $T_j$ number pairs (Eq.\ref{eq:dfSjpar}) in successive Quadrant I columns (Eq.\ref{eq:j+x}).  

\item at successive ($x=-i,y=0$) coordinates on the $x$-axis in Quadrant II successive greedy branch root pairs being cyclic roots of red $L$eftward cotrees $t_j$. At successive $(x=-i,y)$ coordinates in Quadrant II breadth-first ordered number pairs $\sm{x=i\\n\in t_i}$ in red $L$eftward cotrees $t_{i}$ (Fig.\ref{fig:ltree}bdf). Serial breadth-first traversals of red $L$eftward generations $s_i$ / cotrees $t_i$   ($\tcr{\to}\tcr{\to}\tcr{\to}\cdots$, Eq. \ref{eq:dfsi}). match parallel, random access, arrows  ($\tcr{\downarrow}\tcr{\downarrow}\tcr{\downarrow}\cdots$ in successive Quadrant II columns (Eq.\ref{eq:dfsipar})  (Eq.\ref{eq:i-x}). 
    \end{itemize}

\newpage

\noindent The Cartesian half graph $C_H$ 
(Figs.\ref{fig:hamiltonCo},\ref{fig:hamiltonBa}) is left-right-symmetric about its $y$-axis.  $C_H$ is commutative/confluent as an infinity of paths (commutative computation sequences) comprising row arrows $x\leftrightarrow x\pm1$ and column arrows $y\leftrightarrow y\pm1$ is available for an infinity of paths of consecutive arrows from one node (confluing) to another. $C_H$ identifies (formally) as a 4-regular graph of which the (redundant) coordinates $(0,-i)$ on its negative $y$-axis enable for $i=1,2,3,\dots$ non-crossing, non-confluing Hamiltonian tour paths from Quadrant II back to Quadrant I  $(x,-i)\to(0,-i)\to(0,i+1)$, and back to the trivial root $(0,-i)\;\to\;(0,-i+i)$.

\begin{enumerate}
\item The 4-regular Cartesian half graph enables a Hamiltonian tour in increasingly wider quadrangles from and (confluing) back to the trivial root $\cabroot$ at its origin (0,0) as the tour can choose at each coordinate $(x,y)$ two out of four arrows to adjacent nodes  $(x,y)\leftrightarrow(x\pm1,y\pm1)$ (Figs.\ref{fig:hamiltonCo},\ref{fig:hamiltonBa}). 
\item the cotree density test $d(T_{\ge0})\stackrel{?}{=}d(S_{\ge0})$ (Eqs.\ref{eq:cdp},\ref{eq:dSge0}) whether the cumulative periodic $d$ensity $d(T_{\ge0})$ of blue $U$pward successor number pairs from Quadrant I columns and of red $L$eftward successor number pairs from Quadrant II columns equals the $d$ensity of set $S_{\ge0}$ of $a\nmid$a$n_{\textup{odd}}$+b numbers: 
\begin{itemize}
    \item 
requires densities per blue $U$pward cotree column in Quadrant I and per $L$eftward cotree column in Quadrant II, displayed for a few cotree columns in its top row (Figs.\ref{fig:hamiltonCo},\ref{fig:hamiltonBa})
    \item giving the cumulative density of $a-1$ subseries of $U$pward cotree densities and their sum, displayed to the right of its top row (Figs.\ref{fig:hamiltonCo},\ref{fig:hamiltonBa}), and the cumulative density of $1$ series of $U$pward cotree densities, displayed to the left of its top row (Figs.\ref{fig:hamiltonCo},\ref{fig:hamiltonBa}) \qed
\end{itemize}

\end{enumerate}

\newpage
\subsection {Eq. series 4: conjugate inner arrows cotrees, conjugate outer paths to cotrees}
\label{sec:34eqs}

\begin{flalign} 
&\noindent\textrm{Tree $T_{\ge0}$ provides \textbf{automorphism graph }}\; \textrm{Aut} \left(T_{\ge0},\left[\sm{\textcolor{red}{0}\\L},\,\sm{\textcolor{blue}{1}\\U}\right]\right) \textrm{. Nodes: }T_{\ge0} \textrm{, its subtrees, cotrees}\label{eq:auttge0}&&
\end{flalign}
\vspace{-0.7cm}
\begin{align}
&\textrm{leftward path:}\notag\\
&\;\textrm{ leftward subtree}
&\sm{\textcolor{red}{0}^{i=1,2,3,\dots}\\L^i}:\;
&T_{\ge 0}\left(S_{\ge 0},\left[\sm{\textcolor{red}{0}\\L},\sm{\textcolor{blue}{1}\\U}\right]\right) \;\textcolor{red}{\;\to}\quad\;\; t_{0;\ge{i}}\left(s_{0;\ge{i}},\left[\sm{\textcolor{red}{0}\\L_{\textcolor{white}{i}}},\sm{\textcolor{blue}{1}\\U_i}\right]\right) 
&\textrm{(Fig.\ref{fig:ltree}ace)} 
\label{eq:lsuar}\\ 
&\quad \textrm{conjugate arrows}
& &\left[\sm{\textcolor{red}{0}^{\textcolor{white}{i}}\\L^{\textcolor{white}{i}}}=\sm{0^{-i}0 0^i\\L^{-i}LL^i},\quad\;\; \sm{\textcolor{blue}{1}^{\textcolor{white}{i}}\\U_i^{\textcolor{white}{i}}} =\sm{0^{-i}10^i\\L^{-i}UL^i} \right]
\nonumber \\ 
&\textrm{upward path:}\notag\\
&\;\textrm{ upward subtree}
&\sm{\textcolor{blue}{1}^{j=1,2,3,\dots}\\U^j}:\;
&T_{\ge 0}
\left(S_{\ge 0},\left[\sm{\textcolor{red}{0}\\L},\sm{\textcolor{blue}{1}\\U}\right]\right)\; 
\textcolor{blue}{\;\to}\quad\; T_{\ge{j}}
\left(S_{\ge{j}},\left[\sm{\textcolor{red}{0}\\L_{\textcolor{black}{j}}},
\sm{\textcolor{blue}{1}_{\textcolor{white}{i}}\\U_{\textcolor{white}{i}}}\right]\right)
&\textrm{(Fig.\ref{fig:utree}ace)}
\label{eq:usuar}\\ 
&\quad  \textrm{conjugate arrows}
& &\left[ \sm{\textcolor{red}{0}^{\textcolor{white}{i}}\\L_{\textcolor{black}{j}}^{\textcolor{white}{i}}} =\sm{1^{-j}01^j\\U^{-j}LU^j},\quad\;\;
\sm{\textcolor{blue}{1}_{\textcolor{white}{i}}^{\textcolor{white}{i}}\\U_{\textcolor{white}{i}}}=\sm{1^{-j}11^j\\U^{-j}UU^j}\right]
\quad\quad 
\nonumber \\
&\textrm{upw.arrow}\to\textrm{leftw.path}& &\textrm{or, leftward path}\to\textrm{conjugate upward arrow:}\notag\\
&\;\textrm{ red \textit{L}eftward cotrees }t_i
&\sm{\textcolor{blue}{1}\textcolor{red}{0}^{i=1,2,3,\dots}\\UL^i}:\;
&T_{\ge 0}\left(S_{\ge 0},\left[\sm{\textcolor{red}{0}\\L},\sm{\textcolor{blue}{1}\\U}\right]\right)
\textcolor{red}{\;\to}\quad 
t_{i}\left(s_{i},\left[\sm{0_{1;i}\\L_{1;i}},\sm{1_i\\U_i}\right]\right)
&\textrm{(Fig.\ref{fig:ltree}bdf)}
\label{eq:lcoa}\\ 
&
&\sm{\textcolor{red}{0}^{i}\textcolor{blue}{1}_i\\L^iU_i}:\;
&T_{\ge 0}\left(S_{\ge 0},\left[\sm{\textcolor{red}{0}\\L},\sm{\textcolor{blue}{1}\\U}\right]\right)
\textcolor{red}{\;\to}\quad 
t_{i}\left(s_{i},\left[\sm{0_{1;i}\\L_{1;i}},\sm{1_i\\U_i}\right]\right)
\notag\\ 
&\quad \textrm{red conjugate arrows}
& &\left[ \sm{0_{1;i}\\L_{1;i}} =\sm{(0^{-i} 1^{-1})0(1^10^i)\\(L^{-i} U^{-1})L(U^1L^i)},\quad\;\; \sm{1_i\\U_i}=\sm{0^{-i}10^i\\L^{-i}UL^i} \right]
&\textrm{(Fig.\ref{fig:ltree}bdf)}
\label{eq:lcoar}
\\
&\quad\Rightarrow\textrm{Theorem \ref{the:cartesian}:}&&\textrm{breadth-first $t_i$ pairs $(-x=i,y)$ in $QII\subset H_C$}
\label{eq:i-x}
\\
    &\textrm{leftw.arrow}\to\textrm{upw.path}& &\textrm{or, upward path}\to\textrm{conjugate leftward arrow:}\notag\\
&\;\textrm{ blue \textit{U}pward cotrees }T_j
&\sm{\textcolor{red}{0}\textcolor{blue}{1}^{j=1,2,3,\dots}\\LU^j}:\;
&T_{\ge 0}\left(S_{\ge 0},\left[\sm{\textcolor{red}{0}\\L},\sm{\textcolor{blue}{1}\\U}\right]\right)
\textcolor{blue}{\;\to}\quad T_j\left(S_{j},\left[\sm{0_j\\L_j},\sm{1_{1;j}\\U_{1;j}}\right]\right) &\textrm{(Fig.\ref{fig:utree}bdf)} 
\label{eq:ucoa}
\\& &\sm{\textcolor{blue}{1}^{j}\textcolor{red}{0}_j\\U^jL_j}:\;
&T_{\ge 0}\left(S_{\ge 0},\left[\sm{\textcolor{red}{0}\\L},\sm{\textcolor{blue}{1}\\U}\right]\right)
\textcolor{blue}{\;\to}\quad T_j\left(S_{j},\left[\sm{0_j\\L_j},\sm{1_{1;j}\\U_{1;j}}\right]\right)
\notag
\\
&\quad\textrm{blue conjugate arrows}
& &\left[ \sm{0_j\\L_j}=\sm{1^{-j}01^j\\U^{-j}LU^j},\quad\sm{1_{1;j}\\U_{1;j}}=
\sm{(1^{-j}0^{-1})1(0^1 1^j)\\(U^{-j}L^{-1})U(L^1 U^j)}\right]
&\textrm{(Fig.\ref{fig:utree}bdf)}
\label{eq:ucoar} 
\\
&\quad\Rightarrow\textrm{Theorem \ref{the:cartesian}:}&&\textrm{breadth-first $T_j$ pairs $(x=j,y)$ in $QI\subset H_C$}\label{eq:j+x}
\end{align}

\noindent In the decades before and after the turn to the 20th century, 'isomorphism' of geometric shapes came to extend their Euclidean 'congruence'. The informal definition of \textit{isomorphic} graphs is that all their nodes (vertices) and all their arrows (edges) can be draped over each other so that they completely cover each other. Subsection \ref{sec:34eqs} on Eqs.\ref{eq:auttge0}-\ref{eq:j+x} covers \textit{conjugate inner arrows} in four infinite series of isomorphic leftward subtrees (Fig.\ref{fig:ltree}ace), leftward cotrees (Fig.\ref{fig:ltree}bdf), upward subtrees ( Fig.\ref{fig:utree}ace), and upward cotrees (Fig.\ref{fig:utree}bdf) between indirectly connected nodes via the Cayley a$n$+b tree $T_{\ge0}$ itself. \textit{Conjugate outer root paths} to leftward cotrees (Eq.\ref{eq:dfsiconjug}) and upward cotrees (\ref{eq:dfSjconjug}), also defined in terms of Max Dehn’s isomorphism formula \citep[p.11]{RN19} (Eqs.\ref{eq:lsuar}-\ref{eq:j+x}), also define the automorphism graph Aut($T_{\ge0},[L,U]$) (Fig.\ref{fig:aut}). 

\subsubsection {The automorphism graph with two root paths to each cotree}

\begin{figure}[t!]
\caption{Collatz $a=3$ subtree-cotree ladders automorphism graph Aut$\left(T_{\ge 0},[L,U]\right)$}
\label{fig:aut}
\begin{center}
   \includegraphics[width=0.9\textwidth]{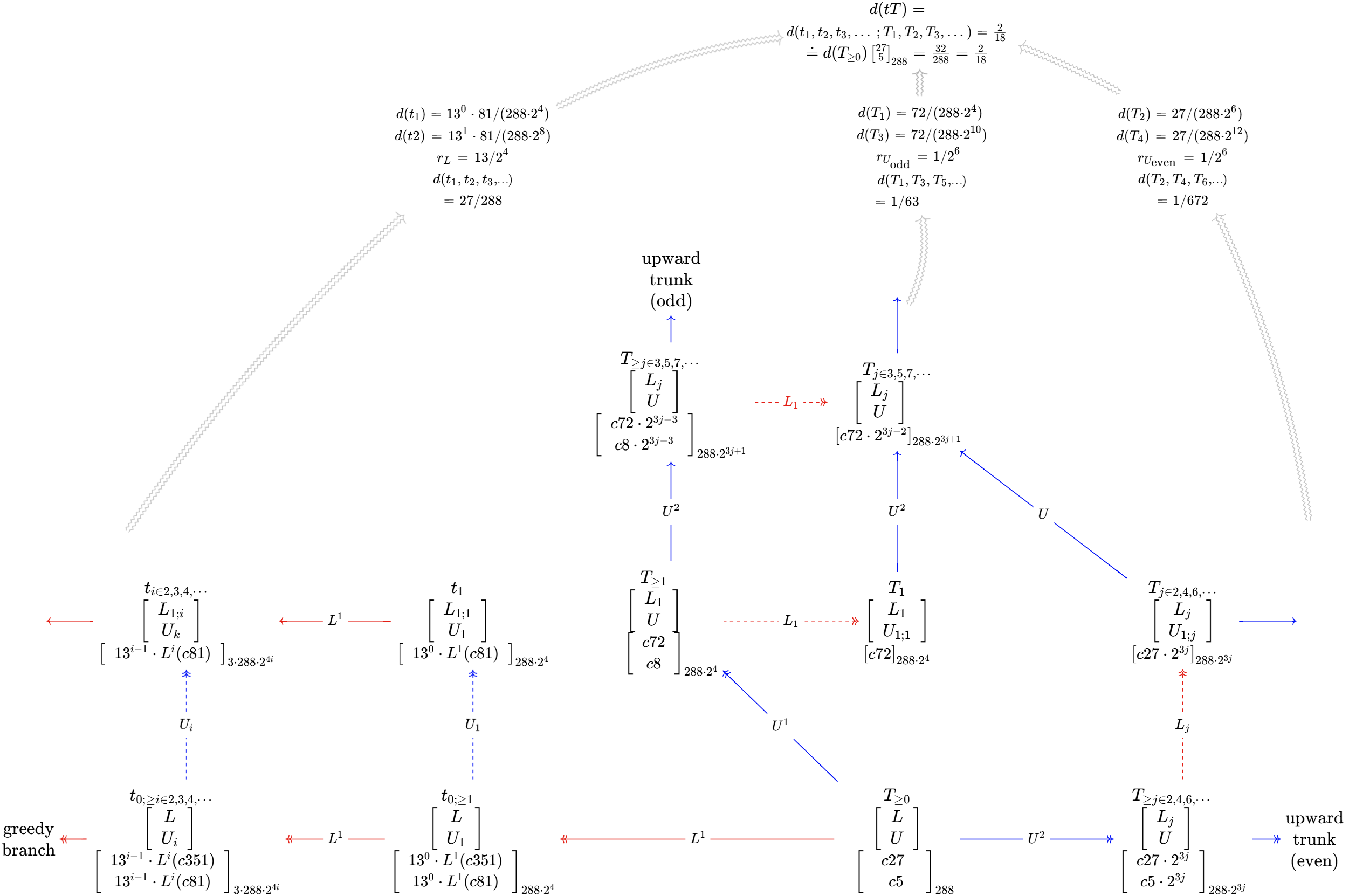}
\end{center}
\centering\begin{minipage}{0.9\linewidth}
\noindent \textbf{Legend \textbar }  The Collatz $3n+1$ automorphism graph Aut$(T_{\ge 0},[L,U])$ has as nodes isomorphic trees (Figs.\ref{fig:tge0}a, \ref{fig:ltree}, \ref{fig:utree}), and as arrows $L-$ and $U-$functions connecting them (Eqs.\ref{eq:lsuar}, \ref{eq:lcoar},\ref{eq:usuar}, \ref{eq:ucoar}). It shows $a=3$ \textit{infinite Jacob ladders}: one of \textit{L}eftward subtrees and cotrees (Fig.\ref{fig:ltree}), one of odd \textit{U}pward subtrees and cotrees and one of even $U$pward subtrees and cotrees  (Fig.\ref{fig:utree}). The 3-ladder Jacob's graph appears to be an infinite version of the finite 2-wing \textit{butterfly} graph known from later visualisations of the isomorphism lemma formulated by Zassenhaus \cite {RN3153}.

$\quad$Arrows vs. brackets: outer vs. inner automorphism. Red vs. blue arrows: leftward vs. upward arrows.  Single vs. double arrowheads ($\to $ vs. $\twoheadrightarrow $: ordinary vs. conjugate cotree generation). Lowercase $t$ vs uppercase $T$: leftward  vs. upward subtrees and cotrees. Superscript vs. subscript, solid vs. dashed arrows: iterates $L^i$ and $U^j$ vs. conjugates $L_j$ and $U_i$. Table \ref{tab:classes} enumerates the sets of congruence classes  $c27,c5,c72,c8,c351$ and $c81$ of the trees. Notations like ${13}^{i-1}\cdot U^i(c81)$ indicate that${13}^{i-1}\cdot 81$ subclasses modulo 54, to which Def.\ref{def:Lw} applies, are obtained by applying $L$ to the previous ${13}^{i-2}\cdot 81$ subclasses. The annotations with cotree densities $d\left(t_1\right),\ d\left(T_1\right)$ and $d\left(T_2\right)$ reveal the cumulative cotree density $2/18$, which is equal to the periodic density of all branching numbers $d(S_{\ge 0})={d(\left[4,16\right]}_{18})=\frac{2}{18}$. 
The top line in the graph shows this cotree density proof $\Bigl( d(tT)= \frac{27}{288}+\frac{1}{63}+\frac{1}{672}=\frac{32}{288}=\frac{2}{18}\Bigr)\doteq \Bigl( d(S_{\ge0}[4,16]_{18})=\frac{2}{18}\Bigr)$ (Eq.\ref{eq:cdp}).$\quad\textbf{\textbar}$ 
\end{minipage}\end{figure}

\normalsize{}

Based on his lectures on Cayley color graphs \cite{peiferCayleyDehn}, Max Dehn (1878-1952) introduced in 1911 the transformation formula $U_1=U^{+1}TU^{-1}$, in diagrammatic, leftmost-innermost, notation $U_1=U^{-1}TU^{+1}$, generalized to $U_i=U^{-i}TU^{+i}$. 
The start node and the arrowhead node of arrow $T$ in a primal graph both having a path $U^i$ to nodes also included in the isomorphic graph, in which arrow $T$ is transformed to a conjugate arrow defined as the roundabout path $U_i=U^{-i}TU^{+i}$ via the start node and the arrowhead node of arrow $T$ in the primal graph. A Collatz-Cayley binary tree $T_{\ge0}\subseteq G_{\ge0}$ (Eq.\ref{eq:tleg}) has two biperiodic primal arrows $T$: the red $L$eftward arrow $\sm{\textcolor{red}{0}\\L=fg^p}$ and the blue $U$pward arrow $\sm{\tcb{1}\\U=g^q}$. 
Four different iterators (Eq.\ref{eq:dtT}) specify inner arrows in cotrees (Eq.\ref{eq:cdp}). Conjugate paths of arrows are not specified in this way in `algebraic graph theory' starting from an infinite binary tree without a cyclic root \cite{RN11,RN20,grig}.

\begin{figure}[h!]
\caption{Red \textit{L}(eftward) subtrees $t_{0;\ge i=1,2,3}$ and same-\textit{s}-generation \textit{L} cotrees $t_{i=1,2,3}$}
\label{fig:ltree}
\noindent \centerline{\includegraphics[width=0.85\textwidth]{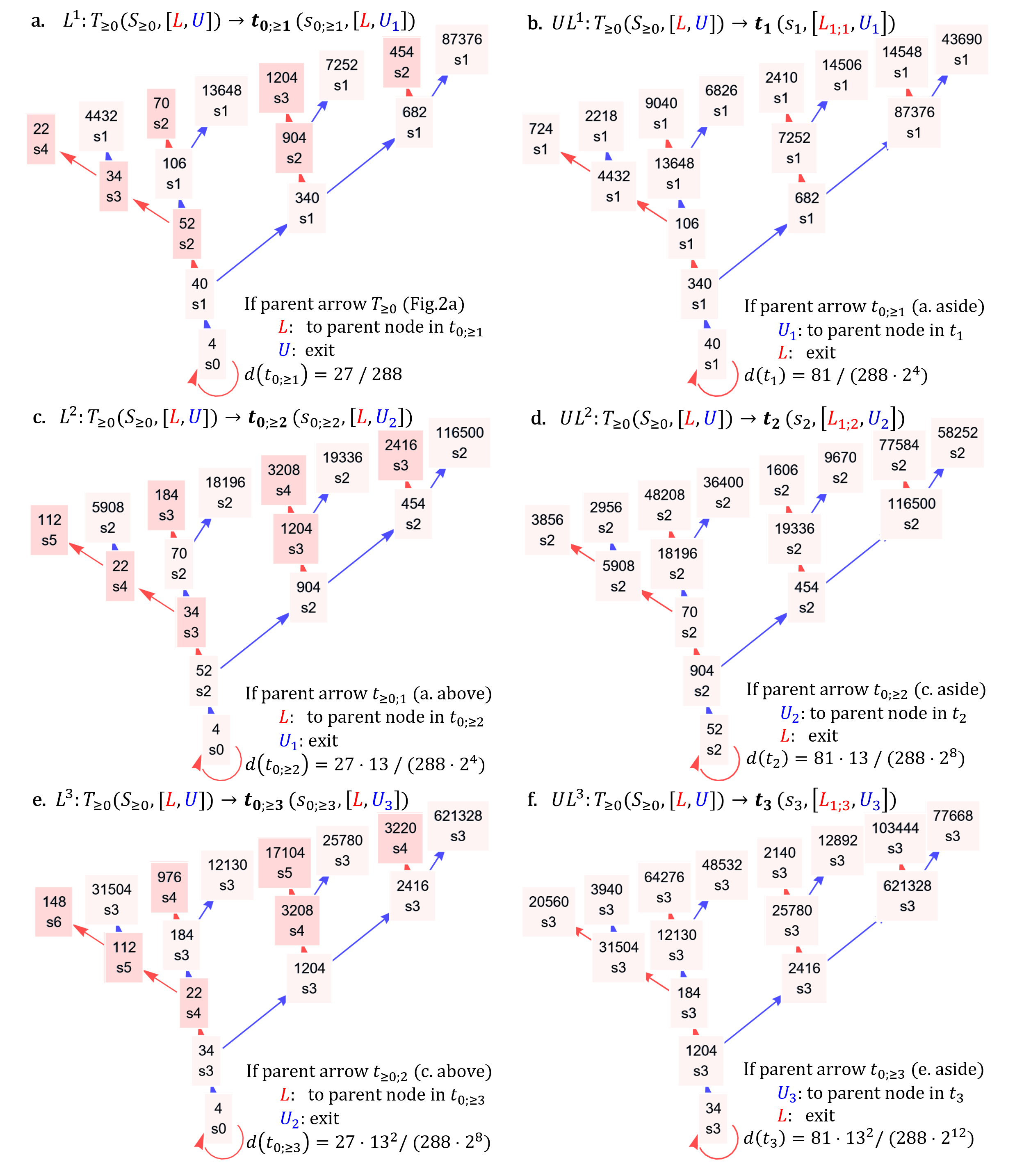}}
\noindent \textbf{Legend \textbar } Leftward moves $L^{i=1,2,3,\dots }$ of tree $T_{\ge 0}$ (Fig.1b) give leftward subtrees with nested generations (a) $s_{0;\ge 1}$ in $t_{0;\ge 1}$, (c) $s_{0;\ge 2}$ in $t_{0;\ge 2}$, (e) $s_{0;\ge 3}$ in $t_{0;\ge 3}$. Moves $UL^{i=1,2,3,\dots }=L^{i=1,2,3,\dots }U_i$ give leftward cotrees with disjoint generations (b) $s_1$ in $t_1$, (d) $s_2$ in $t_2$, (f) $s_3$ in $t_3$. Leftward subtrees and cotrees obtain conjugate upward arrows  $U_{i=1,2,3,..}=L^{-i}UL^i$; cotrees obtain conjugate leftward arrows $L_{1;i=1,2,3,\dots} =(L^{-i}U^{-1})L(UL^i)$. Assessment periodic $d$ensities: Section \ref{sec:cc}.   

$ \begin{array}{@{}c@{}}
256\stackrel{L}{\textcolor{red}{\rightarrow}}\\ 
{\textcolor{blue}{\uparrow} }^{L}\\ 
16\stackrel{L}{\textcolor{red}{\rightarrow}}\\ 
T_{\ge 0} \end{array}
\begin{array}{@{}c@{}}
340\stackrel{L}{\textcolor{red}{\rightarrow}}\\ 
{\textcolor{blue}{\uparrow} }^{U_1}\\ 
40\stackrel{L}{\textcolor{red}{\rightarrow}}\\ 
t_{0;\ge 1} \end{array}
\begin{array}{@{}c@{}}
904\stackrel{L}{\textcolor{red}{\rightarrow}}\\ 
{\textcolor{blue}{\uparrow} }^{U_2} \\ 
52\stackrel{L}{\textcolor{red}{\rightarrow}}\\ 
t_{0;\ge 2} \end{array}
\begin{array}{@{}c@{}}
1204\stackrel{L}{\textcolor{red}{\rightarrow}}\cdots\\ 
{\textcolor{blue}{\uparrow} }^{U_3}   \\ 
34\stackrel{L}{\textcolor{red}{\rightarrow}}\cdots\\ 
t_{0;\ge 3} \end{array}
\begin{array}{c@{}}\text{;}\; \;\;\;
 \\ 
 \\ 
 \\ 
\end{array}
\begin{array}{@{}c@{}}
40\stackrel{U}{\textcolor{blue}{\rightarrow}} \\ 
{\textcolor{red}{\uparrow} }^{L_{\ }} \\ 
16\stackrel{U}{\textcolor{blue}{\rightarrow}} \\ 
T_{\ge 0} \end{array}
\begin{array}{@{}c@{}}
160\stackrel{L}{\textcolor{red}{\rightarrow}} \\ 
{\textcolor{red}{\uparrow} }^{{L_1}_{\ }} \\ 
256\stackrel{L}{\textcolor{red}{\rightarrow}} \\ 
T_{\ge 1} \end{array}
\begin{array}{@{}c@{}}
106_{s1}\stackrel{L}{\textcolor{red}{\rightarrow}} \\ 
{\textcolor{red}{\uparrow} }^{L_{1;1}} \\ 
340_{s1}\stackrel{U}{\textcolor{red}{\rightarrow}} \\ 
t_{1} \end{array}
\begin{array}{@{}c@{}}
70_{s2}\stackrel{L}{\textcolor{red}{\rightarrow}} \\ 
{\textcolor{red}{\uparrow} }^{L_{1;2}} \\ 
904_{s2}\stackrel{L}{\textcolor{red}{\rightarrow}} \\ 
t_{2} \end{array}
\begin{array}{@{}c@{}}
184_{s3}\stackrel{L}{\textcolor{red}{\rightarrow}}\cdots  \\ 
{\textcolor{red}{\uparrow} }^{L_{1;3}} \\ 
1204_{s3}\stackrel{L}{\textcolor{red}{\rightarrow}}\cdots  \\ 
t_{3} \end{array}
\begin{array}{c@{}}
\;\;\quad\quad\quad\textbf{\textbar} 
\end{array}
$
\end{figure}

\begin{figure}[h]
\caption{Blue \textit{U}(pward) subtrees $T_{\ge j=1,2,3}$ and  and same-\textit{S}-generation \textit{U} cotrees $T_{j=1,2,3}$}
\label{fig:utree}
\noindent \centerline{\includegraphics*[width=0.85\textwidth]{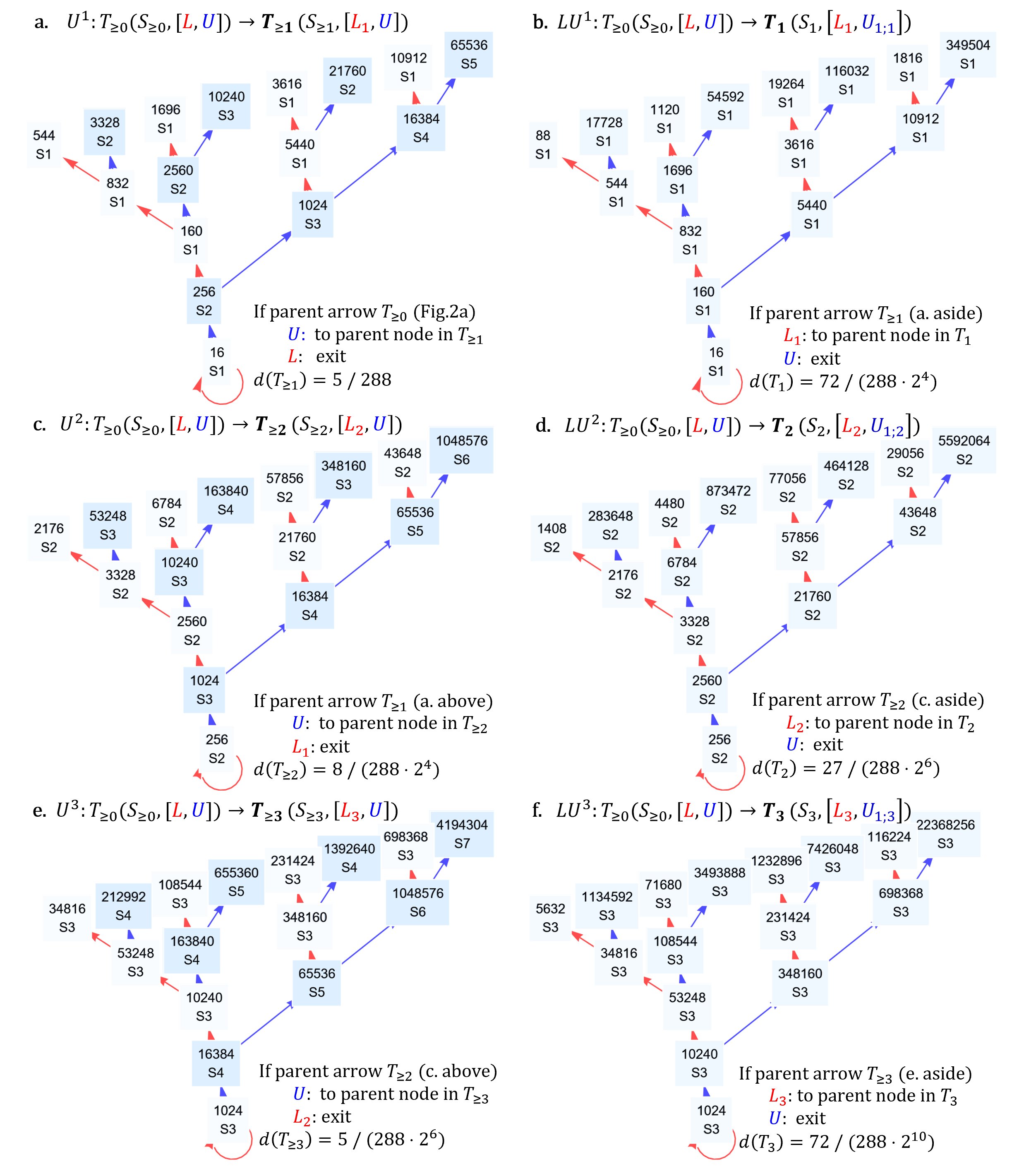}}

\noindent \textbf{Legend \textbar }
Upward moves $U^{j=1,2,3,\dots }$ of tree $T_{\ge 0}$ (Fig.1b) give upward subtrees with nested generations (a) $S_{\ge 1}$ in $T_{\ge 1}$, (c) $S_{\ge 2}$ in $T_{\ge 2}$, (e) $S_{\ge 3}$ in $T_{\ge 3}$. Moves $LU^{j=1,2,3,\dots }$ give upward cotrees with disjoint generations (b) $S_1$ in $T_1$, (d) $S_2$ in $T_2$, (f) $S_3$ in $T_3$. Upward subtrees and cotrees obtain conjugate leftward arrows  $L_{j=1,2,3,..}=U^{-j}LU^j$; cotrees obtain conjugate upward arrows $U_{1;j=1,2,3,\dots }=(U^{-j}L^{-1})U(LU^j)$. Assessment periodic $d$ensities: Section \ref{sec:cc}.  
$ \begin{array}{@{}c@{}}
40\stackrel{U}{\textcolor{blue}{\rightarrow}}\\ 
{\textcolor{red}{\uparrow} }^{L}\\ 
16\stackrel{U}{\textcolor{blue}{\rightarrow}}\\ 
T_{\ge 0} \end{array}
\begin{array}{@{}c@{}}
160\stackrel{U}{\textcolor{blue}{\rightarrow}}\\ 
{\textcolor{red}{\uparrow} }^{L_1}\\ 
256\stackrel{U}{\textcolor{blue}{\rightarrow}}\\ 
T_{\ge 1} \end{array}
\begin{array}{@{}c@{}}
2560\stackrel{U}{\textcolor{blue}{\rightarrow}}\\ 
{\textcolor{red}{\uparrow} }^{L_2} \\ 
1024\stackrel{U}{\textcolor{blue}{\rightarrow}}\\ 
T_{\ge 2} \end{array}
\begin{array}{@{}c@{}}
10240\stackrel{U}{\textcolor{blue}{\rightarrow}}\cdots\\ 
{\textcolor{red}{\uparrow} }^{L_3}   \\ 
16384\stackrel{U}{\textcolor{blue}{\rightarrow}}\cdots\\ 
T_{\ge 3} \end{array}
\begin{array}{c@{}}\text{;}\; \;\;\;
 \\ 
 \\ 
 \\ 
\end{array}
\begin{array}{@{}c@{}}
256\stackrel{L}{\textcolor{red}{\rightarrow}} \\ 
{\textcolor{blue}{\uparrow} }^{U_{\ }} \\ 
16\stackrel{L}{\textcolor{red}{\rightarrow}} \\ 
T_{\ge 0} \end{array}
\begin{array}{@{}c@{}}
340\stackrel{U}{\textcolor{blue}{\rightarrow}} \\ 
{\textcolor{blue}{\uparrow} }^{{U_1}_{\ }} \\ 
40\stackrel{U}{\textcolor{blue}{\rightarrow}} \\ 
t_{0;\ge 1} \end{array}
\begin{array}{@{}c@{}}
5440_{S1}\stackrel{U}{\textcolor{blue}{\rightarrow}} \\ 
{\textcolor{blue}{\uparrow} }^{U_{1;1}} \\ 
160_{S1}\stackrel{U}{\textcolor{blue}{\rightarrow}} \\ 
T_{1} \end{array}
\begin{array}{@{}c@{}}
21760_{S2}\stackrel{U}{\textcolor{blue}{\rightarrow}} \\ 
{\textcolor{blue}{\uparrow} }^{U_{1;2}} \\ 
2560_{S2}\stackrel{U}{\textcolor{blue}{\rightarrow}} \\ 
T_{2} \end{array}
\begin{array}{@{}c@{}}
348160_{S3}\stackrel{U}{\textcolor{blue}{\rightarrow}}\cdots  \\ 
{\textcolor{blue}{\uparrow} }^{U_{1;3}} \\ 
10240_{S3}\stackrel{U}{\textcolor{blue}{\rightarrow}}\cdots  \\ 
T_{3} \end{array}
\begin{array}{c@{}}
\quad\quad\textbf{\textbar}\\ 
 \\ 
 \\ 
\end{array}
$
\end{figure}

The automorphism graph Aut($T_{\ge0},[L,U]$) shows $a$-Jacob's subtree-cotree-ladders, 1 red $L$eftward subtree-cotree ladder, and $a-1$ blue $u$pward subtree-cotree ladders  (Fig.\ref{fig:aut}).
Rungs connect the stile of red leftward subtrees of higher generation $s$-successors generated by $L^{i=1,2,3,\dots}(T_{\ge0})$) (Fig. \ref{fig:ltree}ace)), with the stile of red leftward, same-generation cotrees generated by $UL^{i=1,2,3,\dots}(T_{\ge0})$) iterations (Figs.\ref{fig:ltree}bdf). The $a-1$ infinite ascending ladders have rungs connecting isomorphic, blue upward subtrees of higher generation $S$-successors generated by $U^{j=1,2,3,\dots}(T_{\ge0})$ (Figs. \ref{fig:utree}ae,c$\cdots$) with blue upward, same-generation, cotrees generated by $LU^{j=1,2,3,\dots}(T_{\ge0})$ iterations (Figs. \ref{fig:utree}bf,d$\cdots$). 

Dehn's \textit{isomorphism} formula defines conjugate inner arrows of subtrees and cotrees and conjugate outer arrows to leftward cotrees ($UL^{i=1,2,3,\dots}(T_{\ge0})=L^{i=1,2,3,\dots}U_i(T_{\ge0})$) and upward cotrees ($LU^{j=1,2,3,\dots}(T_{\ge0})=U^{j=1,2,3,\dots}L_j(T_{\ge0})$) (Eqs.\ref{eq:lcoa},\ref{eq:ucoa}). 
The leftward and upward conjugate outer paths $L^iU_i:T_{\ge0}\to t_i$ and $U^jL_j:T_{\ge0}\to T_j$ to generate from the Cayley tree $T_{\ge0}$ leftward cotrees $t_{i=1,2,3,\dots}$ and 
upward cotrees $T_{j=1,2,3,\dots}$ via leftward subtrees $t_{0;i\ge1,2,3,\dots}$ and upward subtrees  $T_{j\ge1,2,3,\dots}$ visualize leftward subtree (Fig.\ref{fig:ltree}ace) - cotree (Fig.\ref{fig:ltree}bde) ladder graphs and upward subtree (Fig.\ref{fig:utree}ace) - cotree (Fig.\ref{fig:utree}bde) ladder graphs. The subtree-cotree rungs are same-generation nodes within subtrees, connected within cotrees by a conjugate bypass via ancestor arrows in $T_{\ge0}$.   

The paths in its automorphism graph Aut$(T_{\ge0},[\sm{\textcolor{red}{0}\\L=fg^p},\sm{\textcolor{blue}{1}\\U=g^q}])$ (Fig.\ref{fig:aut}, Eq.\ref{eq:auttge0}) show four iterated functions generating four infinite series of generated isomorphic graphs---upward subtrees, leftward cotrees, leftward subtrees and upward cotrees---each of them with unique conjugate inner arrows, as specified in Equations \ref{eq:lsuar}-\ref{eq:ucojcm}). 
\begin{align}
      U^{j=1,2,3,\cdots} &\textrm{ onto upward subtrees } \quad T_{\ge j} &\textrm{(Eq.\ref{eq:usuar}, Fig.\ref{fig:utree}ace)} \notag\\
      UL^{i=1,2,3,\cdots} &\textrm{ onto leftward cotrees } \quad\quad t_{i} \textrm{, red nodes generation }\;\; s_{i}  &\textrm{(Eq.\ref{eq:lcoa}, \ref{eq:dfsi}, \ref{eq:dfsicon}, Fig.\ref{fig:ltree}bdf)}\notag\\
        L^{i=1,2,3,\cdots} &\textrm{ onto leftward subtrees } \;\; t_{0;\ge i} &\textrm{(Eq.\ref{eq:lsuar}, Fig.\ref{fig:ltree}ace)} \notag\\
      LU^{j=1,2,3,\cdots} &\textrm{ onto upward cotrees } \quad\quad T_{j} \textrm{, blue nodes generation } S_{j}  &\textrm{(Eqs.\ref{eq:ucoa},\ref{eq:dfSj},\ref{eq:dfSjcon}, Fig.\ref{fig:utree}bdf)}\notag \qed
\end{align}

Dehn's transformation formula implies that leftward cotrees and their node generations $s_i$ can be generated from the primal Cayley tree $T_{\ge0}$  both by a blue \textit{U}pward arrow followed by a red \textit{L}eftward path (Eq.\ref{eq:dfsi}) and by a red \textit{L}eftward path and a terminal conjugate blue \textit{U}pward arrow $U_i$ (Eq.\ref{eq:dfsiconjug}).
Upward cotrees and their node generations $S_j$ can be generated both 
by a red \textit{L}eftward arrow followed by a blue \textit{U}pward path (Eq.\ref{eq:dfSj}), and by a blue \textit{U}pward path and a terminal conjugate red \textit{L}eftward arrow $L_j$ (Eq.\ref{eq:dfSjconjug})  
\begin{align}
& L^i U_i && = L^i (L^{-i}UL^i)  && = UL^i &&\textrm{(both by Eqs.}\ref{eq:lcoa},\ref{eq:dfSj}\textrm{ and Eqs.}\ref{eq:lcoar},\ref{eq:dfSjconjug})\\
& U^j L_j && = U^j (U^{-j}LU^j)  && = LU^j &&\textrm{(both by Eqs.}\ref{eq:ucoa},\ref{eq:dfsi}\textrm{ and Eqs.}\ref{eq:ucoar},\ref{eq:dfsiconjug})
\end{align} 
\noindent For each unique cotree number pair, the unique entangled different colour, same generation (great-great-...-nephew $\leftrightarrow$ great-great-...-niece) number pair is traced by inverse functions acting on the two different colours. One function traces from blue, upward cotrees $T_i$ entangled number pairs in red, leftward cotrees $t_i$. The inverse function traces from from red, leftward cotrees $t_i$ entangled number pairs in blue, upward cotrees $t_i$. `Entanglement' comes from quantum entanglement experiments with two spin colours.
\begin{align}
&L_{-i} U^{-i}  L^i U_i:&&T_i\to t_i &&\textrm{(Figs.\ref{fig:utree}bdf} \;\;\to\;\; \textrm{Figs.\ref{fig:ltree}bdf}) \label{eq:L-iU-i}\\
&U_{-i} L^{-i}  U^i L_i:&&t_i\to T_i &&\textrm{(Figs.\ref{fig:ltree}bdf}\;\;\to\;\; \textrm{Figs.\ref{fig:utree}bdf)}.
\end{align}

\clearpage
\FloatBarrier

\clearpage
\FloatBarrier

\subsection {Eq. series 5: $a$n+b conjecture $a$ geometric series: cotree-density root-trunk test}
\label{sec:35eqs}

\begin{align}
&\textrm{Tree $T_{\ge0}$ provides an \textit{inscribed, non-isomorphic, periodic} congruence (residue, rest)  classes graph} \quad\quad\quad\quad\quad\quad\quad\quad\quad\quad\quad\quad & \textrm{ } \nonumber 
\\&\quad[T_{\ge0}]_{2a^2v},\;
\sm{\textrm{ of purple, red or blue, $L-$ or $U-$, congruence classes as \textit{nodes} and \;\;\;\quad\quad\quad}\\ 
\textrm{ of red and blue \textit{arrows} [\textcolor{red}{f},\textcolor{blue}{g}] generating purple, $L-$ or $U-$. successor classes }} \label{eq:gcc}  
\end{align}
\vspace{-0.6cm}

% \scriptsize per block works
% \begin{align}
% &\textrm{Collatz/Conway a\textit{n}+b function; }c\,=\,a\cdot1+b,\;\;g:\,n\to2n,\;\;g^{-1}:\,n\to n/2 &
% &&\textrm{Collatz }v=1
% &&\textrm{Bařina }v=2 \notag
% \\& \quad\quad\quad\quad\quad\quad a\in\; \textrm{Mersenne primes }& && a=2^2-1=3 && a=2^3-1=7 \notag
% \\& \quad\quad\quad\quad\quad\quad \scriptstyle{ f:\,n\to (n-1)/3} & && \scriptstyle{f:\,n\to (n-1)/3} && f:\,n\to (n\pm1)/7 \notag
% \\& \quad\quad\quad\quad\quad\quad f:\,n\to (n-b)/a& && f^{-1}:\;n\to 3n+1 && f^{-1}:\;n\to 7n\pm1 \nonumber
% \\& \quad\quad\quad\quad\quad\quad g:\,n\to 2n; g^{-1}:\,n\to n/2& && f^{-1}:\;n\to 3n+1 && f^{-1}:\;n\to 7n\pm1 \nonumber
% \end{align}

\begin{align}
&\sm{\textrm{Collatz/Conway a\textit{n}+b function; }c\,=\,a\cdot1+b,\;\;g:\,n\to2n,\;\;g^{-1}:\,n\to n/2\quad\quad\quad\quad\quad\;\;\;
\\\quad{\sm{a\quad \in\; \textrm{Mersenne primes }
3,7,31,127\quad\quad\quad\quad\quad\quad\quad\quad\quad;
\\\quad\sm{f\quad:\;\;n\;\to\; (n-b)/a\quad\quad\quad\quad\quad\quad\quad\quad\quad\quad\quad\quad\quad\quad;\;
\\f^{-1}\;\;:\;\;n\;\to\; a\,n+b\quad\quad\quad\quad\quad\quad\quad\quad\quad\quad\quad\quad\quad\quad\quad\; }}}} 
&
&&\sm{\textrm{Collatz }v=1
\\{\sm{a=2^2-1=3
\\\sm{f:\,n\to (n-1)/3\\f^{-1}:\;n\to 3n+1}}}}
&&\sm{\textrm{Bařina }v=2\\{\sm{a=2^3-1=7
\\\sm{f:\,n\to (n\pm1)/7\\f^{-1}:\;n\to 7n\pm1}}}} \nonumber
\end{align}

\begin{align}
&\quad\textrm{c}, \textrm{ cyclic root }c=a\cdot1+1
&
&&4=2^2
&&8=2^3
\\&\quad\textrm{v}, \textrm{ count subfunctions subclasses odd numbers $[1]_2\stackrel{\textup{bin}}{=}[1]_{10}$}
&
&&1
&&2
\\&\quad\quad\quad \textrm{provable $T_L$ tape conjectures: }v={\textup{log}_2(c)-1} 
\label{eq:vnum}
\\&\quad\textrm{min}_{q}, \textrm{ iterations to 1st blue/\textit{U} successor}\nonumber
\\&\quad\quad\quad \textrm{min}_{q}=\min\left(\,q\,\in\, g^q(n)=U(n)\,\right)\;=\;\frac{c-2}{\textup{log}_2(c)-1}
&
&&2
&&3 
\label{eq:minq}
\\&\quad\mathrm{max}_{pq}, \textrm{ iterations to purple red/\textit{L} or blue/\textit{U} successor}\nonumber
\\&\quad\quad\quad\mathrm{max}_{pq}=\max\left(\sm{p\in\; fg^p(n)=L(n)\\\;\;q\in\;g^q(n)=U(n)\;\;}\right)
& &&{4} &&{7}
\label{eq:maxpq}
\\&\quad\mathrm{lcm}_{o}, \textrm{ perfect periodicity odd numbers indivisible } \nmid a \nonumber
\\&\quad\quad\quad \mathrm{lcm}_{o}\;\;\,=\;2^v a^1\;=\;c/2\,(c-1) 
& &&\sm{2^13^1=\\6} &&\sm{2^2 7^1=\\28}
\label{lcmo}
\\&\quad\mathrm{lcm}_{pq}, \textrm{ period purple, red/\textit{L} or blue/\textit{U}, $n$\,=\,a\,n$_\textup{odd}$+b numbers} \nonumber
\\&\quad\quad\quad \mathrm{lcm}_{pq}\;=\;2^v a^2\;=\;c/2\;(c-1)^2 
& &&\sm{2^1 3=\\18} &&\sm{2^2 7^2=\\196}
\label{eq:lcmpq}
\\&\quad n_o,\, (c-2)(\log_2(c)-1) \textrm{ odd classes }a\nmid[c_o]_{\textup{lcm}_o=c/2(c-1)} & 
&&\sm{(3-1)\cdot1=2\; c_o\\\left[1,5\right]_{6}}  
&&\sm{12\;c_o\;\;
\\\begingroup \setlength\arraycolsep{1pt}
\scriptsize{
\begin{bmatrix} 1 & \,5 & \,9 
\\              13 & \,17 & \,25 
\\              3 & \,11 & \,15 
\\              19 & \,23 & \,27 
\end{bmatrix}_{2^27=28}} \endgroup }
\label{eq:oclass}
\\&\quad n_b = n_o \textrm{ even $c_b\,=\,a c_o + b$ branching classes $[c_b]_{\mathrm{lcm}_{pq}=c/2(c-1)^2} $} 
& 
&&\sm{(3-1)\cdot1=2\;c_b\\\left[4,16\right]_{2^13^2=18}}  
&&\sm{(7-1)\cdot2=12 c_b\;\;
\\\begingroup \setlength\arraycolsep{1pt}
\scriptsize{
\begin{bmatrix} 8 & 36 & 64 
\\              92 & 120 & 176 
\\              20 & 76 & 104 
\\              132 & 160 & 188 
\end{bmatrix}_{\sm{2^27^2\\=196}}}\endgroup}
\label{eq:bclass}
\\&\quad d(S_{\ge0})\textrm{, $d$ensity purple, red or blue, classes }c_b=ac_o+b  
& && && \nonumber
\\&\quad\quad\quad d(S_{\ge0})\;=\;\frac{n_b}{\textrm{lcm}_{pq}}\;=\;\frac{(c-2)\;(\log_2(c)-1)}{c/2\;(c-1)^2}  
& &&\frac{2}{18}&&\frac{12}{196} \label{eq:dSge0}
\\&\quad\mathrm{lcm}_{LU} \textrm{, period red/\textit{L} \& blue/\textit{U}  successors} \nonumber
\\&\quad\quad\quad\mathrm{lcm}_{LU}=\mathrm{lcm}_{pq}\cdot2^{\mathrm{max}_{pq}}=a^2\,2^{v+\textrm{max}_{pq}}
& && &&
\nonumber
\\&\quad\quad\quad\quad\quad\;\;=
\log_2(c-1)\frac{c^3}{2}(c-1)^2
& &&\sm{\sm{18\cdot2^4\;=\\288\;\;=\\\textrm{2 gross }}}&&\sm{196\cdot2^7\;=\\25088\;\;\;\,}
\label{eq:lcmLU}
\\&\quad n_{bs}, (a^2-a)v \textrm{ subcl. $[c_{bs}]_{\mathrm{lcm}_{LU}=\mathrm{lcm}_{pq}\cdot a=2^va^3}$}\,; f^{-1}:c_{bs}\to c_o
&
&&\sm{(9-3)\cdot1=6\; c_{bs}
\\\begingroup \setlength\arraycolsep{1pt}
\scriptsize{
\begin{bmatrix} 4 & 22 & 40 
\\              16 & 34 & 52 
\end{bmatrix}_{\sm{18a\\=54}}}\endgroup}  
&&\sm{(49-7)\cdot2=84\;c_{bs}
\\\begingroup \setlength\arraycolsep{1pt}
\tiny{
\begin{bmatrix} 8 & 36 & \cdots 
\\              \cdots & 1296 & 1352
\\              20 & 76 & \cdots
\\              \cdots & 1336 & 1364
\end{bmatrix}_{\sm{196a\\1372}}}\endgroup}  
\label{eq:bsclass} 
\end{align}

\begin{align}
&
&
&&\sm{\textrm{Collatz }v=1
\\{\sm{a=2^2-1=3
\\\sm{f:\,n\to (n-1)/3\\f^{-1}:\;n\to 3n+1}}}}
&&\sm{\textrm{Bařina }v=2\\{\sm{a=2^3-1=7
\\\sm{f:\,n\to (n\pm1)/7\\f^{-1}:\;n\to 7n\pm1}}}} \nonumber
\\&U^1(T_{\ge0})=\left[
T_{\ge0}\right]_{\textrm{lcm}_{LU}}\cdot\; \textup{\textbf{R}}\;=\;
\left[T_{\ge1}\right]_{\textrm{lcm}_{LU}\cdot2^{\textrm{lcm}_{pq}}}\;\;
&
&& \scriptsize{
\begin{aligned}
& T_{\ge0}\begin{bmatrix} 4 \\ 1 \end{bmatrix}_{288} \;\cdot\;
\begingroup \setlength\arraycolsep{1pt}\begin{bmatrix} 0 & 1 \\ 1 & 0 \end{bmatrix} =\endgroup \\
&T_{\ge1}\begin{bmatrix} 1\\4\end{bmatrix}_{288\cdot2^4}
= T_{\ge0}\begin{matrix} \begin{bmatrix}1\\4\\1\\4\\1\\4\end{bmatrix}_{288a\cdot2^4} \end{matrix} 
\end{aligned}}
&&{\footnotesize{\sm{ \textrm{Section }\\\ref{sec:ccba} }}}
\label{eq:Tge1} 
\\&\;\; {U}^j \textrm{ vector }T_{\ge j}\left[\sm{n_{bL}\\n_{bU}}\right]_{\textup{dim}=(n_b)}^{j}; \textrm{\textit{ U} matrix }  [\textbf{R}]_{(n_b\times n_b)}& && && \notag
\\&L(T_{\ge0})=
\left[T_{\ge0}\right]_{a\cdot\textrm{lcm}_{LU}}\cdot\; \textup{\textbf{r}}\;=\;
\left[t_{0;\ge1}\right]_{\textrm{lcm}_{LU}a\cdot2^{2 min_q}}\quad
&
&& \scriptsize{
\begin{aligned}
&\;\; T_{\ge0}\begin{matrix} \begin{bmatrix}4\\1\\4\\2\\8\\8\end{bmatrix}_{288} \end{matrix} \;\cdot\;
\begingroup \setlength\arraycolsep{1pt}
\begin{bmatrix} 1 & 1 & 1 & 0 & 0 & 0 
\\              1 & 1 & 1 & 0 & 0 & 0
\\              0 & 0 & 0 & 1 & 1 & 1
\\              1 & 1 & 1 & 0 & 0 & 0 
\\              1 & 1 & 1 & 0 & 0 & 0
\\              0 & 0 & 0 & 1 & 1 & 1
\end{bmatrix} =\endgroup \\
&t_{0;\ge1}
\begin{matrix} \begin{bmatrix}15\\15\\15\\12\\12\\12\end{bmatrix}_{288a\cdot2^4}
\end{matrix}
=t_{0;\ge1}\begin{bmatrix}15\\12\end{bmatrix}_{288\cdot2^4}
\end{aligned}
}
&&{\footnotesize{\sm{\textrm{Section}\\\ref{sec:ccba}} }}
\label{eq:tge0} \\
&\;\;{L}^j \textrm{ vector }t_{0;\ge i}\left[\sm{n_{bsL}\\n_{bsU}}\right]_{\textup{dim}=(n_{bs})}^{i};\; \textrm{\textit{ L} matrix }  [\textbf{r}]_{(n_{bs}\times n_{bs})}& && && \notag
\end{align}

\begin{align}
&n_L=v^2a^3\textrm{ red \textit{L}-classes in Cayley a\textit{n}+b tree }T_{\ge0} & 
&& 3^3=27 
&& 2^27^3=1372
\label{eq:sgener1}
\\& A^+=(a(a+1)+1),\; A^+v \textrm{ leftward successors multiplier }  & 
&&13\cdot1=13 
&&57\cdot2=114
\label{eq:Apm}
\\&n^{i;\text{ for }i\ge1}_s\;=\;v\,A^{+i}\,s_{i-1}\textrm{ generation $s_i$ classes in cotree }t_i & 
&& n^i_s=13^i \,n^{i-1}_s 
&& n^i_s=114^i\,n^{i-1}_s
\label{eq:sgeneri}
\\&n_U=v^2A^-\textrm{ blue \textit{U}-classes in a\textit{n}+b tree }T_{\ge0},\, A^-=(a(a-1)-1)  & 
&& 1^2\cdot5=5 
&& 2^2\cdot41=164 \label{eq:Am}
\\&T_{\ge0}[\sm{n_L\\n_U}]_{\text{lcm}_{LU}}=
  \bigg[{{\sm{n_L\;\,=\;a^3v^2 \quad\textrm{ red }c_L\textrm{ classes }\\
             n_U=\,A^-v^2 \;\; \textrm{ blue }c_U \textrm{ classes }}}}\bigg]_{\mathrm{lcm}_{pq}\,\cdot\,2^{\mathrm{max}_{pq}}} 
&  && \bigl[{\sm{27\; c_L\\5\;c_U}}\bigl]_{288} 
&& \bigl[{\sm{1372\;c_L\\164\;c_U}}\bigl]_{25088} 
\label{eq:luclass}
\\&UL^i=L^iU_i:\;T_{\ge0}\left[\sm{n_L\\n_U}\right]_{\mathrm{lcm}_{LU}} \textcolor{red}{\;\to}
t_i\left[n_s^i=vA^{+i}n_s^{i-1}\right]_{\mathrm{lcm}_{LU}\;\cdot\;2^{\,2i\; \mathrm{min}_{q}}} & 
    && \left[\sm{n_s^i=\\13^in_s^{i-1}}\right]_{\sm{288\cdot\\2^{4i}}}
    && \left[\sm{n_s^i=\\114^in_s^{i-1}}\right]_{\sm{25088\cdot\\2^{6i}}}
\label{eq:lcoicm1}
\\&LU^j=U^jL_j:\;T_{\ge0}\left[\sm{n_L\\n_U}\right]_{\mathrm{lcm}_{LU}}\textcolor{blue}{\;\to\;}T_j\left[n_{S_{j=\ell(a-1)+\alpha}}\right]_{\mathrm{lcm}_{LU}\,\cdot\;2^{j\,\mathrm{max}_{pq}}}&
    &&  \sm{n\;_S^j=n_U \cdot\\ 3 \cdot 13^{k-1}} && \sm{\textrm{Section}
\label{eq:ucojcm}}
\\&\quad n_{S_j}\textrm{ generation $S_j$ classes,}\;\;\ell=0,1,2,\cdots,\;\;\alpha=1,\cdots,a-1 \nonumber
\\&\mathrm{per}_{T}=a\cdot \textrm{min}_q,\;\; 2^{\textrm{per}_{T}} \textrm{ periodicity by ${+(a-1)}$ iterations}  & 
\notag 
\\&\quad\ U^{a-1}:T_j(c_r\in S_j)\textcolor{blue}{\to} T_{j+(a-1)}(c_r\cdot2^{\textrm{per}T}\in S_{j+(a-1)}) &
&&\sm{c_r\,\cdot\,2^{3\cdot2=6} \textrm{ by}\\{+(3-1)=+2}}
&&\sm{c_r\,\cdot\,2^{7\cdot3=21}\textrm{ by}\\{+(7-1)=+6}}\label{eq:perT}
\\&\mathrm{per}_{tT}=\mathrm{lcm}({\mathrm{max}_{pq},\mathrm{per}_{T}}) \textrm{, periodicity by $+\frac{\textrm{per}_{tT}}{\textrm{per}_{T}}(a-1)$ iter.} & 
&&\sm{c_{r} \,\cdot\,2^{\mathrm{lcm}(4,6)=12}\\+\frac{12}{6}(3-1)=+4}
&&\sm{c_{r} \,\cdot\,2^{\mathrm{lcm}(7,21)=21}\\+\frac{21}{21}(7-1)=+6}\label{eq:pertT}
\end{align}

\begin{align}
& &
&&\sm{\textrm{Collatz}\\{\sm{a=2^2-1=3\\f:n\to (n-1)/3}}}
&&\sm{\textrm{Bařina}\\{\sm{a=2^3-1=7\\f:n\to (n\pm1)/7}}} \nonumber
\\&\quad d(tT), \textrm{ periodic red/\textit{t} + blue/\textit{T} cotree density} \notag
\\&\quad\quad\textrm{4 iterators:}\notag
\\&\quad\quad\quad i = 1,2,3,\dots\; \;\;\quad\textrm{geometric series \textit{t-}cotree densities}\notag
\\&\quad\quad\quad \ell=0,1,2,\dots\;\;\, \quad\textrm{period arithmetic progression}\notag
\\&\quad\quad\quad\alpha=1,\dots,a-1\;\;\;\; a-1\textrm{ geometric series \textit{T}-cotrees} \notag  
\\&\quad\quad\quad
j = \ell(a - 1) + \alpha\;\textrm{ geometric series \textit{T}-cotree densities} \notag
\\&\quad\quad\textrm{specifying: } 1+(a-1)=a\textrm{ geometric series of }\notag
\\&\quad\quad\quad\textrm{ infinite arithmetic cotree density progressions}\notag
\\&\quad\quad d(tT)=
\sum\limits_{i=1}^{\infty} d(s_i)+
\sum_{\alpha=1}^{a-1}\;\sum_{\ell=0}^{\infty}  d(S_{j=\ell(a-1)+\alpha}) \notag
\\&\quad\quad  \textcolor{white}{d(tT)}=
\sum\limits_{i=1}^{\infty} 
\left[n_s^i=vA^+n_s^{i-1}\right]_{\mathrm{lcm}_{LU}\;\cdot\;a^{\,2i\; \mathrm{min}_{q}}}
+  & && \frac{27}{288} + && \frac{1372}{25088}+\notag
\\&\quad\quad  \textcolor{white}{d(tT)=}
\sum_{\alpha=1}^{a-1}\;\sum_{\ell=0}^{\infty} d[S_j]_{\textrm{lcm}_{LU}\cdot 2^{\alpha\min_q}} & && \frac{1}{2^6-1} + &&\frac{12288+\cdots+1}{2^{21}-1}+
\notag
\\&\quad\quad  \textcolor{white}{d(tT)=}
& && \frac{3}{288\cdot7} &&\frac{28}{25088(2^7-1)377}
\notag
%cotree density proof%
\\&\quad\quad  \textcolor{white}{d(tT)}=\frac{d(t_1)}{1-d(t_2)/d(t_1)}+ 
\notag
\\&\quad\quad  \textcolor{white}{d(tT)=}
\sum_{\alpha=1}^{a-1}\;\frac{d(T_\alpha)}{1-d(T_{\alpha+a})/d(T_\alpha)}\quad  \textrm{ (`cotree density')} & &&\frac{27+5}{288}&&\frac{1372+164}{25088}  \label{eq:dtT}
%subtree density proof%
\\&\quad\quad  \textcolor{white}{d(tT)}=
\frac{\;n_L+n_U\;}{\textrm{lcm}_{LU}}=
\frac{\;a^3v^2+A^+v^2\;}{a^2\,v\,U(c)} \;
\textrm{ (`subtree density')} & &&\frac{27+5}{288}&&\frac{1372+164}{25088}  \label{eq:dsub}
\\&\quad\Rightarrow\textrm{Theorem \ref{the:decidab}:} \textrm{ a\textit{n}+b-conjecture proven}\stackrel{?}{-}\textrm{refuted by:}\notag
\\&\quad\quad\;\;\textrm{\textbullet\;`cotree density proof' / `subtree density proof'} 
 \notag
\\&\quad\quad\quad\;\;\; d(tT)\stackrel{?}{=} d(S_{\ge0})\textrm{, Hamilton tour }S_{\ge0}=\bigcup_{n=1}^{\infty}\textrm{a}n+\textrm{b?}
& &&\stackrel{?}{=} \frac{2}{18}&&\stackrel{?}{=}\frac{12}{196}\label{eq:cdp}
\\&\quad\quad\;\;\textrm{\textbullet\;`root successor test' }T_{\ge 0}\stackrel{?}{=}G_{\ge 0}\textrm{ on period }\ell=0: \notag
\\&\quad\quad\quad\;\;\;\, a-1\textrm{ even }\ell=0\textrm{ numbers } n^{\textrm{even }}_{\ell=0}&
&& 4 \xrightarrow{\text{?}} \;\;\sm{22,\\40\,}\;
&& 8 \xrightarrow{\text{?}} \;\;\sm{204,400,596,\\792,988,1184} \notag
\\&\quad\quad\quad\;\;\;\, a-1 \textrm{ odd }\ell=0\textrm{ numbers } f:n^{\textrm{even }}_{\ell=0}\to (n\pm b)/a&
&& 1 \xrightarrow{\text{?}} \;\;\sm{\,7,\\13\,}\;
&& 1 \xrightarrow{\text{?}} \;\;\sm{29,57,85,\\113,141,169} 
\label{eq:rst}
\\&\quad\quad\;\;\textrm{\textbullet\;{`root-trunk test' }}T_{\ge 0}\stackrel{?}{=}G_{\ge 0} \textrm{ on period }\ell=0\notag
\\&\quad\quad\quad\;\;\,\textrm{ number theoretic: }c \sm{?\\\textrm{ power-of-2}} \textrm{,}\; U(c)\; \stackrel{?}{\in}\; (\ell=0)
& &&\sm{c=3\cdot1+1\;\stackrel{!}{=}\;2^2\\U(c)=2^4\stackrel{!}{=}0+2^4}
&& \sm{c=7\cdot1+1\;\stackrel{!}{=}\;2^3\\U(c)=2^6\;\stackrel{!}{=}\;0+2^6} \notag 
\\&\quad\quad\quad\;\;\,\textrm{ graph theoretic: }\sm{\textrm{one trunk }\\\textrm{not a shrub ?}}\textrm{,} \sm{\textrm{lowest trunk node}\\  U(c)\;\stackrel{?}{\leq}\;(\textup{lcm}_{pq}=(c-1)^2)\frac{c}{2}}
& &&\sm{\textrm{Fig.\ref{fig:tge0}: 1 trunk?}
\\ \sm{U(c)=2^4\;\stackrel{!}{\leq}3^22}\\16\;\stackrel{!}{\leq}\;18}
&&\sm{\textrm{Fig.\ref{fig:tree7pm1}: 1 trunk?}
\\ \sm{U(c)=2^6\;\stackrel{!}{\leq}\;7^22^2\\64\;\stackrel{!}{\leq}\;196}} 
\label{eq:rtt}
\end{align}

\noindent Subsection \ref{sec:35eqs} rests on the observation by Emmy Noether \cite{RN18} that a `module' consists of congruence classes with a Least Common Multiple periodicity (LCMp), having isomorphic `submodules' of congruence classes whose LCMps are multiples of its LCMp. The congruence class graph with purple, either red or blue, $a\nmid n$=n$_{\textup{odd}}$+b numbers in an $an+b$ Cayley tree `module' $T_{\ge0}\subseteq{G_{\ge0}}$ (Figs.\ref{fig:cc3m1},\ref{fig:cc3p1},\ref{fig:cc7pm1}), shows $g$-path lengths in leftward $L=fg^p$ paths defining the $L$eftward function (Defs.\ref{def:Lw},\ref{def:Lw7}) and in upward $U=g^q$ paths defining the $U$pward function (Defs.\ref{def:Uw},\ref{def:Uw7}).  The maximum $g$-path length in leftward $L=fg^p$ paths and upward $U=g^q$ paths shows the Least Common multiple periodicity $\textrm{lcm}_{\textup{pq}}=\frac{c}{2}(c-1)^2=2^va^2$ (Eq.\ref{eq:lcmpq}) of purple, either red or blue $a\nmid$ an$_{\textup{odd}}$+b numbers (Eq:\ref{eq:bclass}): $2^13^2=18$ for Collatz $3n+1$ and $2^27^2=196$ for Bařina $7n\pm1$. The number of either blue or red classes amounts to $n_b=(c-2)(\log_2(c)-1)=(a-1)v$:  $(3-1)1=2$ for Collatz $3n+1$  and $(7-1)2=12$ for Bařina $7n\pm1$. The periodic density of red or blue numbers amounts to $\frac{n_o}{\textrm{lcm}_{\textup{pq}}}$: $\frac{2}{18}$ for Collatz $3n+1$  and $\frac{12}{196}$ for Bařina $7n\pm1$ (Eq.\ref{eq:dSge0}). 

To fit on a single page, the congruence class graph $G_{CC}$ for Bařina $7n\pm1$ is not depicted with its periodicity of $2^va^2=196$ (Eq.\ref{eq:lcmpq}), but with the `perfect' periodicity of yellow odd numbers $(2^va^2)/a=28$,  showing $v=2$ cyclic $g$-paths of $g$-path lengths $2$ or $2+(v+1)=5$, $3$ or $3+(v+1)=6$, or, $4$ or $4+(v+1)=7$ (Fig.\ref{fig:cc7pm1}). Depicting the congruence class graphs for Collatz 3$n$+1 (Fig.\ref{fig:cc3p1}) and of 3$n$-1 (Fig.\ref{fig:cc3m1}) similar to Fig.\ref{fig:cc7pm1} would have shown $v=1$ cyclic $g$-path of $g$-path lengths $1$ or $1+(v+1)=3$, or, $2$ or $2+(v+1)=4$.

The light-purple subclasses of purple classes at the start and arrowhead of $L=fg^p$ paths via same-periodicity $\textrm{lcm}_{\textup{pq}}=2^va^3=\frac{c}{2}(c-1)^2$
yellow odd classes (Eq.\ref{eq:oclass}) have a multiple \textrm{LCMp} of $\textrm{lcm}_{\textup{pq}}=2^va^3=\frac{c}{2}(c-1)^3$. The number of light-purple subclasses amounts to $n_{bs}=(c-2)(c-1)(\log_2(c)-1)=(a-1)av$ (Eq.\ref{eq:bsclass}):  $(3-1)3=6$ for Collatz $3n+1$ and $(7-1)2=12$ for Bařina $7n\pm1$. The $\textrm{lcm}$-periodicity of red $L$eftward and blue $U$pward $a\nmid$ an$_{\textup{odd}}$+b numbers amounts to
$\textup{lcm}_{\textup{LU}}=\log_2(c-1)\frac{c^3}{2}(c-1)^2=$ (Eq.\ref{eq:lcmLU}): $288=2$ gross for Collatz $3n+1$ and $25088$ for Bařina $7n\pm1$.

Each arrowhead node with a same numbered site above it has one red $L$eftward child and hence a unique infinite series of unique $U$ward successors, one in each generation $S_{j=1,2,3,\dots}$ in blue $U$pward cotrees $T_j$  to the left of its site (Figs.\ref{fig:utree}bdf) by paired depth-first functions $\sm{01^j\\LU^j}:$   (Eqs.\ref{eq:dfSj},\ref{eq:ucoa},\ref{eq:ucojcm}).   
Each arrowhead node with a same numbered site above it has one blue $U$pward child and hence a unique infinite series of unique $L$eftward successors, one in each generation $s_{i=1,2,3,\dots}$ in red $L$eftward cotrees $t_i$  to the right of its site (Figs.\ref{fig:ltree}bdf) by paired depth-first functions $\sm{10^i\\UL^i}:$   (Eqs.\ref{eq:dfsi},\ref{eq:lcoa},\ref{eq:lcoicm1}).  

$U$pward subtrees $T_{\ge j}$ (Eq.\ref{eq:usuar},Fig.\ref{fig:utree}ace) and $L$eftward subtrees $t_{0,\ge i}$ (Eq.\ref{eq:lsuar},Fig.\ref{fig:ltree}ace) include number pairs in all future generations $S_{\ge j}$ respectively $s_{\ge i}$, giving non-converging geometric series sums of subtree densities. Cotrees of same-generation successors (Fig.\ref{fig:ltree}bdf,\ref{fig:utree}bdf, in geneaology terms of nephews and nieces, great-nephews and great-nieces, great-great-nephews and great-great-nieces share the least common multiple periodicity $\textrm{lcm}_{pq}$ of corresponding subtrees (Fig.\ref{fig:ltree}ace,\ref{fig:utree}ace) of all-future-generation-successors, but have converging geometric cotree density series sums \textit{d}ensity $0<d<1$ (Eq.\ref{eq:dtT}). 

The cumulative cotree density by converging geometric series sums of red $L$eftward and blue $U$pward cotrees $d(tT)$ is given by four different iterators (Eq.\ref{eq:dtT}) depending on parameters $A^+=a^2+a+1$ (Eq.\ref{eq:Apm}) and $A^-=a^2-a-1$ (Eq.\ref{eq:Am}). Red $L$eftward cotrees $t_{i>1}$ comprise a multiple of $A^+v^2=(a^2+a+1)v^2$ $s_i$-congruence classes of the number of $s_{i-1}$ congruence classes (Eq.\ref{eq:Apm}): $13$ for Collatz 3$n$+1, $114$ for Bařina $7n\pm1$. They are generated from a number of $a^3v^2$ $L$eftward congruence classes in Cayley tree $T_{\ge0}$ itself: $3^31^2=27$ for Collatz 3$n$+1, $2^73^2=1372$ for Bařina $7n\pm1$. All blue $U$pward cotrees $T_j$ have the same number of $A^-v^2=(a^2-a-1)v^2$ $U$pward congruence classes in $U$pward congruence class generations $S_j$ (Eq.\ref{eq:Am}): $5$ for Collatz 3$n$+1, $164$ for Bařina $7n\pm1$.

The Cayley tree $T_{\ge0}\subseteq{G_{\ge0}}$ (Theorem \ref{the:caytree}) is a biperiodic number graph `module', generating through its $L$eftward depth-first iterator $i=1,2,3,\dots$ and upward depth-first iterator $j=1,2,3,\dots$ isomorphic:
\begin{itemize}
    \item red $L$eftward subtrees $t_{0;\ge i}$ (Fig.\ref{fig:ltree}ace, Eq.\ref{eq:lsuar}) and cotrees $t_i$ (Fig.\ref{fig:ltree}bdf, Eq.\ref{eq:lcoar}) with $L$eftward $s$-generations $s_i$-nodes (Eq.\ref{eq:dfsi}), having 
    \begin{itemize}
        \item a multiple LCMp $m_i$ of $s_i$ generation congruence classes of the LCMp of the previous generation $s_{i-1}$ congruence classes, initiated from the LCMp $2^va^2$ of tree $T_{\ge0}$ (Eq.\ref{eq:lcoicm1}), 
    \item a multiple of $A^+v^2=(a^2+a+1)v^2$ $s_i$-congruence classes of the number $s_{i-1}$ congruence classes, initiated from a number of $a^3v^2$ $L$eftward congruence classes in $T_{\ge0}$    
    \end{itemize}
    \item blue $U$pward subtrees $T_{\ge j}$ (Fig.\ref{fig:utree}ace, Eq.\ref{eq:usuar}) and cotrees $T_j$ (Fig.\ref{fig:utree}bdf, Eq.\ref{eq:ucoar}) with $U$pward $S$-generations $S_j$-nodes (Eq.\ref{eq:dfSj}), having by \ref{eq:ucojcm}) 
        \begin{itemize}
        \item a multiple LCMp $m_j$ of $S_j$ generation congruence classes of the LCMp of the previous generation $S_{j-1}$ congruence classes, initiated from a number of $A^-v^2=(a^2-a-1)v^2$ $U$pward congruence classes   (Eq.\ref{eq:ucojcm}), 
    \item the same number of $A^-v^2=(a^2-a-1)v^2$ $U$pward congruence classes in $U$pward congruence class generations $S_j$ of nodes in $T_j$-cotrees as in tree $T_{\ge0}$, albeit with a multiple LCMperiodicity.   
    \end{itemize}
    \end{itemize}

Equations \ref{eq:minq}-\ref{eq:lcmLU},\ref{eq:bsclass}-\ref{eq:rtt} recur in the cotree density test (Eq. \ref{eq:cdp}) and root-successor test (Eq.\ref{eq:rtt}), yielding Lean-certified proven Collatz 3$n$+1 and Bařina 7$n\pm$1 conjectures (Tables \ref{tab:cdp3p1} and \ref{tab:cdp7pm1}), and refuted $3n-1$ and $5n+1$ conjectures (Theorem \ref{the:decidab}).

The `cotree density' proof or refutation of an a$n$+b conjecture rests on $a$ geometric series of successive isomorphic cotrees, where the periodicities $m_i$ and $m_j$ of the disjoint congruence classes in successive $L$eftward and $U$pward cotrees are scaled up to least common multiple periodicities of higher generation cotrees (Eq.\ref{eq:cdp}). The $a+1$ iterators in the Cotree density formula (Eq.\ref{eq:dtT}) reflect the automorphism graph Aut($T_{\ge0},[L,U]$) with 1 iterator counting congruence classes in $L$eftward cotrees, $a-1$ series $0,\dots,a-1$ of different upward cotrees, and 1 iterator counting periods $\ell=0,1,2,\dots$ in congruence classes / arithmetic progressions $\ell m+r$. Since each congruence class $[r]_m$ is an infinite arithmetic progression $\ell m+r$ for $\ell=0,1,2,\dots$ and $r=0,\dots,m-1$ a proposition proven or refuted for a class $[r]_m$ (by the cotree density test, Eq.\ref{eq:cdp}), is also proven respectively refuted by period $\ell=0$ number $r$ (by the root-trunk-test, Eq. \ref{eq:rtt}, bottom row Table \ref{fig:TLEE}).

The root successor test (Theorem \ref{the:decidab}, Eq.\ref{eq:rst}) assesses whether numbers in the lowest period $\ell=0$ from $T_{\ge0}$ congruence classes are taken on board in the Hamiltonian tour, $7$ and $13$ in the case of the proven Collatz $3n+1$ conjecture and the refuted $3n-1$ conjecture (Theorem \ref{the:decidab}). 
The root-trunk test assesses whether the lowest trunk node is a period $\ell=0$ number $(U(c)=c^2\in\{\ell=0\})<2a^2v$ taken on board. This is proven for Collatz $3n+1$ (Figs.\ref{fig:tge0},\ref{fig:cc3p1}) and Bařina $7n\pm1$ (Figs.\ref{fig:tree7pm1},\ref{fig:cc7pm1}) but refuted for a shrub rather than a 1-trunk tree in which $c=a\cdot1+b$ not a power of $2$ (Fig.\ref{fig:tree5p1} $5n+1$) and/or a lowest trunk node $(U(c)>c^2)\in\{\ell\ge0\} > 2a^2v$ (Fig.\ref{fig:tree3m1} $3n-1$, Eq.\ref{eq:rtt}). Lowest trunk nodes $c^2=4^2=16<18$ for Collatz $3n+1$ and  $c^2=8^2=64<196$ for Bařina $7n\pm1$, included in the title, provide Lean-certified conjecture proofs (Tables \ref{tab:cdp3p1},\ref{tab:cdp7pm1}, \ref{fig:Gccc},\ref{fig:TLEE}).

The proven Collatz $3n+1$ function yields an LCM$_{LU}$ (Eq.\ref{eq:lcmLU}) of its red $s$- and blue $S$-generations of a dozen days of two dozen hours, $12 \cdot 2 \cdot 12 = 288$ (Eq.\ref{eq:lcmpq}), being the $3\cdot2^4$ multiple of the perfect number $6$ \cite{klelkies2019}, thereby reflecting our Sumerian-Babylonian $60=12\cdot5$-base clocks of minutes and degrees, daily hours, and yearly months. The Bařina $7n\pm1$ function with an LCM$_{LU}$ of $2^9 7^2 = 25088$ (Eq.\ref{eq:lcmpq}), being the $7\cdot2^7$ multiple of the perfect number $28$ \cite{klelkies2019} of a lunar period of approximately $4$ weeks of $7$ days, reflecting our Jewish-Christian-Islamic weekly calendar enabling short-term planning by day- and moonlight.

\subsubsection{Collatz $3n+1$: Hamiltonian tour cotree density proof and root-trunk-proof}
\label{sec:cc}
\begin{figure}[hbt!]
\caption{Collatz 3\textit{n}+1 congruence classes graph 
${G}_{cc}=({\left[0,\cdots,17\right]}_{18},[f,g])$ }
\label{fig:cc3p1}
\vspace{-0.4cm}
\begin{center}
\includegraphics[width=0.67\textwidth]{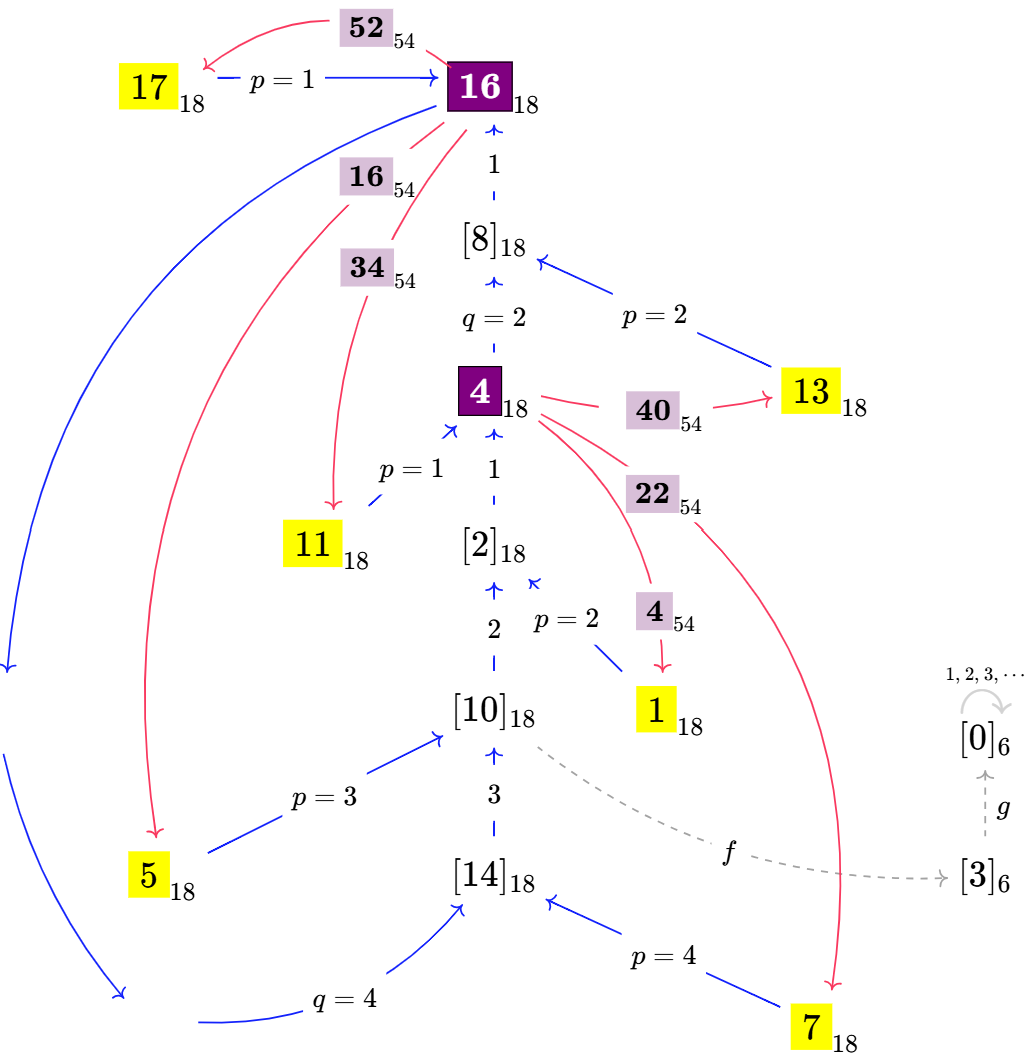}
\end{center}
\vspace{-0.2cm}
\centering\begin{minipage}{0.9\linewidth}
\noindent\textbf{Legend: } Nodes: $a-1=2$ even, purple, branching classes with periodicity/modulus $m=2a^2v=18$ (Eq.\ref{eq:lcmpq}), e.g. $\colorbox{violet}{\textcolor{white}{16}}_{18}$, and $\,a(a-1)=6\,$ odd, yellow, classes, e.g.
$\colorbox{yellow}{\textcolor{black}{5}}_{18}$, (Eq.\ref{eq:oclass}) 
from $a(a-1)=6$ even, lightpurple, branching parent subclasses with periodicity $2a^3v=54$ (Eq.\ref{eq:bsclass}).
Red arrow $\textcolor{red}{\xrightarrow{\textcolor{black}{f:\;n\,\to\,(n-1)/3}}}$ (Eqs.\ref{eq:rootPairs},\ref{eq:oclass}) to odd yellow child annotated with its lightpurple (here lightgray) parent subclass, e.g.
$\colorbox{violet}{\textcolor{white}{16}}_{18}
\textcolor{red}{\xrightarrow{\colorbox{lightgray}{\textcolor{black}{[34]}}}}
\colorbox{yellow}{\textcolor{black}{11}}_{18}$.
Blue arrow $\textcolor{blue}{\xrightarrow{\textcolor{black}{g:\;n\,\to\,2n}}}$, either in blue yellow-to-purple path with path length  $\textcolor{blue}{\xrightarrow{\textcolor{black}{p=\cdots, \textrm{ defining } L=fg^p}}}$, e.g. 
$\colorbox{yellow}{\textcolor{black}{5}}_{18}\textcolor{blue}{\xrightarrow{\textcolor{black}{p=3}}}\colorbox{violet}{\textcolor{white}{4}}_{18}$, 
or in blue purple-to-purple path with path length  
$\textcolor{blue}{\xrightarrow{\textcolor{black}{q=\cdots, 
\textrm{ defining } U=g^q}}}$, e.g.  
$\colorbox{violet}{\textcolor{white}{4}}_{18}
\textcolor{blue}{\xrightarrow{\textcolor{black}{q=2}}}
\colorbox{violet}{\textcolor{white}{16}}_{18}$. 
The $3n+1$ tree proves the Collatz $3n+1$ conjecture by the root-trunk-test (Eq.\ref{eq:rtt}), as it is a one-trunk-tree with a power-of-2 root $c=3n+1=2^2$, and with a  trunk in period $\ell=0$: $(U(c)=2^{2+(q=2)}=16)=((\ell=0)\cdot18+16)\;<\;(m=18)$.$\quad\textbf{\textbar}$
\end{minipage}\end{figure}

The congruence class (residue class, rest class) equations using Least Common multiple periodicities ($\textrm{lcmp}$)  from Section \ref{sec:35eqs} (Eqs.\ref{eq:gcc}-\ref{eq:rtt}) are specified in this section for the proven Collatz $3n+1$ (Figs.\ref{fig:tge0},\ref{fig:cc3p1}) and Bařina $7n\pm1$ (Figs.\ref{fig:tree7pm1},\ref{fig:cc7pm1}) conjectures and in between for the refuted $3n-1$ conjecture (Figs.\ref{fig:tree3m1},\ref{fig:cc3m1}). 

\noindent$L$eftward paths starting with a red arrow generate, by e.g.
$L=fg^p:\colorbox{magenta}{\textcolor{white}{16}}_{54}\textcolor{red}{\to}
[[16]_{54}-1)/3\cdot 2^{p=3}=\colorbox{pink}{\textcolor{black}{88}}]_{288}$, $a^3=27$ 
red $L$eftward Cayley graph classes (Eqs.\ref{eq:lcmLU},\ref{eq:luclass}), e.g. 
$\colorbox{pink}{\textcolor{black}{88}}_{288}$, of branching red $L$eftward Collatz $3n+1$ nodes (Fig.\ref{fig:tge0}a).
$U$pward paths without a red arrow generate, e.g. by
$U=g^q:\colorbox{violet}{\textcolor{white}{4}}_{18}
\to[[4]_{18}\cdot 2^{q=2}=\colorbox{cyan}{\textcolor{black}{16}}]_{288}$,
$a^2-a-1=5$ blue $U$pward Cayley graph classes (Eqs.\ref{eq:lcmLU},\ref{eq:luclass}), e.g. $\colorbox{cyan}{\textcolor{black}{16}}]_{288}$, of branching blue $U$pward Collatz $3n+1$ nodes (Fig.\ref{fig:tge0}a). 

\paragraph{The leftward function $L$ and the upward function $U$}

\noindent Incoming blue $g$-arrows with path lenght $p$ after a first red arrow (Fig.\ref{fig:cc3p1}) specify the leftward function $L=fg^p$ (Def.\ref{def:Lw}) for six light-purple nodes indicating the branching subclasses $c_{bs}={\left[4,22,40,\ \ 16,34,52\right]}_{54}$ (Eq.\ref{eq:bsclass}) yielding $[n_{bs}]_{288}=2^{\max_{fg}-p}=2^{4-p}$ branching classes $[c_{bs}]_{288}=[4,1,4,\;2,8,8]$  (Def. \ref{def:Lw}).

\newpage
\begin{definition} \textit{Leftward function} $L:S_{\ge 0}\textcolor{red}{\to} s_{0;\ge 1};\;\;n\textcolor{red}{\to} (nfg^p=(n-1)/3\cdot 2^p);\;\; n_{bs}=2^{\max_{fg}-p}=2^{4-p}$ 
\label{def:Lw}
\end{definition}
\small{
\begin{tabular}
{p{0.25in}p{0.02in}p{0.01in}p{0.96in}p{0.01in}p{0.25in}p{0.01in}p{2.4in}p{0.01in}p{0.1in}}
 \hline 
$S_{\ge 0}$\newline  & $p$\textit{} & : & $L(n)$\textit{} & ; & $s_{0;\ge 1}$ & $=$ & ${\left[c27\right]}_{288}$ $=$ ${\left[c27\left(15,12\right)\right]}_{288}$ & ; & $n_{bs}$\textit{} \\ \hline 
${\left[4\right]}_{54}$\textit{} & 2 & : & $nfg^2=4\left(n-1\right)/3$\textit{}& ; & ${\left[4\right]}_{72}$\textit{} & $=$\textit{} & ${\left[4,\ 76,\ 148,\ 220\right]}_{288}\subset\left[4\right]_{18}$\textit{} & ; & 4 \\  
${\left[22\right]}_{54}$\textit{} & 4 & : & $nfg^4=16(n-1)/3$\textit{}& ; & ${\left[112\right]}_{288}$\textit{} & $=$\textit{} & ${\left[112\right]}_{288}\subset\left[4\right]_{18}$\textit{} & ; & 1\textit{} \\  
${\left[40\right]}_{54}$\textit{} & 2 & : & $nfg^2=4\ (n-1)/3$\textit{} & ; & ${\left[52\right]}_{72}$\textit{} & $=$\textit{} & ${\left[52,\ 124,\ 196,\ 268\right]}_{288}\subset\left[16\right]_{18}$\textit{} & ; & 4\textit{} \\ \hdashline  
${\left[16\right]}_{54}$\textit{} & 3 & : & $nfg^3=8\ (n-1)/3$\textit{}& ; & ${\left[40\right]}_{144}$\textit{} & $=$\textit{} & ${\left[40,\ 184\right]}_{288}\subset\left[4\right]_{18}$\textit{} & ; & 2\textit{} \\  
${\left[34\right]}_{54}$\textit{} & 1 & : & $nfg^1=2\ (n-1)/3$\textit{} & ; & ${\left[22\right]}_{36}$\textit{} & $=$\textit{} & ${\left[22,58,94,130,166,202,238,274\right]}_{288}\subset\left[4\right]_{18}$\textit{} & ; & 8\textit{} \\  
${\left[52\right]}_{54}$ & 1 & : & $nfg^1=2\ (n-1)/3$\textit{}& ; & ${\left[34\right]}_{36}$\textit{} & $=$\textit{} & ${\left[34,70,106,142,178,214,250,286\right]}_{288}\subset\left[16\right]_{18}$\textit{} & ; & 8\textit{} \\ \hline 
\end{tabular}
}
\vspace{3mm}
\normalsize{}

\noindent The powers $p=[[2,4,2],[3,1,1]$ for the six subclasses modulo $[[4,22,40],[16,34,52]]_{54}$ yield $a^3v^2=3^3\cdot1^2=27$ (Eq.\ref{eq:luclass}) leftward output classes with intrinsic periodicities $[[72=2^3 3^2, 288=2^5 3^2,  72=2^3 3^2],  [144=2^4 3^2, 36=2^2 3^2]$. These intrinsic periodicities have as their Least Common Multiple periodicity LCMp $=2^5 3^2=288$ (Eq.\ref{eq:lcmLU}). The leftward subclasses vector~$[c_{bs}]_{288}=[4,1,4,2,8,8]$ indicates that in total $a^3=27$ (Eq.\ref{eq:luclass}) leftward subclasses occur in each successive LCM$_{LU}=288$ period (Eq.\ref{eq:lcmLU}). Definition \ref{def:Lw} shows the meaning of the notation $c27(15,12)$, which is that each successive LCM$_{LU}=288$ comprising 27 $L$eftward numbers comprises  $4+1+2+8=15$ class $[4]_{18}$ subclasses and $4+16=12$ class $[16]_{18}$ subclasses. 

The upward function $U=g^q$ (Def.\ref{def:Uw}) is specified differently for the two branching classes $S_{\ge 0}={\left[4,16\right]}_{18}$. In Fig.\ref{fig:cc3p1}, 
the $p$-number on the first $g:n\to 2n$ arrow in an incoming blue path of $q$ successive blue $g$ arrows shows the power $q$ to which $2$ is raised in $U=g^p: n\to n\cdot 2^q$. 

\begin{definition} \label{def:Uw} \textit{Upward function} $U:\;S_{\ge 0}\textcolor{blue}{\to} S_{\ge 1};\;\;n\textcolor{blue}{\to} n\cdot 2^q$
\end{definition}
\small{
\begin{tabular}{p{0.3in}p{0.1in}p{0.05in}p{0.85in}p{0.1in}p{0.4in}p{0.1in}p{1.6in}p{0.02in}p{0.12in}} \hline 
$S_{\ge 0}$\newline  & \textit{q} & : & $U(n)$ & ; & $S_{\ge 1}$  & $=$ & ${\left[c5\right]}_{288}$ $={\left[c5\left(1,4\right)\right]}_{288}$   & ; & $n_{b}$  \\ \hline 
${\left[4\right]}_{18}$ & 2 & : & ${ng}^2=4n$ & ; & ${\left[16\right]}_{72}$ & $=$ & ${\left[16,\ 88,\ 160,\ 232\right]}_{288}\subset\left[16\right]_{18}$ & ; & $4$ \\ 
${\left[16\right]}_{18}$ & 4 & : & ${ng}^4=16n$ & ; & ${\left[256\right]}_{288}$ & $=$ & ${\left[256\right]}_{288}\subset\left[4\right]_{18}$ & ; & $1$ \\ \hline 
\end{tabular}
}
\normalsize{}

\noindent The powers $q=2$ and $q=4$ specifying the upward function $U=g^q$ for the classes $[4]_{18}$ and $[16]_{18}$ by upward lengths of paths without red arrows (Fig.\ref{fig:cc3p1}) yield upward output classes $[16]_{72}$ and $[256]_{288}$ with intrinsic periodicities $72=2^3 3^2$ and $288=2^5 3^2$. The upward classes vector~$[c_{b}]_{288}=[4,1]$ indicates that in total $A^+=(a^2-a-1)v^2=(9-3-1)1=5$ (Eq.\ref{eq:Am}) upward classes occur in each successive LCM period of $\textrm{lcm}_{}\cdot2^{\max_{fg}}=288$ (Eq.\ref{eq:lcmLU}) natural numbers. The notation $c5(1,4)$ means that branching class $[4]_{18}$ comprises $4$ and branching class $[16]_{18}$ comprises $1$ upward class.  The notation $S_{\ge1}=c5(1,4)=[[256], [16,88,160,232]_{288}$ means that in upward subtree $T_{\ge1}$ one subclass, $[256]_{288}$, is a subclass of $[4]_{18}$, while $4$ subclasses belong to branching class $[16]_{18}$. The periodicity expansion by one upward iteration, denoted as $\theta_{U1}=288/18=2^{2\min_q}=2^4$ (Eqs.\ref{eq:minq},\ref{eq:maxpq}), is obtained by dividing the periodicity of upward output numbers, which is $288$ (Eq.\ref{eq:lcmLU}), by the periodicity of branching number arguments, which is $18$ (Eq.\ref{eq:lcmpq}). The periodicity expansion by $a-1=2$ upward iterations is invariably $U^2:n\textcolor{blue}{\to}2^6n=64n$, as can be seen from the blue upward cycle length of six $g:n\textcolor{blue}{\to}2n$ blue arrows in Figure \ref{fig:cc3p1}, connecting the two purple classes to themselves. Applying $U^2(n):n\textcolor{blue}{\to}64n$ to all $3n+1$ numbers generates $a-1=c-2=2$ infinite geometric series (Eq.\ref{eq:ucojcm}) by successive \textit{odd upward iterations} $LU^{1,3,5,...}$ and successive \textit{even upward iterations} $LU^{2,4,6,...}$ in blue upward $S$-generation classes (Eq.\ref{eq:dtT}).

Fig.\ref{fig:cc3p1} also shows the sequences of $f$, $g$,$f^{-1}$, and $g^{-1}$-arrows that specify the composite $F$orward function $F$ for a walk of numbers from the $16$ non-branching classes to a number from the two purple classes (Def.\ref{def:Fw}). Arrows $g^{-1}$ and $f^{-1}$ that bring numbers immediately closer to the trivial root show up in the paths of number classes $[0,3]_6=[\{0,6,12\},\{3,9,15\}]_{18}$ divisible by $3$.

\newpage
\begin{definition} 
\label{def:Fw}
Forward function $F:{\left[0,1,2,3,\ \ 5\dots 15,\ \ 17\right]}_{18}\to S_{\ge 0};\ n\to F(n)$
\end{definition}
\small{
\begin{tabular}{p{1.4in} p{0.1in} p{0.1in} p{0.7in} p{0.1in} p{1.4in} p{0.4in}} \hline $n\in $
non-branching class  & $r$\textit{} & : & \textit{} $F(n)$& & & $S_{\ge 0}$ \\ \hline 
${\left[3\right]}_6\ ={\left[3,9,15\right]}_{18}$\textit{} & 2 & : & ${nf}^{-1}g^2$\textit{}& $=$ & $12n+4$\textit{} & ${\left[4\right]}_{18}$ \\ 
${\left[0\right]}_{6}={\left[0,16,12\right]}_{18}$ & 2 & : & ${ng}^{-i}\ f^{-1}g^2$& $=$ & $12n\;\,/\;\,2^{i=1,2,3,\dots}+4$ & ${\left[4\right]}_{18}$ \\ 
${\left[2,8,11,17\right]}_{18}$\textit{} & 1 & : & $ng^1$\textit{}& $=$ & $2n$\textit{} & ${\left[4,16\right]}_{18}$ \\ 
${\left[1,10,13\right]}_{18}$ & 2 & : & $ng^2$\textit{}& $=$ & $4n$\textit{} & ${\left[4,16\right]}_{18}$ \\  
${\left[5,14\right]}_{18}$ & 3 & : & $ng^3$& $=$ & $8n$ & ${\left[4\right]}_{18}$ \\ 
${\left[7\right]}_{18}$ & 4 &  & $ng^4$& $=$ & $16n$ & ${\left[4\right]}_{18}$ \\ \hline 
\end{tabular}
}
\normalsize{}
\vspace{3mm}

\noindent 
\begin{lemma} 
\label{lem:walk1}
Every number in each of the $16$ non-branching Collatz $3n+1$ classes ${\left[0,\dots ,3,5,\dots ,15,17\right]}_{18}$ can be taken on board in the Hamiltonian tour from and back to the trivial root $\cabroot$ taking on board numbers from QI- and QII-coordinates in the Cartesian Upper Half graph $C_H$. \emph{The forward function $F$ (Def.\ref{def:Fw}) specifies for all non-branching Collatz numbers $n$ a walk, via a $[10]_{18}$ number, to a Cayley tree branching number $F(n)\in[4]_{18}\subset (S_{\ge 0}=[4,16]_{18})$ on a $QI$- or $QII$ coordinate (Theorem \ref{the:cartesian}).} \qed
\end{lemma}

\begin{table}[!htbp]
\caption{Cotree congruence classes tape: 2 a\textit{n}+b proofs by 1 red and \textit{a}-1 blue cotree density columns}
\label{fig:Gccc}
\noindent \centerline{\includegraphics[width=0.88\textwidth]{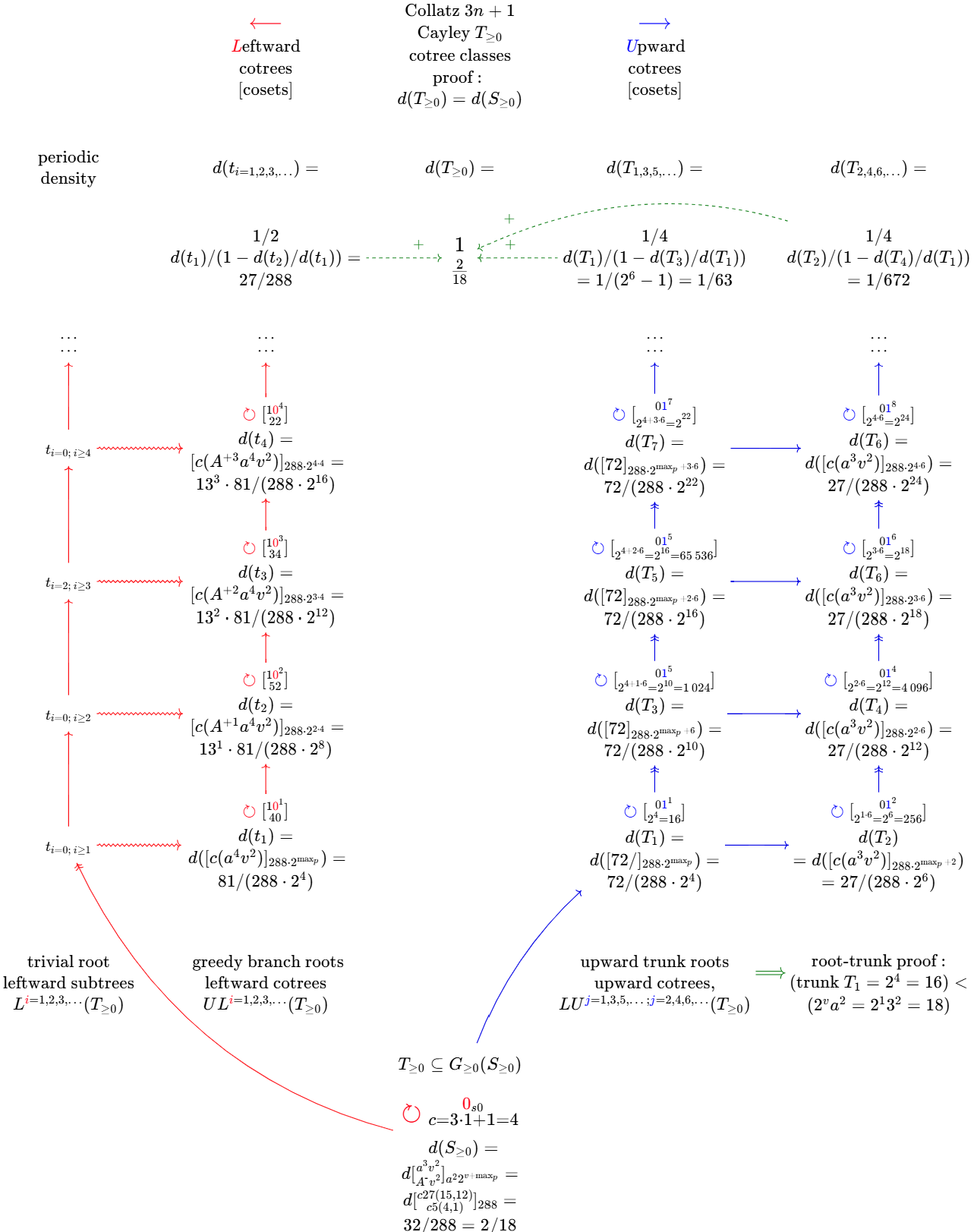}}
\centering\begin{minipage}{0.9\linewidth}
\noindent \textbf{Legend \textbar }\small{ The Collatz 3$n$+1 cotree congruence classes graph $T_{CCC}$ shows $a=3$ geometric series of cotree densities, 1 red series of cotree densities $d(t_i)$ generated by $UL^{i=1,2,3,\dots}$ iterations and $a-1=2$ blue series of cotree densities: an odd upward series $d(T_{j=1,3,5,\dots})$ generated by $LU^{j=1,3,5,\dots}$ iterations, and an even upward series $d(T_{j=2,4,6,\dots})$ generated by $LU^{j=2,4,6,\dots}$ iterations. The generation of red leftward cotrees $UL_i(T_{\ge0})\to t_i$  via red leftward subtrees $L_i(T_{\ge0}) \to t_{0;i\ge0}$ (subsection \ref{sec:34eqs} gives a multiple of congruence classes in successive cotrees $t_i$ of $A^{+i-1}a^4v^2$ classes (subsection \ref{sec:34eqs}, Fig.\ref{fig:ltree}, Eqs.\ref{eq:lcmLU},\ref{eq:Apm})  classes, in which $a=3$, $v=1$, $A^{+1}=a^2+a+1=13$ and $A^-=a^2-a-1=5$. The red, the odd blue and even blue geometric cotree series densities amount to a cumulative cotree density of $\frac{27}{288}+\frac{1}{63}+\frac{1}{672}=\frac{32}{288}=\frac{2}{18}$ which is the density of all number nodes $d(S_{\ge0})=d[\sm{a^3v^2\\A^-v^2}]_{a^22^{v+max_{pq}}}=32/288=2/18$ 
(Eq.\ref{eq:dSge0}) generated by the $3n+1$ function, proving the Collatz $3n+1$ conjecture. This implies the root-trunk-proof that the lowest trunk number $U(c)=16$ in cotree $T1$ must be a period $\ell=0$ number $0.18+16=16$ lower than its periodicty $2a^2v=18$.} \textbf{\textbar}
\end{minipage}\end{table}

\paragraph{Periodicities of Collatz $3n+1$ leftward and upward subtrees and cotrees}

\begin{corollary}
\label{cor:dens2} LCM \textit{periodicities of leftward and upward subsets and cosets} (Figs.\ref{fig:ltree}, \ref{fig:utree})
\vspace{-0.5em}
\begin{align}
{\theta}_{U2}=2^{3\cdot 2},&\textrm{ LCMp even upward subtrees: } &U^{j=2,4,6,\cdots} &:[S_{\ge 0}]_{288}\textcolor{blue}{\to} [S_{\ge j}]_{288\cdot 2^{3j}}   \label{eq:esLcmp}\\
&\;\;\textrm{ \textcolor{white}{LCMp even upward }cotrees: } &LU^{j=2,4,6,\cdots} &:[S_{\ge 0}]_{288}\textcolor{blue}{\to} [S_j]_{288\cdot 2^{3j}}   
\label{eq:ecLcmp}\\
{\theta}_{U1}=2^4,{\theta}_{U2}=2^{3\cdot 2}, &\textrm{ LCMp odd upward subtrees: } &U^{j=1,3,5,\cdots} &:[S_{\ge 0}]_{288}\textcolor{blue}{\to} [S_{\ge j}]_{288\cdot 2^{3j+1}} 
\label{eq:osLcmp} \\
 &\textrm{\textcolor{white}{  LCMp odd upward }\;\;\,cotrees: } &LU^{j=1,3,5,\cdots} &:[S_{\ge 0}]_{288}\textcolor{blue}{\to} [S_j]_{288\cdot 2^{3j+1}}  
 \label{eq:ocLcmp} \\
3{\theta }_R=2^4,\;\;\; &\textrm{ LCMp leftward subtrees:} &L^{i=1,2,3,\cdots} &:[S_{\ge 0}]_{288}\textcolor{red}{\to} [s_{0;\ge i}]_{3\cdot 288\cdot 2^{4i}} 
\label{eq:lsLcmp}\\
 &\textrm{ \textcolor{white}{LCMp leftward }\, cotrees: } &UL^{i=1,2,3,\cdots} &:[S_{\ge 0}]_{288}\textcolor{red}{\to} [s_i]_{3\cdot 288\cdot 2^{4i}} 
 \label{eq:lcLcmp}
\end{align}
\end{corollary}

\noindent Table \ref{tab:classes} lists the four sets $c351,c81,c72$ and $c8$ of the congruence classes of the nodes of the first two subtrees $t_{\ge 1}$ and $T_{\ge 1}$. These four sets underlie the classes of all their successor subtrees and successor cotrees. Periodic \textit{cotree density} is calculated as the number of classes relative to their Least Common Multiple periodicity. The automorphism graph (Fig.\ref{fig:aut}) is annotated with the cotree densities  $d\left(T_1\right)=72/(288\cdot 2^4)$, $d\left(T_2\right)=27/(288\cdot 2^6)$ and $d\left(t_1\right)=81/(288\cdot 2^4)$. The notation ${\left[c72\right]}_{288\cdot 2^4}$ for the classes of $T_1$ elucidates that $d\left(T_1\right)=72/(288\cdot 2^4)$ indeed. The geometric sum formula $s=a/(1-r)$ underlies the cumulative densities $1/672+1/63=(A^-=a^2-a-1=5)/288$ of even and odd upward cotree numbers, and $(a^3=27)/288$ of leftward cotree numbers (Eq.\ref{eq:lcmLU}), giving the Lean-certified `cotree density proof' $d(tT)=d(S_{\ge0})$ (Eq.\ref{eq:cdp}) of the Collatz $3n+1$ conjecture $(5+27)/288=2/18$ (Eq.\ref{eq:dSge0} below, Table \ref{tab:cdp3p1}).

Table \ref{fig:Gccc} shows how a Collatz-Conway a$n$+b conjecture for which the cyclic root $c=a\cdot1+1$ is a power-of-2 and $c^2$ is the lowest trunk number, exemplified by the Collatz 3$n$+1 conjecture, is proven both by the cotree-density-test on an infinite series of periods $\ell=0,1,2,\dots $ of $a\nmid$ an$_{\textup{odd}}$+b numbers and the root trunk-test on the finite subset of period $\ell=0$ numbers.

\newpage 
\begin{table}
  \caption{Six base sets of Collatz $3n+1$ congruence classes of subtrees and cotrees}\label{tab:classes}
\end{table}

\small{
\noindent $S_{\ge 0}= \left[4,16\right]_{18}$, argument classes $U$, upward branching classes vector to lcm-p$=$288 is~$[n_b]_{288}=\left[4,1\right]$
}

\noindent $S_{\ge0}=\left[4,16\right]_{18}=\left[4,22,40,16,34,52\right]_{54}$,
$\mathrm{argument\ subclasses\ \ }L\mathrm{,\ }$subclasses vector to lcm-p is ~$[n_{bs}]_{288}=\left[4,2,1,8,4,8\right]$

\vspace {6mm}

\noindent $T_{\ge 0}{\left[c27(15,12),c5(1,4)\right]}_{288}$, leftward and upward classes$\ T_{\ge 0}\left[s_{\ge 0},S_{\ge 1}\right]$ by argument classes ${\left[4,16\right]}_{18}$

\vspace {2mm}

${\left[c27(15,12)\right]}_{288}$ = [[ 4, 22, 40, 58,  76,  94, 112,130,148,166,184,202,220,238,274],

\hspace{3.6cm}  [34,52, 70,106,124,142,178,196,214,250,268,286]$]_{288}$

${\left[c5(1,4)\right]}_{288}$ = [[256],[16,88,160,232$]_{288}$  

\vspace {5mm}

\noindent $T_{\ge 1}{\left[c72,c8\right]}_{288\cdot 2^4}$ upward coset $S_1$ and subset $S_{\ge 2}$ in the first (odd) upward subtree $T_{\ge 1}\left[S_1,S_{\ge 2}\right]$ 

\vspace {2mm}

$\left[c72\right]_{288\cdot 2^4}$ = [ 
     16,     88,  160,    232,   304,   376,   448,   520,   544,   592,   664,   736,

\noindent  \hspace{2.4cm}  808,   832,  880,    952, 1096, 1120, 1168, 1240, 1312, 1384, 1456, 1528,

\noindent  \hspace{2.4cm} 1600, 1672, 1696, 1744, 1816, 1888, 1960, 1984, 2032, 2104, 2248, 2272,

\noindent  \hspace{2.4cm} 2320, 2392, 2464, 2536, 2608, 2680, 2752, 2824, 2848, 2896, 2968, 3040,

\noindent  \hspace{2.4cm} 3112, 3136, 3184, 3256, 3400, 3424, 3472, 3544, 3616, 3688, 3760, 3832,

\noindent  \hspace{2.4cm} 3904, 3976, 4000, 4048, 4120 ,4192, 4264, 4288, 4336, 4408, 4552, 4576 $]_{288\cdot 2^4}$

$\left[c8\right]_{288\cdot 2^4}$ = \hspace{0.25cm} [256,  1024, 1408, 2176, 2560, 3328, 3712, 4480 $]_{288\cdot 2^4}$

\vspace {6mm}

\noindent $t_{0;\ge 1}{\left[3c351,3c81\right]}_{288\cdot {3\cdot 2}^4}$ leftward subset $s_{0;\ge 2}$ and coset $s_1$ in the first leftward subtree $t_{0;\ge 1}[s_{0;\ge 1},s_1$]

\vspace {2mm}
\footnotesize{
$\left[c351\right]_{288\cdot 2^4}$ = [ 
       4,     22,     34,     52,     70,     76,     94,   112,   124,   130,   142,   148,   166,   178,   

\noindent  \hspace{2.4cm}    184,   196,   214,   220,   238,   268,   274,   286,   292,   310,   322,   328,   358,

\noindent  \hspace{2.4cm}    364,   382,   400,   412,   418,   430,   436,   454,   466,   472,   484,   502,   508,   526,

\noindent  \hspace{2.4cm}   556,   562,  574,    580,   598,   610,   628,   646,   652,   670,   688,   700,  706,

\noindent  \hspace{2.4cm}    718,   742,  754,    760,   772,   790,   796,   814,   844,   850,   862,   868,   886,   898, 

\noindent \hspace{2.4cm}   904,   916,  934,    940,   958,   976,   988,   994, 1006, 1012, 1030, 1042, 1048,

\noindent \hspace{2.4cm}   1060, 1078, 1084, 1102, 1132, 1138, 1150, 1156, 1174, 1186, 1204, 1222, 1228, 1246,

\noindent \hspace{2.4cm}   1264, 1276, 1282, 1294, 1300, 1318, 1330, 1336, 1348, 1366, 1372, 1390, 1420,

\noindent \hspace{2.4cm}   1426, 1438, 1444, 1462, 1474, 1480, 1510, 1516, 1534, 1552, 1564, 1570, 1582, 1588,

\noindent \hspace{2.4cm}   1606, 1618, 1624, 1636, 1654, 1660, 1678, 1708, 1714, 1726, 1732, 1750, 1762,

\noindent \hspace{2.4cm}   1780, 1798, 1804, 1822, 1840, 1852, 1858, 1870, 1894, 1906, 1912, 1924, 1942, 1948,

\noindent \hspace{2.4cm}   1966, 1996, 2002, 2014, 2020, 2038, 2050, 2056, 2068, 2086, 2092, 2110, 2128,

\noindent \hspace{2.4cm}   2140, 2146, 2158, 2164, 2182, 2194, 2200, 2212, 2230, 2236, 2254, 2284, 2290, 2302,

\noindent \hspace{2.4cm}   2308, 2326, 2338, 2356, 2374, 2380, 2398, 2416, 2428, 2434, 2446, 2452, 2470,

\noindent \hspace{2.4cm}   2482, 2488, 2500, 2518, 2524, 2542, 2572, 2578, 2590, 2596, 2614, 2626, 2632, 2662,

\noindent \hspace{2.4cm}   2668, 2686, 2704, 2716, 2722, 2734, 2740, 2758, 2770, 2776, 2788, 2806, 2812,

\noindent \hspace{2.4cm}   2830, 2860, 2866, 2878, 2884, 2902, 2914, 2932, 2950, 2956, 2974, 2992, 3004, 3010,

\noindent \hspace{2.4cm}   3022, 3046, 3058, 3064, 3076, 3094, 3100, 3118, 3148, 3154, 3166, 3172, 3190,

\noindent \hspace{2.4cm}   3202, 3208, 3220, 3238, 3244, 3262, 3280, 3292, 3298, 3310, 3316, 3334, 3346, 3352,

\noindent \hspace{2.4cm}   3364, 3382, 3388, 3406, 3436, 3442, 3454, 3460, 3478, 3490, 3508, 3526, 3532,

\noindent  \hspace{2.4cm}  3550, 3568, 3580, 3586, 3598, 3604, 3622, 3634, 3640, 3652, 3670, 3676, 3694, 3724,

\noindent  \hspace{2.4cm}  3730, 3742, 3748, 3766, 3778, 3784, 3814, 3820, 3838, 3856, 3868, 3874, 3886,

\noindent \hspace{2.4cm}   3892, 3910, 3922, 3928, 3940, 3958, 3964, 3982, 4012, 4018, 4030, 4036, 4054, 4066,

\noindent \hspace{2.4cm}   4084, 4102, 4108, 4126, 4144, 4156, 4162, 4174, 4198, 4210, 4216, 4228, 4246,

\noindent \hspace{2.4cm}  4252, 4270, 4300, 4306, 4318, 4324, 4342, 4354, 4360, 4372, 4390, 4396, 4414, 4444,

\noindent  \hspace{2.4cm} 4450, 4462, 4468, 4486, 4498, 4504, 4516, 4534, 4540, 4558, 4588, 4594, 4606 $]_{288\cdot 2^4}$ }

$\left[c81\right]_{288\cdot 2^4}$ = [ 
    40,      58,   106,   202,   250,   340,   346,   394,   490,   538,   616,   634,   682,   724,

\noindent  \hspace{2.4cm}    778,   826,   922,   970, 1066, 1114, 1192, 1210, 1258, 1354, 1402, 1492, 1498,

\noindent \hspace{2.4cm}   1546, 1642, 1690, 1768, 1786, 1834, 1876, 1930, 1978, 2074, 2122, 2218, 2266, 2344,

\noindent \hspace{2.4cm}   2362, 2410, 2506, 2554, 2644, 2650, 2698, 2794, 2842, 2920, 2938, 2986, 3028,

\noindent \hspace{2.4cm}   3082, 3130, 3226, 3274, 3370, 3418, 3496, 3514, 3562, 3658, 3706, 3796, 3802, 3850,

\noindent  \hspace{2.4cm} 3946, 3994, 4072, 4090, 4138, 4180, 4234, 4282, 4378, 4426, 4432, 4522, 4570 $]_{288\cdot 2^4}$   

\FloatBarrier
\normalsize{}

\paragraph {Congruence classes and periodic densities of even upward subtrees and cotrees}
Table \ref{tab:Lnba4} specifies for the Collatz $3n+1$ function the red $L$-arrow to  $L$eftward child number pairs (Eq.\ref{eq:redChild}), which is the first arrow in paths to all higher-generation $L$eftward subtrees (Eq.\ref{eq:lsuar}) and in paths to all $U$pward same-generation cotrees (Eq.\ref{eq:ucoa}).  Table \ref{tab:Lnba4}  shows the emergence in the first leftward subtree $t_{0;\ge1}$ of $\Sigma n_b=a^4$ first generation $L$eftward $s_1$-classes with a periodicity of $2^{5+(\max_{fg}=4)}3^2=4608$, resulting from the $L$eftward matrix multiplication function $L:T_{\ge0}([c_{bs}]_{18a=54})\cdot T_{\ge0}
([C_{bs}]_{54})=t_{0;\ge1}([c_{b}]_{18})$. The matrix multiplication output of $a=3$ times $15$ and $12$ leftward classes mounts up to $3(15+12)=3^4=81$ subclasses with a Least Common Multiple periodicity of $2^{5+(\max_{fg}=4)}3^2=4608$. 

\begin{table}[H]
\caption{$3n+1$ multiplication $L:T_{\ge0}([c_{bs}]_{18a=54})\cdot T_{\ge0}
([C_{bs}]_{54})=t_{0;\ge1}([c_{b}]_{18})$ gives $n_{b}=a^4$ classes $[c_{b}]_{4608}$}
\label{tab:Lnba4}
\[
\left(\begin{array}{c c|c}  
T_{\ge0}:  &  &  \\  
{[c_bs]_{288\cdot2^{\textup{lcm}_{pq}}=4608}} & {[c_{bs}]_{18a}}  &{[n_{bs}]_{18\cdot2^{\textup{max}_{pq}}=288}}\\ \hline 
 &[4]_{54}&4  \\ 
{[4]_{18}} & {[22]_{54}}&2 \\ 
 &{[40]_{54}} &4\\  \hdashline
 &{[16]_{54}} &8\\ 
{[16]_{18}}& {[34]_{54}}&1 \\  
 & {[52]_{54}}&8 \\ \hline 
\end{array}\right) 
\cdot
\left(\begin{array}{c c|c c c : c c c}  
T_{\ge0}:C_{bs}  &        &        & {[4]_{18}} &        &        & {[16]_{18}} &  \\  
         &        & {[4]_{54}} & {[22]_{54}} & {[40]_{54}} & {[16]_{54}} & {[34]_{54}} & {[52]_{54}}\\ \hline 
         & {[4]_{54}}  & 1 & 1 & 1 & 0 & 0 & 0 \\ 
{[4]_{18}} & {[22]_{54}} & 1 & 1 & 1 & 0 & 0 & 0 \\ 
         & {[40]_{54}} & 0 & 0 & 0 & 1 & 1 & 1 \\ \hdashline 
         & {[16]_{54}} & 1 & 1 & 1 & 0 & 0 & 0 \\ 
{[16]_{18}}& {[34]_{54}} & 1 & 1 & 1 & 0 & 0 & 0 \\  
         & {[52]_{54}} & 0 & 0 & 0 & 1 & 1 & 1 \\ \hline
\end{array}\right)
=
\] 
\end{table}
\vspace{-\baselineskip}
\begin{table}[H]
\[\textcolor{white}{
\left(\begin{array}{c c|c}  
  &  &  \\  
 &   &   \\ \hline 
 &{[4]_{54}}&4  \\ 
\end{array}\right) 
\cdot==}
\textcolor{white}{XXXX}
\left(\begin{array}{c | c|c c c : c c c}  
t_{0;\ge1}: &     &        & {[4]_{18}} &        &        & {[16]_{18}} &  \\  
 &\Sigma n_{bs}=a^4   & {[4]_{54}} & {[22]_{54}} & {[40]_{54}} & [16]_{54} & {[34]_{54}} & {[52]_{54}}\\ \hline 
{[c_{bs}]_{288\cdot2^{\textup{max}_{pq}}=4608}}    & (15+12)a=81  & 15 & 15 & 15 & 12 & 12 & 12 \\ 
\end{array}\right)=
\] 
\end{table}
\vspace{-\baselineskip}
\begin{table}[H]
\[\textcolor{white}{xxxxxxxxx
\left(\begin{array}{c c|c}  
 T_{\ge0}:C_{bs} &  &  \\  
 &   &   \\ \hline 
 &{[4]_{54}}&4  \\ 
\end{array}\right) 
\cdot==}
\left(\begin{array}{c c|c : c}  
t_{0;\ge1}:   &        &   {[4]_{18}} &  {[16]_{18}}     \\  
  & \Sigma n_b=a^4       & & \\ \hline 
 {[c_{b}]_{288\cdot2^{\textup{max}_{pq}}=2^93^2=4608}}  & 27a=81  & 15a=45 & 12a=36  \\ 
\end{array}\right)\textcolor{white}{}\,
\] 
\end{table}
\noindent Table \ref{tab:UnbsAm} specifies for the Collatz $3n+1$ function the blue $U$-arrow to  $U$pward child number pairs (Eq.\ref{eq:blueChild}), which is the first arrow in paths to all higher-generation $U$pward subtrees (Eq.\ref{eq:usuar}) and in paths to all $L$eftward same-generation cotrees (Eq.\ref{eq:lcoa}). Table \ref{tab:UnbsAm} specifies the emergence in the first $U$pward subtree $T_{\ge1}$ of $\Sigma n_{bs}=A^-a=(a^2-a-1)\cdot3=15$ first generation upward $S_1$-subclasses with a periodicity of $2^{5+(\max_{fg}=4)}3^2=4608$, given by the Upward matrix multiplication function $L:T_{\ge0}([c_{bs}]_{18a=54})\cdot T_{\ge0}
([C_{bs}]_{54})=t_{0;\ge1}([c_{b}]_{18})$. The matrix multiplication output of $a=3$ times $4$ and $1$ upward subclasses $c_{bs}$ mounts up to $n_{bs}=3\cdot5=15$ upward subclasses, with a Least Common Multiple periodicity of $2^{5+(\max_{fg}=4)}3^2a=4608$. 

\begin{table}[H]\caption{$3n+1$ multiplication $U:T_{\ge0}([c_b]_{18})\cdot T_{\ge0}([C_b]_{18a=54})=T_{\ge1}([c_{bs}]_{54})$ gives $n_{bs}=A^-a$ subclasses  $[c_{bs}]_{288}$}
\label{tab:UnbsAm}
\[
\left(\begin{array}{c|c}  
T_{\ge0}:  &   \\  
{[c_b]_{\textup{lcm}_{pq}=18}} &{[n_{b}]_{18\cdot2^{max_{fg}}=288}      }\\ \hline 
{[4]_{18}}&4  \\ 
{[16]_{18}} &1 \\ 
\end{array}\right) 
\cdot
\left(\begin{array}{c c|c c}  
T_{\ge0}:C_{b}  &        &  {[4]_{18}}      & {[16]_{18}}  \\  
  &        &  & \\ \hline 
  & {[4]_{18}}  & 0 & 1  \\ 
  & {[16]_{18}} & 1 & 0  \\ 
\end{array}\right)
=
\] 
\end{table}
\vspace{-\baselineskip}
\begin{table}[H]
\[\textcolor{white}{
\left(\begin{array}{c| c}  
 T_{\ge0}:C_{bs} &   \\  
 &     \\ \hline 
 &{[4]_{54}}  \\ 
\end{array}\right)
\cdot==}
\left(\begin{array}{c c|c  c}  
T_{\ge1}:   &        &   {[4]_{18}} &  {[16]_{18}}     \\  
  & \Sigma n_b=A^-=      & & \\ \hline 
{[c_{b}]_{288\cdot2^{\max_{fg}}=4096}} &  a^2-a-1=5 & 4 & 1  \\ 
\end{array}\right)\,
=\textcolor{white}{XXXXXXl}
\] 
\end{table}
\vspace{-\baselineskip}
\begin{table}[H]
\[\textcolor{white}{
\left(\begin{array}{c c|c}  
  &  &  \\  
 &   &   \\ \hline 
 & & \\ 
\end{array}\right) 
}
\textcolor{white}{XX}
\left(\begin{array}{c | c|c c c : c c c}  
T_{\ge1}: &     &        & {[4]_{18}} &        &        & {[16]_{18}} &  \\  
 &\Sigma n_{bs}=A^-a=   & {[4]_{54}} & {[22]_{54}} & {[40]_{54}} & [16]_{54} & {[34]_{54}} & {[52]_{54}}\\ \hline 
{[c_{bs}]_{4096a=12\,288}}   & 5\cdot3=15  & 4 & 4 & 4 & 1 & 1 & 1 \\ 
\end{array}\textcolor{white}{XX}
\right)
\] 
\end{table}
\noindent Table \ref{tab:Lcdp} shows the generation of leftward cotrees $t_1,t_2,t_3,\dots$ by iterations $UL^{i=1,2,3,\dots}(T_{\ge0})$ from the binary tree $T_{\ge0}$---or iterations $L^{i=1,2,3,\dots}(T_{\ge1})$ from the first upward subtree $U_{\ge1}$ (Eq.\ref{eq:lcoa})---subclasses (Fig.\ref{tab:UnbsAm}). Their cumulative
$L$eftward cotree density $d(t_1,t_2,t_3,\dots)=a^3/{288}=27/{288}$ rests on   $A^+=((a+1)a+1)v^2=13$ more leftward congruence classes at each Leftward iteration.

\begin{table}[H]\caption{Collatz $3n+1$ $L$eftward cotree density $t_1,t_2,t_3,\dots$  ($a=3,\;v=1,\;A^+=((a+1)a+1)v^2=13$)}
\label{tab:Lcdp}
\[
\begin{array}{c |c : c c c c c| c c : c}  
i     &  &  1  & 2 & 3 & 4&\dots &  d \\ 
\textup{tree}  & T_{\ge0} &  t_1  &  t_2 &  t_3 & t_4 & \dots & (t_{i=1,2,3,\dots})  \\ \hline 
{[4]_{54}} &  4 & 15a=45 & 15aA^+=585 & 15aA^{+2}=7605 & 15aA^{+3}=98865&\dots & \\         
{[22]_{54}}&  1 & 15a & 15aA^+ & 15aA^{+2} &15aA^{+3}&\dots & \\  
{[40]_{54}}&  4 & 15a & 15aA^+ & 15aA^{+2} & 15aA^{+3}&\dots & \\ \hdashline        
{[16]_{54}}&  8 & 12a=36 & 12aA^+=468 & 12aA^{+2}=6084 & 12aA^{+3}=79092&\dots & \\  
{[34]_{54}}&  1 & 12a & 12aA^+ & 12aA^{+2} & 12aA^{+3}&\dots & \\
{[52]_{54}}&  8 & 12a & 12aA^+ & 12aA^{+2} & 12aA^{+3}&\dots & \\         \hdashline
\Sigma&  a^3\text{\scriptsize{=27}} & a^5\text{\scriptsize{=243}} & a^5A^+\text{\scriptsize{=3\,159}} & a^5A^{+2}\text{\scriptsize{=41\,067}} & a^5A^{+3}\text{\scriptsize{=533\,871}}&\dots & \\         \hline
\textup{LCMp}&  \textrm{\scriptsize{lcm}}_{LU}=2^53^2 &
2^{5+\textrm{max}_{pq}=9}a^3 & 2^{5+2\textrm{max}_{pq}=13}a^3 &2^{5+3\textrm{max}_{pq}=17}a^3  &2^{11+4\textrm{max}_{pq}=21}a^3 & & \\  
\textup{density }d&  \frac{27}{288} &  \frac{a^4=81}{288\cdot2^{4\cdot1}}  &  \frac{a^4A^+=1\,053}{288\cdot2^{4\cdot2}} & \frac{41\,067}{3\cdot288\cdot2^{4\cdot3}}& \frac{533\,871}{3\cdot288\cdot2^{4\cdot4}} &\cdots & \frac{a^3}{288}\,=\,\frac{27}{288}\\ 
\end{array} 
\] 
\end{table}

\noindent Table \ref{tab:Lcdp} shows that for the Collatz $3n+1$ function---with $v=1$, $a=3$, $A^+=(a+1)a+1=13, \text{lcm}_{pq}=4,\text{lcm}_{LU}=288$---the periodic density of leftward cotree $t_i$ is specified as 
\[d(t_i)=\frac{a^4 A^{i-1} v^2}{\text{lcm}_{LU}\cdot2^{i\cdot\text{max}_{pq}}}
=\frac{a^4 A^{i-1}}{288\cdot2^{\text{4}i}}
=\frac{3^4 13^{i-1}}{288\cdot2^{4i}}, \text{ cumulative density } 
d(t_{i=1,2,3,\dots})=
\frac{a^3v^2}{\text{lcm}_{LU}}
=\frac{27}{288}
\].

Similar to Table \ref{tab:Lcdp} for leftward cotrees, Table \ref{tab:Ucdp} shows the generation of $a-1=3-1=2$ series of odd upward cotrees $T_1,T_3,T_5,\dots$ and even upward cotrees $T_2,T_4,T_6,\dots$ by upward iterations $LU^{j=1,2,3,\dots}(T_{\ge0})$ from the binary tree $T_{\ge0}$ or upward iterations $U^{j=1,2,3,\dots}(T_{\ge1})$ from the first leftward subtree $t_{0;\ge1}$ with $n_b$ upward classes (Table.\ref{tab:Lnba4}; Eq.\ref{eq:ucoa})). Their cumulative
$U$pward cotree density $d(T_1,T_3,T_5,\dots,T_2,T_4,T_6,\dots)=\frac{A^-=(a-1)a-1}{\text{lcm}_{LU}}=\frac{5}{288}$ rests on $a-1=3-1=2$ geometric series of odd and even upward cotree densities.

\begin{table}[H]\caption{Collatz $3n+1$ $U$pward cotree density $T_1,T_3,T_5,\dots;T_2,T_4,T_6,\dots$  ($a=3,\;A^-=\,(a-1)a-1=5$)}
\label{tab:Ucdp}
\[
\begin{array}{c |c : c c: c c: c| c c : c}  
j     &  &  1  & 2 & 3 & 4&\dots & & & \\ 
\textup{tree}  & T_{\ge0} &  T_1  &  T_2 &  T_3 & T_4 & \dots &  \sm{d(T\\1,3,5,\dots)} & \sm{d(T\\2,4,6,\dots)} & \sm{d(T\\1,2,3,\dots)}   \\ \hline 
{[4]_{54}} &  15 & 4\cdot15 & 15 & 4
\cdot15 & 15&\dots & &\\         
{[16]_{18}}&  12 & 1\cdot12 & 12 & 1\cdot12 &12&\dots & &\\ \hdashline 
\Sigma &  27 & 72 & 27 & 72 &27&\dots & &\\  \hline
\textup{LCMp}&  \textrm{\scriptsize{lcm}}_{LU}=2^53^2 & 2^{5+2\min_q=9}a^2 & 2^{5+a\min_{q}=11}a^2 &2^{9+a\min_q=15}a^2  &2^{11+a\min_q=17}a^2 & &\frac{1}{2^{a\text{ min}_q}-1} 
&\frac{a}{7\;\text{lcm}_{LU}}
&\frac{A^-v^2}{\text{lcm}_{LU}}\\ 
\textup{density }d&  \frac{A^-\,=\;5}{288\cdot2^0} &  \frac{72}{288\cdot2^4}  &  \frac{27}{288\cdot2^{0+6}} & \frac{72}{288\cdot2^{4+6}}& \frac{27}{288\cdot2^{12}} &\cdots &  =\frac{1}{63} &=\frac{1}{672} & =\frac{5}{288}\\  
\end{array} 
\] 
\end{table}
\noindent Table \ref{tab:Ucdp} shows, in line with the Upward $3n+1$ function (Def.\ref{def:Uw}), that each iteration of the Leftward function shifts the distribution of subclasses $[4]_{18}:[16]_{18}$, from $27$ classes $15:12$ in Cayley tree $T_{\ge0}$ and its even upward cotrees $T_2,T_4,\dots$, to $72$ upward classes $4\cdot15:1\cdot12=72$ in odd upward cotrees $T_2,T_4,\dots$.

Table \ref{tab:cdp3p1} summarizes the previous two tables \ref{tab:Lcdp}  and \ref{tab:Ucdp}. The sum of cumulative densities gives the Lean-certified`cotree density proof' of the $3n+1$ conjecture $d(tT)=(27+5)/288=2/18$ (Eq.\ref{eq:cdp}, Tables \ref{tab:cdp3p1},\ref{fig:Gccc},\ref{fig:TLEE}). \noindent The `root successor test' proves the Collatz $3n+1$ conjecture by root paths for period $\ell=0$ numbers $7$ and $13$, of which the absence refutes the $3n-1$ conjecture (Eq.\ref{eq:rst}, Theorem \ref{the:decidab}.6). The Lean-certified `root-trunk test' proves the Collatz conjecture by its power-of-2 trunk $c=a\cdot1+b=2^2$ and its lowest trunk node $U(c)=2^{2+min_q}=2^{2+2}=c^2=16=((\ell=0)\cdot m+16)<(m=18)$ (Eq.\ref{eq:rtt}, Theorem \ref{the:decidab}.7, Tables \ref{fig:Gccc},\ref{fig:TLEE}).
 
\begin{table}[H]
\caption{Collatz  $3n+1$ density, 1 $L$eftward and \(3-1\) $U$pward geometric series of cotree densities}
\label{tab:cdp3p1}
\setlength{\abovedisplayskip}{0pt}
\setlength{\belowdisplayskip}{0pt}
\begin{align} 
\nonumber\\ \toprule
d(t_i)&\textrm{ density red $L$(eftward) cotree} & &\textcolor{white}{=}  \nonumber \\ 
d(t_1)&=\frac{a^4=81}{288\cdot2^{4}};\quad d(t_2)=\frac{a^4A^+=1053}{288\cdot2^{8}}
&  \sum_{i=1}^{\infty} d(t_{i=1,2,3,\cdots})&
=\frac{
    \frac{a^4}{288\cdot2^{4}}}
    {1-\frac{A^+v=13}{2^4}}
=\frac{a^3}{288}
\nonumber \\ 
&\textrm{Cumulative density leftward }3n+1\textrm{ cotrees} & \textcolor{white}{\sum_{i=1}^{\infty} d(t_{i=1,2,3,\cdots})}&=\frac{27}{288}\nonumber \\ \midrule \nonumber 
d(T_j)&\textrm{ density blue $U$(pward) cotree} & &\textcolor{white}{=}  \nonumber \\ 
d(T_1)&=\frac{72}{288\cdot2^{4}};\quad d(T_3)=\frac{72}{288\cdot2^{4+6}} &  
\sum_{i=0}^{\infty} d(T_{i=1,3,5,\cdots})&=\frac{1}{2^{6}-1}=\frac{1}{63} \nonumber  \\
d(T_2)&=\frac{27}{288\cdot2^{6}};\quad d(T_4)=\frac{27}{288\cdot2^{6+6}} &
\sum_{i=0}^{\infty} d(T_{i=2,4,6,\cdots})&=\frac{3}{288\cdot7}=\frac{1}{672}\nonumber \\ 
&\textrm{Cumulative density upward }3n+1\textrm{ cotrees} & \sum_{j=1}^{\infty} d(T_{j=1,2,3,\cdots})&=\frac{A^-v^2=(a-1)a-1=5}{288} \nonumber 
\\ \midrule &\textrm{Cotree density proof Collatz 3$n$+1} &  d(tT) &=\frac{27+5}{288}=\frac{2}{18}=d(S_{\ge0}) & \qed \nonumber
\\  &\textrm{Lean certificate: Collatz$\_$cotree$\_$density} &  & & \nonumber
\\\bottomrule \nonumber
\end{align}
\end{table}
\FloatBarrier

\subsubsection{Bařina $7n\pm1$:  Hamiltonian tour cotree density proof and root-trunk-proof}
\label{sec:ccba}

\begin{figure}[b!]
\caption{Bařina $7n\pm1$ congruence classes graph}
\label{fig:cc7pm1}
\begin{center}
\includegraphics[width=0.9\textwidth]{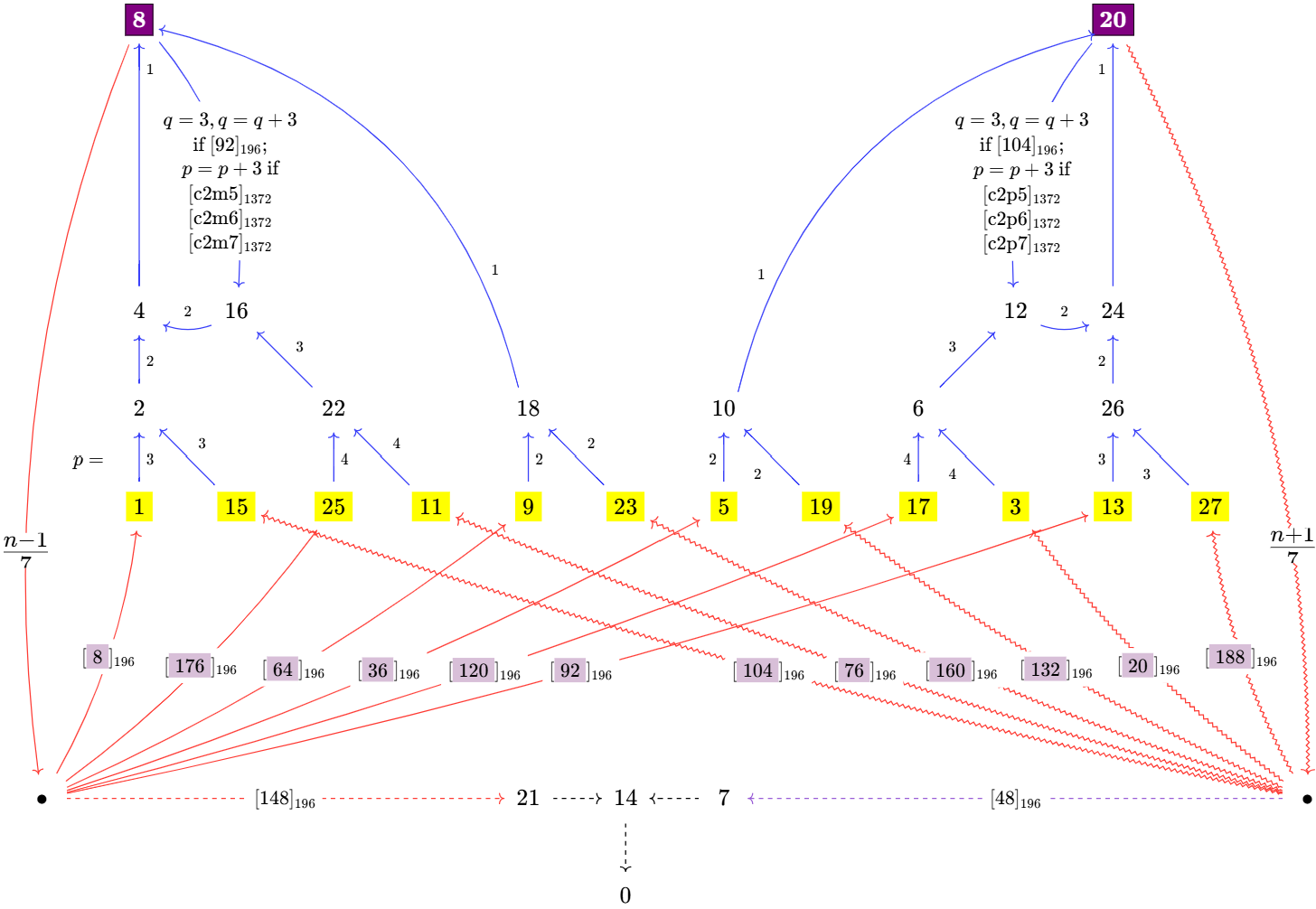}
\end{center}
\vspace{-0.4cm}
\centering\begin{minipage}{0.9\linewidth}
\noindent\textbf{Legend: } 
The condensed $7n\pm1$ congruence graph shows : $2$ even, purple, branching classes with periodicity/modulus $m=(2a^2v)/a=28$ (Eq.\ref{eq:lcmpq}), pictured as
$\colorbox{violet}{\textcolor{white}{8}}_{28}$  and $\colorbox{violet}{\textcolor{white}{20}}_{28}$, and $(a(a-1)v)/a\,=\,12$ odd, yellow, classes, e.g.
$\colorbox{yellow}{\textcolor{black}{25}}_{28}$, (Eq.\ref{eq:oclass}) 
from $(a(a-1)v)/a\,=\,12$ even, lightpurple, branching parent subclasses with periodicity $196$ (Eq.\ref{eq:bsclass}).
Red arrow $\textcolor{red}{\xrightarrow{\textcolor{black}{f:\;n\,\to\,(n\pm1)/7}}}$ (Eqs.\ref{eq:rootPairs},\ref{eq:oclass}) to odd yellow child annotated with its lightpurple (here lightgray) parent subclass, e.g.
$\colorbox{violet}{\textcolor{white}{8}}_{28}
\textcolor{red}{\xrightarrow{\colorbox{lightgray}{\textcolor{black}{[176]}}}}
\colorbox{yellow}{\textcolor{black}{25}}_{28}$.
Blue arrow $\textcolor{blue}{\xrightarrow{\textcolor{black}{g:\;n\,\to\,2n}}}$, either in a blue yellow-to-purple path with path length  $\textcolor{blue}{\xrightarrow{\textcolor{black}{p=\cdots, \textrm{ defining } L=fg^p}}}$, e.g. 
$\colorbox{yellow}{\textcolor{black}{25}}_{28}\textcolor{blue}{\xrightarrow{\textcolor{black}{p=4 \textrm{ or } p+3=7}}}\colorbox{violet}{\textcolor{white}{8}}_{28}$, 
or in a blue purple-to-purple path with path length  
$\textcolor{blue}{\xrightarrow{\textcolor{black}{q=\cdots, 
\textrm{ defining } U=g^q}}}$, e.g.  
$\colorbox{violet}{\textcolor{white}{8}}_{18}
\textcolor{blue}{\xrightarrow{\textcolor{black}{q=3 \textrm{ or } q+3=6}}}
\colorbox{violet}{\textcolor{white}{8}}_{28}$. 
The congruence class graph shows that the Bařina $7n\pm1$ conjecture is proven, since $c=7\cdot1+1=8=2^3$ is a power-of-2, and the number trunk$_q$ of $g$-arrows to reach the trunk number $U(c)=c\cdot2^{\text{trunk}_q}=8\cdot2^{trunk_g}=8\cdot2^3=64\in[8]_{28}$ equals $\min_q=3$ (Eq.\ref{eq:rtt}, Table \ref{fig:TLEE}). 
The $7n\pm1$ tree proves the Bařina $7n\pm1$ conjecture by the root-trunk-test (Eq.\ref{eq:rtt}), as it is a one-trunk-tree with a power-of-2 root $c=7\cdot1+1=2^3$, and a lowest trunk $c^2=8^2$ in period $\ell=0$: $(U(c)=2^{3+(q=3)}=64)=((\ell=0)\cdot196+64)\leq (m=196)$. \quad\textbf{\textbar}
\end{minipage}\end{figure}

The Bařina $7n\pm1$ function (Def.\ref{def:f7pm1}, Fig.\ref{fig:tree7pm1}) has $v=2$ different subfunctions for the odd class $[1]_2$: $n\to7n+1$ for subclass $[1]_4$ and $n\to7n-1$ for subclass $[3]_4$. Its congruence class graph with $a^2v^2=\frac{c}{2}(c-1)^2-196$ classes (Eq.\ref{eq:lcmpq}), of which $a-1\cdot v=12$ are branching classes, is depicted as a congruence classes graph with $m=a^1v^2=28$ classes, of which $12/(a-1)=2$ purple branching classes $\colorbox{violet}{\textcolor{white}{8}}_{28}$  and $\colorbox{violet}{\textcolor{white}{20}}_{28}$ (Fig.\ref{fig:cc7pm1}). Combined branching classes imply combined leftward and upward paths to them, listed in Table \ref{tab:7pm1classes}. Next to $10$ upward paths $c10u3$ from and to the branching nodes $\colorbox{violet}{\textcolor{white}{8,20}}_{28}$ with length $3$, and $2$ upward paths $c2u6$ with length $6$ due to an additional cycle of $\textrm{min}_q=3$ (Eq.\ref{eq:minq}) $g$-arrows. Next to $6$ combinations of $12$ leftward paths of length $2$ ($c12m2$, $c12p2$), $3$ ($c12m3$, $c12p3$) and $4$ ($c12m4$, $c12p4$), $6$ combinations of $2$ leftward paths show up with an additional cycle of $\textrm{min}_q=3$ arrows ($c2m5$, $c2p5$, $c2m6$, $c2p6$, $c2p7$, $c2m7$). The two red leftward classes  $[\textcolor{red}{456},\textcolor{red}{916}]_{a^3v^2=1372}$ are two contracting classes giving repeatedly lower output numbers, which specify,just like contracting Collatz 3$n$+1 class $[52]_54$, the Turing halting proof they converge to the trivial cyclic root rather than to a non-trivial cycle or a divergent trajectory (Section \ref{sec:halting}).

\paragraph{The LCM periodicity of the Leftward and Upward functions}

The $L$eftward and $U$pward functions \(L = f g^p\) (Def.~\ref{def:Lw7}) and \(U = g^q\) (Def.~\ref{def:Uw7}) for the Bařina \(7n\!\pm\!1\) function take as arguments the branching classes \(S_{\ge0}\) sharing powers of \(2\) (\(p\) or \(q\)), with labels assigned in Table~\ref{tab:7pm1classes}.

The $L$eftward function \(L = f g^p\) (Def.~\ref{def:Lw7}) transforms the \(12\) argument classes \([c12p2]_{2^2 7^2}\) with periodicity \(2^2 7^2\) into \(384\) red $L$eftward output $s$-classes \([384\;c12p2]_{2^9 7^2}\) with LCMp (least-common-multiple periodicity) \(2^9 7^2=25\,088\).  
The count \(384 = 12 \cdot (2^9 7^2)/(2^{2+2} 7^2) = 12 \cdot 2^5\) follows from scaling the \(c12p2\) periodicity \(2^2 7^2\) to the LCM-periodicity \(2^9 7^2\); the total count of Leftward numbers per LCMp period over all rows is \(c1372 = c(2^2 7^3) = c(a^3 v^2 = 343 \cdot 4)\) (Eq.~\ref{eq:luclass}, Fig.\ref{fig:hamiltonBa} top row left).

\begin{definition} \label{def:Lw7} \textit{Leftward function} $L:\;S_{\ge 0}\textcolor{red}{\to} s_{0;\ge 1};\;\;n\textcolor{red}{\to} L(n)=2^p\frac{n\pm1}{7}$
\end{definition}
\small{
\begin{tabular}{p{0.55in}p{0.05in}p{0.02in}p{0.85in}p{0.02in}p{0.7in}p{0.1in}p{1.3in}} \hline 
    $S_{\ge 0}$\newline  & \textit{p} & : & $L(n)$  & ; &  $s_{0;\ge 1}$  & $=$ & ${\left[c1372l\right]}_{2^97^2=25\,088} $ \\ \hline 
${\left[c12m2\right]}_{2^27^2}$ & 2 & : & ${fg}^2(n)=2^2\frac{n-1}{7}$ &; &  ${\left[c12m2\right]}_{2^{2+2}7^2}$ & $=$ & ${\left[384\; c12m2\right]}_{2^{9}7^2}$ \\ 
${\left[c2m5\right]}_{2^27^2}$ & 5 & : &${fg}^5(n)=2^5\frac{n-1}{7}$ & ; &  ${\left[c2m5\right]}_{2^{2+5}7^2}$ & $=$ & ${\left[8\; c2m5\right]}_{2^97^2}$  \\ 
${\left[c12p2\right]}_{2^27^2}$ & 2 & : & ${fg}^2(n)=2^2\frac{n+1}{7}$ &; &  ${\left[c12p2\right]}_{2^{2+2}7^2}$ & $=$ & ${\left[384\; c12p2\right]}_{2^{9}7^2}$ \\ 
${\left[c2p5\right]}_{2^27^2}$ & 5 & : &${fg}^5(n)=2^5\frac{n+1}{7}$& ; &  ${\left[c2p5\right]}_{2^{2+5}7^2}$ & $=$ & ${\left[8\; c2p5\right]}_{2^97^2}$  \\   \hdashline
${\left[c12m3\right]}_{2^27^2}$ & 3 & : & ${fg}^3(n)=2^3\frac{n-1}{7}$ &; &  ${\left[c12m3\right]}_{2^{2+3}7^2}$ & $=$ & ${\left[192\; c12m3\right]}_{2^{9}7^2}$ \\ 
${\left[c2m6\right]}_{2^27^2}$ & 6 & : &${fg}^6(n)=2^6\frac{n-1}{7}$ & ; &  ${\left[c2m6\right]}_{2^{2+6}7^2}$ & $=$ & ${\left[4\; c2m6\right]}_{2^97^2}$  \\ 
${\left[c12p3\right]}_{2^27^2}$ & 3 & : & ${fg}^3(n)=2^3\frac{n+1}{7}$ &; &  ${\left[c12p3\right]}_{2^{2+3}7^2}$ & $=$ & ${\left[192\; c12p3\right]}_{2^{9}7^2}$ \\ 
${\left[c2p6\right]}_{2^27^2}$ & 6 & : &${fg}^6(n)=2^6\frac{n+1}{7}$& ; &  ${\left[c2p6\right]}_{2^{2+6}7^2}$ & $=$ & ${\left[4\; c2p6\right]}_{2^97^2}$  \\ 
\hdashline
${\left[c12m4\right]}_{2^27^2}$ & 4 & : & ${fg}^4(n)=2^4\frac{n-1}{7}$ &; &  ${\left[c12m4\right]}_{2^{2+4}7^2}$ & $=$ & ${\left[96\; c12m4\right]}_{2^{9}7^2}$ \\ 
${\left[c2m7\right]}_{2^27^2}$ & 7 & : &${fg}^7(n)=2^7\frac{n-1}{7}$ & ; &  ${\left[c2m7\right]}_{2^{2+7}7^2}$ & $=$ & ${\left[2\; c2m7\right]}_{2^97^2}$  \\ 
${\left[c12p4\right]}_{2^27^2}$ & 4 & : & ${fg}^4(n)=2^4\frac{n+1}{7}$ &; &  ${\left[c12p4\right]}_{2^{2+4}7^2}$ & $=$ & ${\left[96\; c12p5\right]}_{2^{9}7^2}$ \\ 
${\left[c2p7\right]}_{2^27^2}$ & 7 & : &${fg}^5(n)=2^7\frac{n+1}{7}$& ; & ${\left[c2p7\right]}_{2^{2+7}7^2}$ & $=$ & ${\left[2\; c2p7\right]}_{2^97^2}$  \\   \hline
\end{tabular}
}
\normalsize{}

\noindent  The class $c2u6$, indicating a number \textit{c}ount $c2$ of the $u6$ classes $[92]_{196}$ and $[104]_{196}$ for which the power-of-2 $q$ in $U=g^q$ is $U^{2 \min_q=6}$ (Eq.\ref{eq:minq}), shows up in the $7n\pm1$  congruence graph (Fig.\ref{fig:cc7pm1}) as a blue path of length $\min_q=3$ with an additional cycle of $(\min_q=3)$ arrows: for $[92]_{196}$ the additional cycle  $[8]_{196}\textcolor{blue}{\to}[16]_{196}\textcolor{blue}{\to}[4]_{196}\textcolor{blue}{\to}[8]_{196}$. The $U$pward function \(U = g^q\) (Def.~\ref{def:Uw7}) transforms the \(2\) argument classes \([c2u6]_{2^2 7^2}\) with periodicity \(2^2 7^2\) into \(4\) blue $U$pward output $S$-classes \([4\;c2u6]_{2^9 7^2}\) with LCMp \(2^9 7^2=25\,088\). 
The count \(4 = 12 \cdot (2^9 7^2)/(2^{2+2} 7^2) = 12 \cdot 2^5\) follows from scaling the \(c2u6\) periodicity \(2^2 7^2\) to the LCM-periodicity \(2^9 7^2\); the total count of $U$pward numbers per LCMp period over all rows is \(c164 = c(2^2 a(a-1)-1) = c(4\cdot41 = 164)\) (Eq.~\ref{eq:luclass}, Fig.\ref{fig:hamiltonBa} top row right).

\FloatBarrier

\begin{table}
\caption{\textit{c}-sets of $a^2v^2=7^22^2=196$ branching classes by $7n\pm1$ defining $a^3v^2=7^32^2=1372$ subclasses}
\label{tab:7pm1classes}
\begin{align}
&7n\pm1 \textrm{ branching classes } &S_{\ge 0}&=\;\;\left[\;\begin{array}{c} 8,36,76,104,132,176\\20,64,92,120,160,188\end{array}\;\right]_{\,2^27^2=196}  \label{eq:clas7}\\
&\textrm{parents upward classes }&U:S_{\ge0}\to S_{\ge1}&\;:\;[\;c10u3,c2u6\;]_{196} \label{eq:uargclas7}\\
&&c10u3&=\;\;[\;8,36,76,132,176,\;\;20,64,120,160,188\;]_{196}  \notag\\
&&c2u6&=\;\;[\;104,92\;]_{196}  \notag\\
&\textrm{upward classes } &S_{\ge 1}&\,:\;\;  \left[c82u\right]_{196\cdot2^6=12544}\;=\;\left[c164u\right]_{196\cdot2^7=25088} \label{eq:uclas7}\\
&\textrm{parents leftward classes }&L:S_{\ge0}\to s_{0;\ge1}&\;:\;   {\left[ \begin{array}{c}
\mathrm{c12m2},\mathrm{c2m5},\mathrm{c12p2},\mathrm{c2p5},\mathrm{c12m3}, \mathrm{c2m6},\\ 
\mathrm{c12p3},\mathrm{c2p6},\mathrm{c12m4},\mathrm{c2m7},\mathrm{c12p4},\mathrm{c2p7} \end{array}
\right]}_{1372} \label{eq:largclas7}\\
&&c12m2 &=\;\;{\left[ \begin{array}{c}
36,64,232,\textcolor{red}{456},624,652, \\ 
820,848,1016,1044,1212,1240 \end{array}
\right]}_{196\cdot7=1372} \notag\\
&&c2m5&=\;\;{\left[\;260,428\;\right]}_{1372} \notag\\ 
&&c12p2 &=\;\;{\left[ \begin{array}{c}
132,160,328,356,524,552, \\ 
720,748,\textcolor{red}{916},1140,1308,1336  \end{array}
\right]}_{1372} \notag \\
&&c2p5&=\;\;{\left[\;944,1112\;\right]}_{1372} \notag\\
&&c12m3&=\;\;{\left[ \begin{array}{c}
8,92,204,288,400,484, \\ 
596,680,792,876,1184,1268 \end{array}
\right]}_{1372}\notag \\ 
&&c2m6&=\;\;{\left[\;988,1072\;\right]}_{1372}\notag \\ 
&&c12p3&=\;\;{\left[ \begin{array}{c}
104,188,496,580,692,776, \\ 
888,972,1084,1168,1280,1364 \end{array}
\right]}_{1372}\notag \\ 
&&c2p6&=\;\;{\left[\;300,384\;\right]}_{1372}\notag \\ 
&&c12m4&=\;\;{\left[ \begin{array}{c}
120,176,316,372,512,568, \\ 
764,904,960,1100,1156,1296 \end{array}
\right]}_{1372}\notag \\ 
&&c2m7&=\;\;{\left[\;708,1352\;\right]}_{1372}\notag \\ 
&&c12p4&=\;\;  {\left[ \begin{array}{c}
76,216,272,412,468,608, \\ 
804,860,1000,1056,1196,1252 \end{array}
\right]}_{1372}\notag \\ 
&&c2p7&=\;\;{\left[\;20,664\;\right]}_{1372}\notag \\
&\textrm{leftward classes }&s_{0;\ge 1}&\,:\;\;[\;c1372l\;]_{\,196\cdot2^7=25088} \label{eq:lclas7} \\
&&c1372l&=\;\;{\left[\begin{array}{c}[ 8, 20, 36, 76, 92, 104,\\120, 132, 176, 188, 204, 216,\\\textrm{nllc=1348 not listed leftward classes},\\24872, 24884, 24900, 24912, 24956, 24968,\\24984, 24996, 25012, 25052, 25068, 25080 \end{array} \right]}_{25088} \notag
\end{align}
\end{table}
\FloatBarrier

\begin{definition} \label{def:Uw7} \textit{Upward function} $U:\;S_{\ge 0}\textcolor{blue}{\to} S_{\ge 1};\;\;n\textcolor{blue}{\to} U(n)=2^qn$
\end{definition}
\small{
\begin{tabular}{p{0.55in}p{0.05in}p{0.05in}p{0.75in}p{0.02in}p{0.7in}p{0.1in}p{1.3in}} \hline 
    $S_{\ge 0}$\newline  & \textit{q} & : & $U(n)=2^q n$  & ;  & $S_{\ge 1}$  & $=$ & ${\left[c164u\right]}_{2^97^2=25\,088} $ \\ \hline 
${\left[c10u3\right]}_{2^27^2}$ & 3 & : & ${g}^3(n)=2^3n$ &;  & ${\left[c10u3\right]}_{2^{2+3}7^2}$ & $=$ & ${\left[160\;c10u3\right]}_{2^{9}7^2}$ \\ 
${\left[c2u6\right]}_{2^27^2}$ & 6 & : & ${g}^6(n)=2^6n$ & ; & ${\left[c2u6\right]}_{2^{2+6}7^2}$ & $=$ & ${\left[4\;c2u6\right]}_{2^97^2}$  \\ \hline
\end{tabular}
}

\noindent The definitions of the $L$eftward and $U$pward function, giving  leftward subtree and upward subtree generations $s_{0;\ge1}$ and $S_{\ge1}$, imply (by Eqs.\ref{eq:lcmpq},\ref{eq:lcmLU},\ref{eq:bclass},\ref{eq:dSge0},\ref{eq:cdp}) the subtree density proof of the Bařina $7n\pm1$ conjecture (Eq.\ref{eq:sdp7}). 

\begin{align}
d(s_{0;\ge1})+d(S_{\ge1})\quad&\stackrel{!}{=}\quad d(S_{\ge0}) \quad\textrm{ as: }\notag\\
\frac{(a(a-1)-1)v^2=164}{2^va^2\cdot2^{\textrm{maxpq}}=25088} +\frac{a^3v^2=1372}{25088}\quad&\stackrel{!}{=}\quad\frac{(a-1)v=12}{2^va^2=196} \label{eq:sdp7}
\end{align}
The equivalent of the $6x6$ leftward subclass matrix $T$  with periodicity $2\cdot3^3=54$ for the Collatz $3n+1$ conjecture is for the Bařina $7n\pm1$ conjecture a $12x12$ leftward subclass matrix of the classes listed in \ref{tab:7pm1classes},  not depicted here. The equivalent of Table \ref{tab:cdp3p1} summarizing the Lean-certified `cotree density proof' for the Collatz $3n+1$ conjecture is 
Table \ref{tab:cdp7pm1} for the Bařina $7n\pm1$ function.

\begin{table}[H]
\caption{Bařina $7n\pm1$ density of 1 $L$eftward and \(7-1\) $U$pward geometric series of cotree densities}
\label{tab:cdp7pm1}
\setlength{\abovedisplayskip}{0pt}
\setlength{\belowdisplayskip}{0pt}
\begin{align} 
\toprule 
d(t_i)&\textrm{ density red $L$(eftward) cotree} & &\textcolor{white}{=} 
 \nonumber \\ 
d(t_1)&=\frac{196\cdot7^2}{196\cdot2^{6+7\cdot1}}; \quad  
d(t_2)=\frac{196a^2A^+v}{196\cdot2^{6+7\cdot2}}=\frac{196\cdot5586}{196\cdot2^{20}} &
\sum_{i=1}^{\infty} d(t_{i=1,2,3,\dots})
&
=\frac{
    {\frac{196\cdot7^2}{196\cdot2^{6+7\cdot1}}}  
} {1-\frac{A^+v=114}{2^7}}
=\frac{196\cdot7}{196\cdot2^7}
\nonumber  \\ 
&\textrm{Cumulative density leftward }7n\pm1\textrm{ cotrees} & 
\textcolor{white}{sum_{i=1}^{\infty} d(t_{i=1,2,3,\dots})}&=\frac{1372}{25088} \nonumber \\ \midrule \nonumber
d(T_j)&\textrm{ density blue $U$(pward) cotree} & &  \nonumber \\ 
d(T_1)&=\frac{196\cdot8\cdot6}{25088\cdot2^{6}} &  
\sum_{i=0}^{\infty} d(T_{1+6i})&=\frac{12288}{2^{21}-1} \nonumber  \\
d(T_2)&=\frac{196\cdot8\cdot5}{25088\cdot2^{9}} &
\sum_{i=0}^{\infty} d(T_{2+6i})&=\frac{1280}{2^{21}-1}\nonumber \\
d(T_3)&=\frac{196\cdot8\cdot4}{25088\cdot2^{12}} &
\sum_{i=0}^{\infty} d(T_{3+6i})&=\frac{128}{2^{21}-1}\nonumber \\
d(T_4)&=\frac{196\cdot8\cdot3}{25088\cdot2^{15}} &
\sum_{i=0}^{\infty} d(T_{4+6i})&=\frac{12}{2^{21}-1}\nonumber \\
d(T_5)&=\frac{196\cdot8\cdot2}{25088\cdot2^{18}} &
\sum_{i=0}^{\infty} d(T_{5+6i})&=\frac{1}{2^{21}-1}\nonumber \\ 
d(T_6)&=\frac{196\cdot7}{25088\cdot2^{21}} &
\sum_{i=0}^{\infty} d(T_{6+6i})&=\frac{28}{25088\cdot(2^7-1)\cdot337}=\frac{1}{38\,347\,904}\nonumber \\ 
&\textrm{Cumulative density upward }7n\pm1\textrm{ cotrees} & \sum_{j=1}^{\infty} d(T_{j=1,2,3,\cdots})&=\frac{A^-v^2=((a-1)a-1)2^2=164}{25088} \nonumber 
\\ \midrule 
&\textrm{Cumulative density Bařina $7n\pm1$ cotrees}&d(tT) &=\frac{1372+164}{25088}=\frac{12}{196}=d(S_{\ge0}) &\qed \nonumber
\\  &\textrm{Lean certificate: Barina$\_$cotree$\_$density} &  & & 
\\\bottomrule \nonumber
\end{align}
\end{table}
\normalsize{}
\subsubsection{$3n-1$ and $5n+1$: root-successor refutations and root-trunk-refutations}
\label{sec:cc3m1}
\begin{figure}[hbt!]
\caption{The 3\textit{n}-1 congruence classes graph  
${G}_{cc}=({\left[0,\cdots,17\right]}_{18},[f,g])$ }
\label{fig:cc3m1}
\begin{center}
\includegraphics[width=0.7\textwidth]{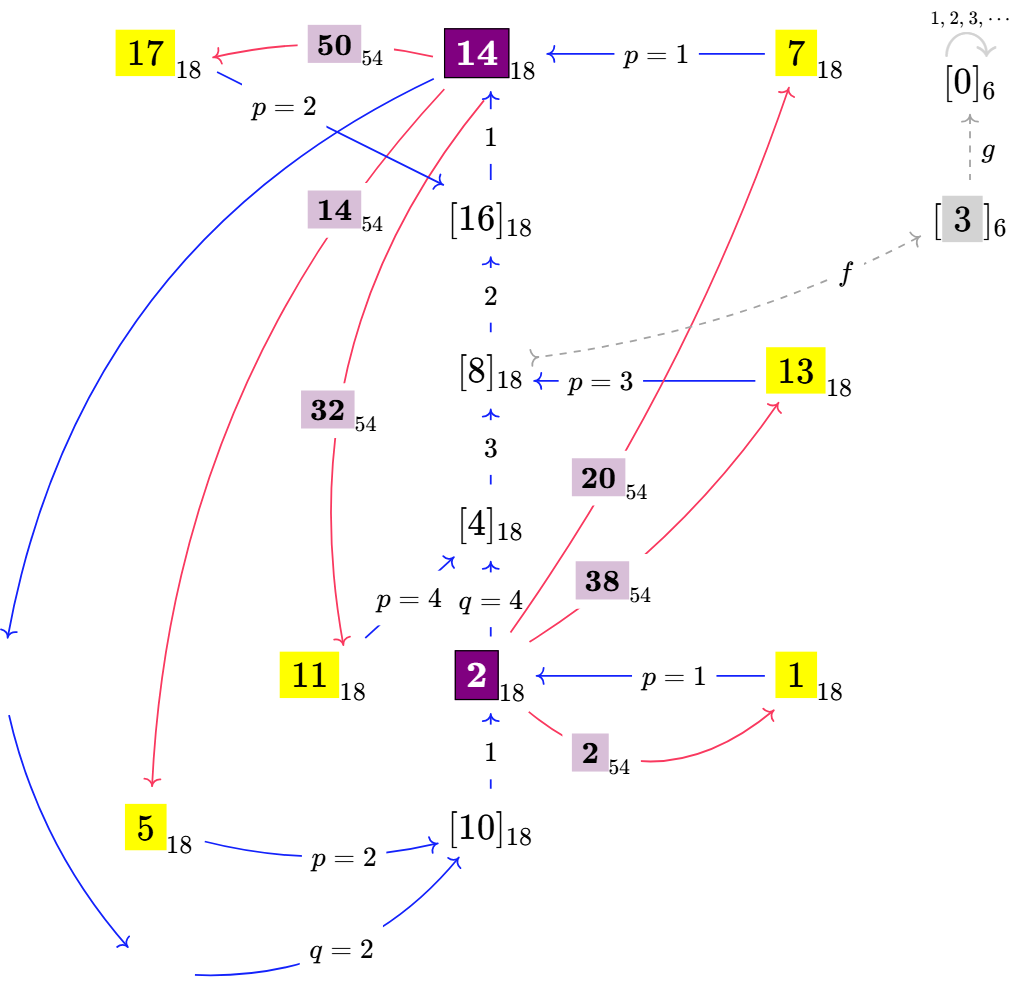}
\end{center}
\centering\begin{minipage}{0.9\linewidth}
\noindent\textbf{Legend: } Nodes: $a-1=2$ even, purple, branching classes with periodicity/modulus $m=2a^2v=18$ (Eq.\ref{eq:lcmpq}), e.g. $\colorbox{violet}{\textcolor{white}{14}}_{18}$, and $\,a(a-1)=6\,$ odd, yellow, classes, e.g.
$\colorbox{yellow}{\textcolor{black}{11}}_{18}$, (Eq.\ref{eq:oclass}) 
from $a(a-1)=6$ even, lightpurple, branching parent subclasses with periodicity $2a^3v=54$ (Eq.\ref{eq:bsclass}).
Red arrow $\textcolor{red}{\xrightarrow{\textcolor{black}{f:\;n\,\to\,(n+1)/3}}}$ (Eqs.\ref{eq:rootPairs},\ref{eq:oclass}) to odd yellow child annotated with its lightpurple (here lightgray) parent subclass, e.g.
$\colorbox{violet}{\textcolor{white}{14}}_{18}
\textcolor{red}{\xrightarrow{\colorbox{lightgray}{\textcolor{black}{[32]}}}}
\colorbox{yellow}{\textcolor{black}{11}}_{18}$.
Blue arrow $\textcolor{blue}{\xrightarrow{\textcolor{black}{g:\;n\,\to\,2n}}}$, either in blue yellow-to-purple path with path length  $\textcolor{blue}{\xrightarrow{\textcolor{black}{p=\cdots, \textrm{ defining } L=fg^p}}}$, e.g. 
$\colorbox{yellow}{\textcolor{black}{11}}_{18}\textcolor{blue}{\xrightarrow{\textcolor{black}{p=3}}}\colorbox{violet}{\textcolor{white}{14}}_{18}$, 
or in blue purple-to-purple path with path length  
$\textcolor{blue}{\xrightarrow{\textcolor{black}{q=\cdots, 
\textrm{ defining } U=g^q}}}$, e.g.  
$\colorbox{violet}{\textcolor{white}{2}}_{18}
\textcolor{blue}{\xrightarrow{\textcolor{black}{q=4}}}
\colorbox{violet}{\textcolor{white}{14}}_{18}$. 
The $3n-1$ tree refutes the $3n-1$ conjecture by the root-trunk-test (Eq.\ref{eq:rtt}), although it is a one-trunk-tree with a power-of-2 root $c=3n-1=2^1$, but it does not have a lowest trunk node $c^2$ in period $\ell=0$:  $(U(c)=2^{1+(q=4)}=32)=((\ell=1)\cdot18+14)\;>\; (m=18)$. $\quad\textbf{\textbar}$
\end{minipage}\end{figure}

The `root successor test' and the `root trunk test' in rows 6 and 7 of Table \ref{tab:roottests} show refutations of the $3n-1$ conjecture (Figs.\ref{fig:tree3m1},\ref{fig:cc3m1}) and of the $5n+1$ conjecture (Fig.\ref{fig:tree5p1}), studied by Kontorovich and Lagarias \cite{RN6}, in column 1 and 2. Columns 3 and 4 show the already proven Collatz $3n+1$ and Barina $7n\pm1$ conjectures.

\begin{figure}[hbt!]
\caption{Proof/refutation of an $an+b$-conjecture: $c=a\cdot1+1$ power-of-$2$? $c^2$ lowest trunk node?}
\label{tab:roottests}
\centering
\small{
\begin{tabular}{p{1.5in}p{0.74in}p{0.9in}p{0.74in}p{1.13in}} \hline 
$an\pm b$   & $3n-1$ & $5n+1$ & $3n+1$ &$7n\pm1$ \\ \hline
$a$ & 3 & 5 & 3 & 7 \\ 
$b$ & $-1$ & $+1$ & $+1$ & $\pm 1$ \\ 
$v=\log_2(c)-1$ \small{subfunctions}\newline{\textcolor{white}{aaa} odd argument (sub)classes} & 1\newline \footnotesize{$[1]_2$} & 1\newline \footnotesize{$[1]_2$} & 1\newline \footnotesize{$[1]_2$} & 2 \newline \footnotesize{$[1,3]_4$} \\ 
\hline
1.\;\;$c=a\cdot1+b$ cyclic root &$2=2^1$&$6\neq2^q$&$4=2^2$&$8=2^3$
\\
\;\;\;\;$U(c)=2^qc$ lowest trunk &$32=2^5>c^2$&$96>c^2$&$16=c^2$&$64=c^2$
\\
2.\;\;$(a-1)v$ odd classes
\newline indivisible by $a$; period $2av$,
\newline per.density $((a-1)v)/2av$ & ${\left[1,5\right]}_6$  \newline \newline $2/6$ & ${\left[1,3,7,9\right]}_{10}$   \newline \newline $4/10$ & ${\left[1,5\right]}_6$   \newline \newline $2/6$ & ${\left[ \begin{array}{c}
1,3,5,9, \\ 
11,13,15,17, \\ 
19,23,25,27 \end{array}
\right]}_{28}$ \newline $12/28$ \newline \\
3.\;\;\small{$(a-1)v$ classes $S_{\ge 0}$, \newline periodicity $2a^2v$,  }
& \footnotesize{${\left[2,14\right]}_{18}$ }
& \footnotesize{${\left[6,16,36,46\right]}_{50}$ }
& \footnotesize{${\left[4,16\right]}_{18}$ }
& \footnotesize{
${\left[ \begin{array}{c}
\mathrm{8,\ 36,\ 64,\ 92,\ } \\ 
\mathrm{120,\ 176,20,76,} \\ 
104,132,160,188 \end{array}
\right]}_{196}$ }\\ 
\small{periodic density  $d(S_{\ge 0})=$}
\newline \small{$((a-1)v)\;/\;(2a^2v)$}\; \newline 
& $2/18$ & $4/50$ & $2/18$ & $12/196$ \\ 
4.\;\;\small{$a(a-1)v$ subclasses $S_{\ge0}$, \newline periodicity $2a^3v$, \newline per.density $(a(a-1)v)/2a^3v$ }
& \footnotesize{${\left[ \begin{array}{c}
2,14,20, \\ 
32,38,50 \end{array}
\right]}_{54}$ } \newline \small{$6/54$}
& \footnotesize{ ${\left[ \begin{array}{c}
6,16,36,\dots  \\ 
\dots ,236,246 \end{array}
\right]}_{250}$ } \newline \small{$20/250$}
& \footnotesize{${\left[ \begin{array}{c}
4,16,22, \\ 
34,40,52 \end{array}
\right]}_{54}$} \newline \newline \small{$6/54$}& 
\footnotesize{${\left[ \begin{array}{c}
\mathrm{8,36,64,\dots } \\ 
\dots ,1352,1364 \end{array}
\right]}_{1372}$} \newline \newline \small{$84/1372$} \newline \\ 
5.\;\;\small{$a$ root class subclasses,}\newline\vspace{1mm}  $a-1$ subclasses next to $[c]_{2a^3v}$ 
& \footnotesize{${\left[2;\;\,20,38\right]}_{54}$} 
& \footnotesize{${\left[ \begin{array}{c}
6;\;\,56,106, 
\\ 156,206 \end{array}
\right]}_{250}$ }
& \footnotesize{${\left[4;\;\,22,40\right]}_{54}$} 
& \footnotesize{$\left[ \begin{array}{c}
8;\;\,204,400,596, 
\\ 
792,988,1184 \end{array}
\right]_{1372}$ }\\ \hline
6.\;\;\small{pass root successor test \newline on $a-1$ root subclass numbers
\newline  {$\quad c\xrightarrow{\text{?}} c\;<\;[c]_{2a^2v}\;<\;2a^3v$} 
\newline\newline\newline\newline\small{or on odd }(n-b)/a} &  \vspace{0.002mm} no: \newline \footnotesize{20, 38} \newline\newline\newline\newline\newline\footnotesize{7, 13} & 
\vspace{0.002mm} no: \newline \footnotesize{56, 106, 156, 206} 
\newline\newline\newline\newline\newline\footnotesize{11, 21, 31, 41}& 
\vspace{0.002mm}yes: \newline\small{
$\sm{10^4\\22\;\,}$\;
$\sm{10^1\\40\;\,}$} 
\newline\newline\newline\newline\newline\footnotesize{7, 13}& 
\vspace{0.002mm}yes:  \newline\small{
$\quad\sm{10^210^7\\\quad\;\,204}$}\quad\quad
$\sm{10^210^{12}\\400}$\newline\small{
$\sm{10^210^{12}10^3\\596}$\;
$\sm{10^210^{8}10^{4}\\792\quad}$}\newline\small{
$\sm{10^110^{3}10^1\\\;988}$\;
$\sm{\;\;10^210^{15}\\\;\;\;1184}$}\newline\newline\footnotesize{29,57,85,113,141,169}\vspace{1mm}
\\ \hline 
7.\;\;\small{pass root-trunk test}\newline on root $c$ and trunk $U(c)$\newline $c \sm{?\\\textrm{ power-of-2}},$\;$U(c) \stackrel{?}{\in} (\ell=0) \newline \textrm{\textcolor{white}{blablablabl}} c^2 \textrm{ lowest trunk node?}$ 
& \footnotesize{$c=2\sm{\textrm{yes,}\\\textrm{ power-of-2}}\newline\sm{U(c)=32 \stackrel{\textrm{no,}}{\;\notin}(\ell=0)\\ 32=1\cdot18+14;\;32>18}$} \newline\newline\textrm{no, refuted}
& \footnotesize{$c=6 \sm{\textrm{no, a shrub,}\\\textrm{not-power-of-2}}\newline\sm{U(c)=96\stackrel{\textrm{no,}}{\;\notin}(\ell=0)\; \\ 96=1\cdot50+46;\; 96>50}$ } \newline\newline\textrm{no, refuted}
& \footnotesize{$c=4\sm{\textrm{yes,}\\\textrm{ power-of-2}};\newline\sm{U(c)=16 \stackrel{\;\textrm{yes,}}{\;\in}(\ell=0)\\ 16=0\cdot18+16;\;16<18}$}  \newline\newline\textrm{yes, proven}
& \footnotesize{$c=8\sm{\textrm{yes,}\\\textrm{ power-of-2}}\newline\sm{U(c)=64 \stackrel{\;\textrm{yes,}}{\;\in}(\ell=0)\\ 64=0\cdot196+64;\;64<196}$ }\newline\newline\textrm{yes, proven}
\\\hline
\end{tabular}
}
\centering\begin{minipage}{0.9\linewidth}
\noindent\textbf{Legend: }
The top row lists the $a$n+$b$ parameters $a$, $b$ and the number $v=\log_2(c)-1$ of subfuctions generating the a$n$+b tree from odd congruence (sub)classes. 
Row 1 lists the cyclic roots $c=a\cdot1+1$ (Eq.\ref{eq:Croot}) determining $v$, and whether the lowest trunk number $U(c)$ (Eq.\ref{eq:Trunk}) is equal to $U(c)=c^2$ or exceeds $U(c)>c^2$.
Row 2 lists the yellow odd congruence classes $n_o$ (Eq.\ref{eq:oclass})
Row 3 lists the purple, red or blue, congruence classes $n_b$ (Eq.\ref{eq:bclass})
Row 4 lists the light-purple, congruence subclasses of purple classes $n_bs$ (Eq.\ref{eq:bclass})
Row 5 lists the period $\ell=0$ numbers in the light-purple root class subclasses with lcm$_{\textup{p(eriodicity)}}=2a^3v$.
Row 6 lists the root successor tests (Eq.\ref{eq:rtt}) whether the $a-1$ root class subclass numbers next to $[c]_{2a^3v}$ in row 5, are root successor numbers. These are the $a-1$ subclass numbers in period $\ell=0$, of subclass numbers $(\ell=0)(m=2a^3v) +\sum_{a=1}^{a-1}{2a^2v\cdot i}$, which are $a-1$ numbers $c<[c]_{2a^2v}<2a^3v$ higher than the leftward trivial root number $c$ and lower than the upper boundary, the lcm$_{\textup{p(eriodicity)}}=2a^3v$ of light-purple numbers.
Row 7 lists the root-trunk tests (Eq.\ref{eq:rtt}) whether the trivial root $c=a\cdot1+b$ is a power-of-2, and whether the Cayley tree is a one-trunk tree with lowest upward trunk node $(U(c)=c^2\in\ell_0)<(m=2a^2v)$.\quad\textbf{\textbar}
\end{minipage}\end{figure}

\begin{theorem}
\label{the:decidab}
A Collatz-Conway $an+b$-conjecture is refuted ($T_{\ge0}\subset G_{\ge0}$ rather than $T_{\ge0}=G_{\ge0}\,$) 
if the trivial cyclic root $c=a\cdot1+1$ is not a power-of-$2$ of a one-trunk tree with lowest upward trunk number $U(c)=c^2$ (Table \ref{tab:roottests} row 2), as shown by the: 
\begin{itemize}
\item `root successor' test: period $\ell=0$ branching numbers $c<n\in[c]_{2a^2}\stackrel{?}{<}2a^3$ (Table \ref{tab:roottests} row 6)
\item `root-trunk' test: period $\ell=0$ trunk number $U(c)\stackrel{?}{=}c^2$ (Table \ref{tab:roottests} row 7),
\end{itemize}
\end{theorem}
\noindent \textit{Proof.} 
If the trivial root $c=a\cdot1+1$ is not a power-of-$2$, then at least one power-of-2, labeled $\gamma$, is not rooted in $c$. This applies for example to $5n+1$ (Fig.\ref{fig:tree5p1}), for which $c=5\cdot1+1=6$,  but $\gamma=2^4=16$ reduces to $3<c$ by $g^{-1}:\gamma\to(\gamma-1)/5=3$, not rooted by $5n+1$ in $c=6$. 

If the lowest trunk number $U(c)$ is not $c^2$, then a lower squared-power-of-2, labeled $2^{\gamma}$ is not rooted in $c$. This applies to $3n-1$ (Figs.\ref{fig:tree3m1},\ref{fig:cc3m1}), for which $c=3\cdot1-1=2$ and $U(c)=2\cdot2^4=2^5=32$, but $2^{\gamma}=2^3=8$ reduces to $3$ by $g^{-1}:2^{\gamma}\to(2^{\gamma}+1)=9/3=3$, not rooted by $3n-1$ in $c=2$. 

\noindent The `root successor test' and the `root trunk test' in rows 6 and 7 of Table \ref{tab:roottests} show that a$n$+b conjectures 3$n$-1 and 5$n$+1 in columns 1 and 2 with a cyclic root $c=a\cdot1+1$ that is not a power of 2 and/or a lowest trunk number $U(c)>c2$ are refuted.\qed

\subsubsection{Infinite Turing/Leibniz tape $an+b$ conjectures:  Hamiltonian tour  proofs and root-trunk-proofs}

Combinations of equations from subsections \ref{sec:31eqs}-
\ref{sec:35eqs} in Tables \ref{tab:roottests}, \ref{fig:Gccc} specify cyclic root $c$-expressible proofs of Collatz-Conway a$n$+b conjectures on the infinite Turing/Leibniz tape $T_L$ (Table \ref{fig:TLEE}, Theorem {\ref{the:TLEE}}). 

Following the $\eta_0$ Conway 1$n$+1 base tree, Table \ref{fig:TLEE} shows in its columns $\eta=0,1,2,\cdots$ four Collatz-Conway a$n$+b trees : $\eta_1$ Collatz $3n+1$, $\eta_2$ Bařina $7n\pm1$, $\eta_3$ Euclid/Euler, and $\eta_5$ Fermat/Gauss. The name Euclid/Euler reflects the observation of Euclid of (Mersenne) prime $31$ paired with perfect number $31\cdot2^4=496$, generalized by Euler, who calculated by hand that $2^{31}-1$ is a Mersenne prime number. The name Fermat/Gauss reflects the compass and straightedge depiction by Gauss (1777-1855) of the regular $17-$polygon, while $17$ is also a Fermat (1607-1685) prime $2^{2^p}+1$.

The rows in Table \ref{fig:TLEE} show parameters and proofs of proven/provable a$n$+b conjectures paired with, expressible in, their power-of-2 cyclic root $c=2^{2^{\eta-1}+1}$ cyclic root $c=a\cdot1+1$, with as penultimate row the cotree-density `Hamiltonian tour' test (Figs.\ref{fig:hamiltonCo},\ref{fig:hamiltonBa},Table \ref{fig:Gccc}) showing cotree-density proofs (Figs.\ref{fig:hamiltonCo},\ref{fig:hamiltonBa},Table \ref{fig:Gccc}) and as last row the root-trunk test (Theorem \ref{the:decidab}, Table \ref{fig:Gccc}) showing root-trunk proofs.

\begin{theorem}
\label{the:TLEE} All $\eta=1,2,3,\dots=\mathbb{N}_1$ Collatz-Conway a$n$+b functions specifying Cayley a$n$+b successor trees $T_{\ge0}$ of the Conway$_{\eta0}$  1$n$+1 tree with cyclic root $c=1\cdot1+1={2^{a-1}+1}={2^{0}+1}=2^1$ (Fig.\ref{fig:treeP1}) on the infinite Turing-Leibniz tape $T_L$ (Table \ref{fig:TLEE}) with Turing-computable cyclic roots $c=a\cdot1+1=2^{2^{\eta-1}+1}$ of proven a$n$+b conjectures---for which the a$n$+b function generated graph $G_{\ge0}=T_{\ge0}$ is the Cayley tree  rather than a superset graph $G_{\ge0}>T_{\ge0}$ (Eq.\ref{eq:tleg})---by, among others, a $c$-expressible:
\begin{itemize} 
\item periods $\ell=0,1,2,\dots$ cotree-density proof of a `Hamiltonian tour' from and back to $c$ (Eq.\ref{eq:dtT})
\item period $\ell=0$ root-trunk proof of a lowest trunk node $c^2$ (Eq.\ref{eq:rtt}).
\end{itemize}
\end{theorem}
\noindent\textit{Proof.} The Turing-Leibniz tape $T_L$ (Tab.\ref{fig:TLEE}) with countably infinite column numbers $\eta=0,1,2,3,4,\dots=\mathbb{N}_0$ paired in row 2 with Turing-computable cyclic roots $c=a\cdot1+1=2^{2^{\eta-1}+1}$  
\begin{itemize} 
\item of Cayley $E=V=F=\mathbb{N}_0$ a$n$+b trees $T_{\ge0}=G_{\ge0}$ (Figs.\ref{fig:tge0},\ref{fig:treeP1},\ref{fig:tree7pm1}) rather than superset graphs $G_{\ge0}\supset T_{\ge0}$ (Figs.\ref{fig:tree3m1},\ref{fig:tree5p1}), as described in Theorem \ref{the:caytree} equations \ref{eq:tleg}-\ref{eq:rtt}
\item 
with unique pairs $\sm{y\\n}$ of $\quad\quad a\nmid n\quad\,=$\;an$_{\textup{odd}}$+b numbers and \newline binary root path numbers 
$\;\;\;y_{\ge0}\in\mathbb{N}_0=\croot\tcr{0}_{s0}\tcb{\to1}_{S1}\;\tok\;1\tcr{0}_{s1}\tcb{\to}1\tcb{1}_{S2}\;\tok\;10\tcr{0}_{s2},\dots$ (Eq.\ref{eq:sorS})
with same-colour red or blue arrows to them in red $L$eftward $s_{\ge0}$-generations in $t_{\ge1}$ cotrees or blue $U$pward $S_{\ge1}$-generations in $T_{\ge1}$ cotrees (Eq.\ref{eq:sorS})
\end{itemize}
paired in row 1 with $a=c-1=a\cdot1+1=2^{2^{\eta-1}}$, in rows 3-9 with $c$-expressible $T_{\ge0}$ features, and in rows 10-11 with the $c$-expressible Lean-certified `cotree density' proof of a Hamiltonian tour from and to $c$ and the `root-trunk' proof of a lowest trunk node $c^2$. Longer row legends in `Legend paired rows tape $T_L$'.  

\newpage
\begin{table}[hbt!]
\caption{Infinite Turing-Leibniz tape $T_L$: columns a$n$+b conjectures, rows $c$-expressible Cayley $T_{\ge0}$ features}
\label{fig:TLEE}
{\small
\[
\setlength{\arraycolsep}{2pt}
\renewcommand{\arraystretch}{1.0}
% \begin{array}{@{}|lccccccc@{}}
\begin{array}{@{}l!{\vrule width .4pt}ccccccc@{}}
\hline\textrm{Cayley tree } T_{\ge0}\quad\quad\quad\quad \text{\, }   
& \text{Conway} & \text{Collatz} & \text{Ba}\check{r}\text{ina} & \text{Euclid/Euler}       
& \cdots & \textrm{Fermat/Gauss} & \cdots 
\\\quad a\textit{n}+b      
& 1n+1 & 3n+1 & 7n\pm1 & 31n\pm b_1b_2
&  & 131\,071n\pm b_1\cdots b_4 &
\\\quad\eta\,\textrm{, Turing-Leibniz tape } T_L\textrm{ column}      
& 0 & 1 & 2 & 3
&\dots  & 5 &\dots
\\
\hline a=c-1  
& 1 & 3=2^2-1 & 7=2^3-1 & 31=2^5-1
& \dots & 131\,071=2^{17}-1 & \dots 
\\
c=a\cdot1+1\text{, cyclic root} &  &  &  & 
& \dots & 2^{17} & \dots 
\\
\quad c_{\eta0}=2^{\eta}+1;\quad\;\; c_{\eta\ge1} = 2^{2^{\eta-1}+1} & 2^1 & 2^2 & 2^3 & 2^5
& \dots & 2^{17} & \dots 
\\
\psi=c-2\text{, \;\#}\,b \text{, \;\# upward cotree series} & 0 & 2 & 2\cdot3 & 2\cdot3\cdot5
& \dots &  2\cdot3\cdot5\cdot17\cdot257
 & \dots 
\\
U(c)=2^{2(v+1)}=c^2\text{, lowest trunk node} & 2^2=4& 2^{2\cdot2}=16 & 2^{3\cdot2}=64 & 2^{5\cdot2}=1\,024   
& \dots & 2^{17\cdot2=34} & \dots
\\
    v\;=\log_2(c)-1\text{, \# odd subfunctions} & 0 & 1=2^0 & 2=2^1 & 4=2^2
        & \dots & 16=2^4 & \dots 
\\
    \quad k=1\cdots v \textrm{, yellow } [1^k]_{{10}^k} \textrm{ setting } b_k &  
    & \sm{[1^1]_{10^1}\\+1\quad} & 
    \sm{[1^1,1^2]_{10^2}\\+1,-1\quad} & 
    \sm{[1^1,1^2,1^3,1^4]_{10^4}\\+1,-1,+?,-?\quad}    & \dots & 
    \sm{[1^1,1^2,\dots,1^{16}]_{10^{16}}\\+1,-1,\dots,-?\quad\;\;}      & \dots 
\\\eta_{\ge1}\in\;\mathbb{N}_1=\log_2(\log_2(c)-1)+1 &  & 1 & 2 & 3
        & \dots & 5 & \dots 
    \\
    \text{lcm}_{o}=\frac{c}{2}\;(c-1)\text{, perfect yellow lcmp}   &  & 2^1(2^2-1)=6 & 2^2(2^3-1)=28 & 2^4(2^5-1)
        & \dots & 2^{16}(2^{17}-1)\;\;\; & \dots 
    \\
    \text{lcm}_{pq }\;=\frac{c}{2}\;(c-1)^2\text{, purple lcmp} &  & 2^1(2^2-1)^2=18 & 2^2(2^3-1)^2=196 & 2^4(2^5-1)^2
        & \dots & 2^{16}(2^{17}-1)^2& \dots
    \\\text{lcm}_{LU}=(\log_2(c)-1)\;\frac{c^3}{2}\; (c-1)^2,&  & \;\;2^5(2^2-1)^2= & \;\;2^9(2^3-1)^2= & 2^{16}(2^5-1)^2
        & \dots & 2^{54}(2^{17}-1)^2& \dots
    \\\quad\quad\quad\;\text{ red blue lcmp} &  & 288 & 25\,088 & 
        &  & &     
    \\\hline
        \text{cotree density test, Hamiltonian tour} &  &  &  & 
        &  & & 
    \\  \quad
        d(tT)\stackrel{?}{=}d(S_{\ge0})\quad\quad & & &  &\frac{2^4(2^5-1)^3+2^4929^1}{2^{16}(2^5-1)^2} & & \frac{2^8(2^{17}-1)^3+\dots}{2^{29}(2^9-1)}^2
    \\  \quad\quad\;\;
        \frac{\;a^3v^2+A^-v^2\;}{\text{lcm}_{LU}}\stackrel{?}{=}\frac{(a-1)v}{\text{lcm}_{pq}}\quad\quad & & \frac{27+5}{288}\stackrel{!}{=}\frac{2}{18} &  \frac{1372+164}{25088}\stackrel{!}{=}\frac{12}{196} &\stackrel{!}{=}\frac{120}{2^4(2^5-1)^2}
        &\dots&  \stackrel{!}{=}\frac{2\,097\,120}{2^{16}(2^{17}-1)^2}& \dots
    \\  \textrm{Lean-certif: infty$\_$Ham$\_$tour$\_$proofs}
   \\\hline
        \text{root-trunk test }\stackrel{?}{\Rightarrow} \text{root-trunk proof} &  &  &  & 
        &  & & 
    \\\quad  (U(c)\stackrel{?}{<}\text{lcm}_{pq})\stackrel{!}{\Rightarrow}
    (c^2\stackrel{!}{<}\frac{c}{2}\;(c-1)^2)
    & & 16\stackrel{!}{<}18 &  64\stackrel{!}{<}196 & 2^{10}\stackrel{!}{<}2^4(2^5-1)^2
        &\dots& 2^{34}\stackrel{!}{<}2^{16}(2^{17}-1)^2 & \dots
    \\\textrm{Lean-certif: infty$\_$root$\_$trunk$\_$proofs}
    \\\hline
\end{array}
\]}
\end{table}
\textbf{Legend paired rows tape $T_L$}: 
\begin{enumerate}
\item $a=c-1$, short legs Pythagorean triangles for $a_{\eta1}=3$, subsequent short legs $a_{\eta\ge2} = 2^{2^{\eta-1}}$ by sum of legs  previous triangle (OEIS A365577), $a_{\eta1}-a_{\eta3}$ Mersenne primes $2^p-1$  
\item $c=a\cdot1+1$, Cayley tree $T_{\ge0}$ Turing-computable cycic roots  $c_{\eta0}=2^{\eta}+1=3;\quad c_{\eta\ge1} = 2^{2^{\eta-1}+1}$ 
\item $\psi=c-2=a-1$, \# number of $b$-parameters, \# number of  
$U$pward cotree paths in Aut($T_{\ge0},[L,U]$),
$\psi_{\eta1}-\psi_{\eta..??}$ primes,
in $\psi_{\eta2}-\psi_{\eta6}$ products 
of $(\psi_{\eta1}=2)$ and 5 successive Fermat primes $2^{2^p}+1$ 
\item $U(c)=2^{2(v+1)}=c^2$, Cayley tree $T_{\ge0}$ lowest trunk node
\item $v\;=\log_2(c)-1$, \# odd subfunctions setting $b$ for $v$ odd subclasses $[b_1,\dots,b_v]_{10^v}$ of odd class $[1]_{10^1}$
\item $\eta_{\ge1}\in\;\mathbb{N}_1=\log_2(\log_2(c)-1)+1$, natural numbers $1,2,3,\dots$ paired in $T_L$-columns with computable cyclic roots $c$ and $c$-expressible features of proven Collatz-Conway a$n$+b-conjectures   
\item lcm$_{o}=\frac{c}{2}(c-1)$, periodicity yellow $a\nmid$ n$_{\textup{odd}}$ numbers,  lcm$o_{\eta1}..o_{\eta3}$ perfect numbers $6$, $28$, $2^4(2^5$-$1)$=$496$
\item lcm$_{pq}=\frac{c}{2}(c-1)^2$, lcm($p$eriodicity) purple, either red or blue, $a\nmid$ an$_{\textup{odd}}$+b, numbers
\item lcm$_{LU}=(\log_2(c)-1)\;\frac{c^3}{2}$, lcm($p$eriodicity) red and blue, $a\nmid$ an$_{\textup{odd}}$+b, Cayley tree $T_{\ge0}$ numbers
\item  $d(tT)\stackrel{?}{=}d(S_{\ge0})$, `cotree density' test Hamiltonian tour $\frac{\;a^3v^2+A^-v^2\;}{\text{lcm}_{LU}}\stackrel{?}{=}\frac{(a-1)v}{\text{lcm}_{pq}}\newline
\quad\quad\Rightarrow$ `cotree density' proof Hamiltonian tour
        $\frac{2(c-2)\log_2(\frac{c}{2})}{(c-1)^2\,c\; }\stackrel{!}{=}
        \frac{2(c-2)(\log_2(c)-1)}{(c-1)^2\, c }$
\item  $U(c)\stackrel{?}{<}\text{lcm}_{pq}{\quad\Rightarrow\quad}$`root-trunk' proof
    $\;\;c^2\stackrel{!}{<}\frac{c}{2}\;(c-1)^2\;\;$ for $c>2+\sqrt3$, for integer $c\ge4$.\qed 
\end{enumerate}

\section{Supplementary proofs}
\label{sec:supplement}
Common proof approaches to Collatz-Conway a$n$+b conjectures, aimed at proving or refuting the Collatz 3$n$+1 conjecture by studying whether the root paths of all separate natural numbers converge to the trivial root $\cabroot$ \cite{RN3} were hitherto unsuccessful, or delivered heuristic proofs. The supplementary proofs below complete three of these, incomplete or heuristic, trivial root convergence proofs: the staircase term rewriting approach, the halting proof approach, and the 'almost all' proof on odd numbers in the Syracuse tree (Lemmas \ref{lem:termrewrite},\ref{lem:halting},\ref{lem:syrproof}). Proof completion is achieved by designating for a full proof (1) a subset of the natural numbers and/or (2) an improved proof delivering equation, using elementary functions specifying Cayley tree $T_{\ge0}$ of a Collatz-Conway a$n$+b conjecture (Theorem \ref{the:caytree}), most notably the paired inverse depth-first functions $\sm{\textcolor{blue}{1^{-j}}\;\textcolor{red}{0}^{-1}\\U^{-j}L^{-1}}:$ and $\sm{\textcolor{red}{0^{-i}}\;\textcolor{blue}{1}^{-1}\\L^{-i}U^{-1}}:\;$ (Eqs.\ref{eq:IdfSj},\ref{eq:Idfsi}).

\subsection{Staircase term rewriting of root successors to (dis)prove an a\textit{n}+b conjecture} 
\label{sec:termrewrites4}

The heuristic staircase \textit{term rewriting} approach \cite{german,klop} offers a fast check of convergence for arbitrarily selected numbers to the trivial root $\cabroot$. 
\begin{lemma}
\label{lem:termrewrite}
The staircase term rewriting heuristic for a selected Collatz-Conway a\textit{n}+b conjecture that the trivial root is ‘always reached, in all probability!’  \cite{klop} is a fast supplementary proof ($T_{\ge0} = G_{\ge0}$) or refutation ($T_{\ge0} \subset G_{\ge0}$) on the finite subset of an+b numbers required for the root-successor test (Theorem \ref{the:decidab}) using the paired inverse depth-first functions $\sm{\textcolor{blue}{1^{-j}}\;\textcolor{red}{0}^{-1}\\U^{-j}L^{-1}}:$ and $\sm{\textcolor{red}{0^{-i}}\;\textcolor{blue}{1}^{-1}\\L^{-i}U^{-1}}:\;$ (Eqs.\ref{eq:IdfSj},\ref{eq:Idfsi}).
\end{lemma}

\noindent\textit{Proof.} Table \ref{tab:klop31} illustrates the fast staircase jump approach for the Collatz 3$n$+1 function and the even number $94 = 3 \cdot 31 + 1$, which is the lowest even number not reached by $20$ breadth-first iterations of the Collatz $3n + 1$ function (Fig. \ref{fig:tge0}b), of which the odd number $f:94\tcr{\to}(94-1)/3=31$ exemplified the heuristic Syracuse staircase proof \cite{klop}. Odd numbers $n=(n_\textup{even}-b)/a$ in Syracuse root paths are written as binary numbers, of which the zeros and ones are pictured in the reverse order as black squares $0:\blacksquare$ and white squares $1:\square$ respectively. The root path to $31$ in the Syracuse tree becomes a Syracuse staircase of $39$ steps with as lowest upwards step the step from $1:\square$ to $5:\square\blacksquare\square$ and as $39^{th}$ step the step from $47:\square\square\square\square\blacksquare\square$ to $31:\square\square\square\square\square$. 

Instead of reaching the trivial root $\cabroot$ by $39$ Syracuse step rewrites, Table \ref{tab:klop31} shows it  is reached by $4$ staircase jumps defined by the paired inverse depth-first function $\sm{\textcolor{red}{0^{-i}}\;\textcolor{blue}{1}^{-1}\\L^{-i}U^{-1}}:\;$ from $L$eftward cotree $t_{26}$, via cotrees $t_4$, $t_7$, $t_1$ to $\cabroot$.  Jumps by $\sm{\textcolor{blue}{1^{-j}}\;\textcolor{red}{0}^{-1}\\U^{-j}L^{-1}}:$ are impossible.

The staircase term rewriting approach in the Cayley tree $T_{\ge0}$ of even numbers $a\nmid n=a n_\textup{odd}+b$, e.g. of Collatz 3$n$+1 number $3\cdot31+1=94$, removes by functions $\sm{\textcolor{red}{0}^{-i}_{s_i}1^{-1}\\L^{-i}U^{-1}(n)}:$ (Eq.\ref{eq:Idfsi}) and $\sm{\textcolor{blue}{1}^j_{S_j}0^{-1}\\U^{-i}L^{-1}(n)}:$ (Eq.\ref{eq:IdfSj}) trailing red (leftward) and blue (upward) same-cotree number pairs, including same-cotree cyclic root number pairs on the greedy branch and upward trunk. Table \ref{tab:klop31} shows the $4$ rewrites by which Eq.\ref{eq:Idfsi} transfers red, $L$eftward, number pairs $\sm{y\\n}$ in cotrees $t_{26}$, $t_4$, $t_7$ and $t_1$. The $4$ staircase jumps correspond to the $26+4+7+(1(U(c)\tcb{\to}c)+1(c\tcr{\to}1)) = 39$ Syracuse staircase steps from $94$ down to $\cabroot$. 
\vspace{-0.3cm}
\begin{table}[h!]
\centering
\begin{tabular}{@{}lccccc@{}}

\parbox[t]{0.40\textwidth}{{\small\textrm{Term rewrite, cotree switch}}}
&
\makebox[0.045\textwidth][c]{0}
&
\makebox[0.045\textwidth][c]{1}
&
\makebox[0.055\textwidth][c]{2}
&
\makebox[0.065\textwidth][c]{3}
&
\makebox[0.075\textwidth][c]{4}
\\

\parbox[t]{0.30\textwidth}{\small{\textrm{cotree}}}
&
&
\makebox[0.045\textwidth][c]{$t_1$}
&
\makebox[0.045\textwidth][c]{$t_7$}
&
\makebox[0.055\textwidth][c]{$t_4$}
&
\makebox[0.085\textwidth][c]{$t_{26}$}
\\

\parbox[t]{0.30\textwidth}{
$\sm{y\quad\;\;\textup{binary root path number }\quad\;\;\\n\quad\;\;
a\,\,\nmid\,\,a n_{\textup{odd}}+b
\quad\;\;\;\text{}\;\;\quad\quad\quad}$
}
&
\makebox[0.045\textwidth][c]{$\sm{\tcr{0}\\4}$}
&
\makebox[0.045\textwidth][c]{$\sm{1\tcr{0}\\40}$}
&
\makebox[0.055\textwidth][c]{$\sm{101\tcr{0}^7\\2308\;}$}
&
\makebox[0.065\textwidth][c]{$\sm{1010^71\tcr{0}^4\\1822}$}
&
\makebox[0.13\textwidth][c]{$\sm{1010^{7}10^{4}1\tcr{0}^{26}\\94}$}
\\

\parbox[t]{0.40\textwidth}{\small{\textrm{cotree root on greedy branch}}}
&
\makebox[0.045\textwidth][c]{$\sm{\tcr{0}\\4}$}
&
\makebox[0.045\textwidth][c]{$\sm{1\tcr{0}\\40}$}
&
\makebox[0.055\textwidth][c]{$\sm{1\tcr{0}^7\\196}$}
&
\makebox[0.065\textwidth][c]{$\sm{1\tcr{0}^4\\22}$}
&
\makebox[0.13\textwidth][c]{$\sm{1\tcr{0}^{26}\\23\,293\,636}$}

\end{tabular}

\caption{$4$ staircase jump (red-red) term rewrites by Eq.~\ref{eq:Idfsi}
confirm root path of $3\cdot31+1=94$}
\label{tab:klop31}
\end{table}
\vspace{-0.3cm}
\noindent Staircase jump term rewrites by Eqs.\ref{eq:IdfSj},\ref{eq:Idfsi} for period $\ell=0$ root successor numbers $a$$\nmid$$a$
n$_{\textup{odd}}$+b $\stackrel{?}{<}$ $2a^3v$ provide a supplementary `root-successor' proof or refutation (Eq.\ref{eq:rst}), equivalent to the `root-successor' tests for the conjectures Collatz 3$n$+1 and Bařina 7$n\pm$1 (proven) and 3$n$-1 and 5$n$+1 (refuted) in Theorem \ref{the:decidab}.

\subsection{Turing halting proof}
\label{sec:halting}

The 'halting proof' is a later term for the question posed by Alan Turing whether an (automated program executing an) iterated function, halts for all possible input numbers in a finite solution, as expected, or could be struck in a 'non-trivial' cyclic loop or a 'diverging trajectory' of calculations with ever higher numbers \cite{turing1936}. The halting problem for a Collatz-Conway a$n$+b function is whether numbers exist for which iterations do not halt in the cyclic trivial root number $c=a\cdot1+b$, but in a \textit{non-trivial cycle} or a \textit{diverging trajectory} with ever higher numbers \cite{RN3}. A non-trivial cycle or a diverging trajectory must have a lowest number $n_X\neq c$, possibly a number $n_X$ beyond computer power.

\begin{lemma}
\label{lem:halting}
The (supplementary) halting proof of a selected Collatz-Conway a\textit{n}+b conjecture that the trivial root is ‘always reached' is a proof ($T_{\ge0} = G_{\ge0}$) on contracting congruence classes --- or infinite arithmetic expressions --- generating repeatedely descending paths that these paths contract to the trivial root $\cabroot$ rather than to the lowest number  $n_x\neq c$ of a diverging trajectory or a non-trivial root.
\end{lemma}

\noindent
\textit{Proof}. Numbers in upward congruence classes never descend to lowest numbers such as $n_X$. Lightpurple congruence subclasses can give descending paths via yellow odd classes, presumably to $n_X$. 

For the Collatz $3n+1$ function (Fig.\ref{fig:cc3p1}) class $[52]_{54}$ is potentially a repeatedly descending congruence class allowing for paths to lower numbers, presumably to $n_X$. For the $7n\pm1$ function the depicted congruence class graph with periodicity $2^va=28$ \ref{fig:cc7pm1} fits on a single page, but the unabridged version with $n_b=7\cdot12=84$ with a LCM$p$-periodicity of lcm$_\textup{pq}=2^va^3=1372$ shows two descending, contracting, congruence classes $[456]_{1372}\subset[64]_{196}\subset[8]_{28}$ and $[916]_{1372}\subset[132]_{196}\subset[20]_{28}$. Successive breadth-first successor levels $k=1,2,3,\dots$ give for these contracting subclasses last nodes of contracting paths higher than first nodes at lower levels, disabling contracting paths to the lowest number $n_X$ of a non-trivial cycle.

\small{
\noindent 
\begin{alignat*}{7}
&\textrm{Contracting subclass} & & \textrm{contracting path}&&\quad k=0&&\;k=1&&\;k=2&&\;k=3&&\;k=4 \\
3n+1:\;&[52]_{54} & \quad\textrm{first node: } & (3^k-1)\cdot 54+52&&\quad\;52&&\quad160&&\quad\;484&&\;1\,456&&\;4\,372 \\
&& \quad\textrm{last node: } & (2^k-1)\cdot 54+52 &&\quad\;52&&\quad106&&\quad\;214&&\;430&&\;862 \notag \\
7n\pm1:\;&[456]_{1\,372} & \quad\textrm{first node: } &1/3 (7^k-1)\cdot 1372+456 &&\quad456&&\;3\,200&&\;22\,408&&\;156\,864&&\;1\,098\,056\notag \\
&&\quad\textrm{last node: } & 1/3 (4^k-1)\cdot 1372+456&&\quad456 &&\;1\,828&&\;7\,316&&\;29\,268 &&\;117\,076\notag \\
7n\pm1:\;&[916]_{1\,372} & \quad\textrm{first node: } &2/3 (7^k-1)\cdot 1372+916  &&\quad916&&\;6\,404&&\;44\,820&&\;313\,732&&\;2\,196\,116 \notag \\
&& \quad\textrm{last node: } & 2/3 (4^k-1)\cdot 1372+916   &&\quad916&&\;3\,660&&\;14\,636&&\;58\,540  &&\;234\,156 \notag\qed \\
\end{alignat*}
}
\normalsize{}

\vspace{-1.3cm}
\subsection{Completing the for `almost all' proof provided by the Syracuse tree}
\label{sec:syrac}
In his 2022 article “Almost all Collatz orbits attain almost bounded values” \cite{RN3051}, Fields medalist Terence Tao takes up the “ultimate challenge” posed by Jeffrey Lagarias, even though his admired graduate teacher John Conway \cite{taoconway} regarded a$n$+b conjectures as undecidable \cite{RN3209}. Tao concedes that “the full resolution of the conjecture remains well beyond the current methods". Alex Kontorovich, who compared with Lagarias the for “almost all” $3n+1$ function (Fig.~\ref{fig:tge0}) and the “not for all” $5n+1$ function (Fig.~\ref{fig:tree5p1}) \cite{RN19}, ventures that Tao's for "almost all” proof may have exceptions and may even be  \textit{undecidable} \cite{RN3154} (but see Lemmas \ref{lem:cantorno}-\ref{lem:wolframno}).

\begin{figure}[H]
\caption{Syracuse tree of odd numbers paired with even Collatz 3$n_{\textup{odd}}$+1 Cayley $T_{\ge0}$ numbers}
\vspace{-0.5cm}
\label{fig:syrtree}
\begin{center}
\includegraphics[width=0.85\textwidth]{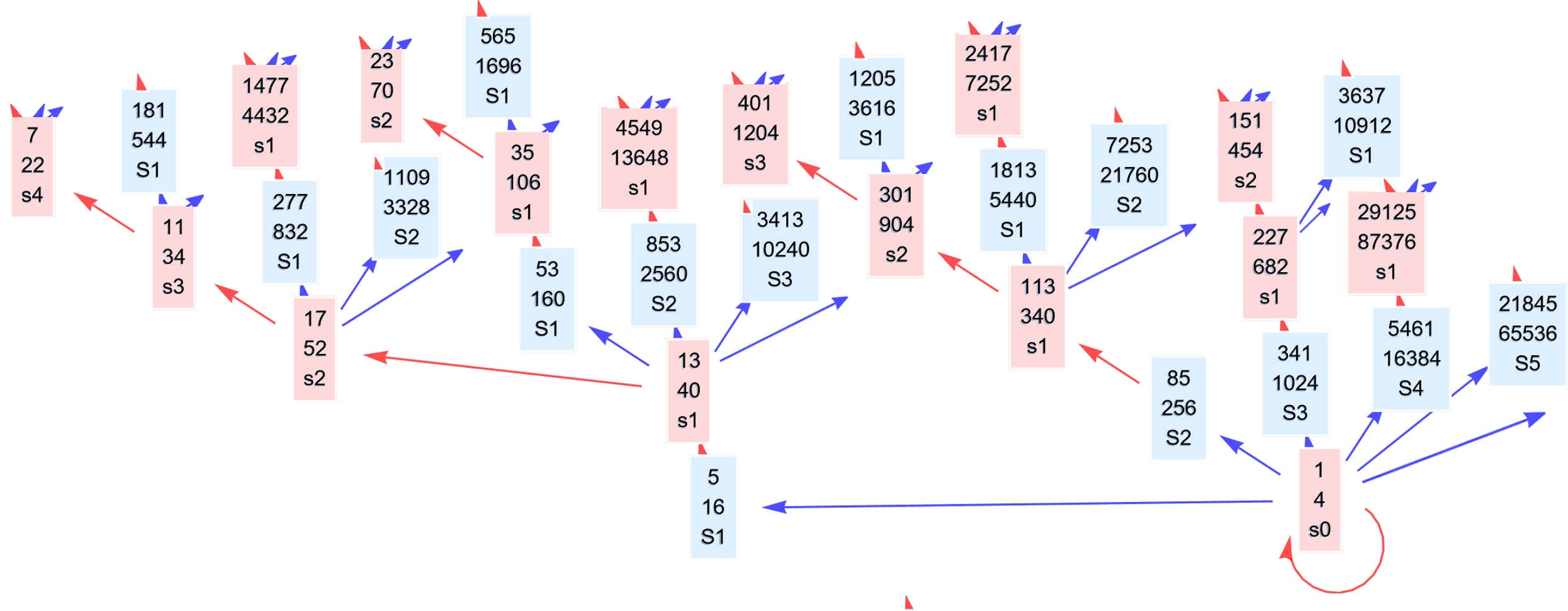}
\end{center}
\end{figure}

Consistent with the Cayley tree $T_{\ge0}$ (Eq.\ref{eq:bsclass}), Tao focuses on the yellow odd numbers (Figs.\ref{fig:gcol},\ref{fig:cc3p1})  from congruence classes $[1,5]_6$ indivisible by $a=3$ \citep[Prop.1.17,p.11]{RN3051}. Fig.\ref{fig:syrtree} of the $3n+1$ Syracuse tree depicts the $2^5=32$ odd numbers from the classes $[1,5]_6$ that are paired with red or blue even purple Collatz $n_{\textup{even}}=3\nmid3n_{\textup{odd}}+1$ classes $S_{\ge0}=[4,16]_{18}$ in the depicted Cayley 3$n$+1 tree (Fig.\ref{fig:tge0}a) , falling apart in $2^4=16$ red $L$eftward and $2^4=16$  blue $U$pward children (Eq.\ref{eq:Am}), paired with their red $s$- or blue $S$-generation. 

Lemma \ref{lem:syrproof} states that Tao's `for almost all' proof can be completed to a full proof, as it can be seen in Figure \ref{fig:syrtree} that only the cyclic root $\cabroot$ is retained when removing all infinite red paths of red number pairs and their red parental arrows including the blue arrows to their blue parental arrowhead nodes.      
\vspace{-1cm}
\begin{lemma}
\label{lem:syrproof}
Tao's `for almost all' proof of the Collatz 3$n$+1 conjecture on its Syracuse tree $T_S$ of odd numbers (Fig.\ref{lem:syrproof}) is completed to a 'for all', proof, by using instead of the Log-log-log-$\cdots$ function the integer-pair pruning function $\sm{\textcolor{red}{0}^{-i=1,2,3,\dots}_{s_i}1^{-1}\\L^{-i}U^{-1}(n)}:\;$ (Eq.\ref{eq:Idfsi}), which retains only the cyclic root $\cabroot$.
\end{lemma}
\noindent\textit{Proof.} Odd Syracuse tree $T_S$ numbers are paired one-by-one by $f:n_{\textup{odd}}\to3n+1$ to even Cayley $a\nmid n_{\textup{even}}=$3n$_{\textup{odd}}$+1 numbers, for which the 3$n$+1 conjecture is proven by the cotree density `Hamiltonian tour' 
proof and the `root-trunk' proof (Eqs.\ref{eq:cdp},\ref{eq:rtt}). Only serial red arrows to red $L$eftward arrowhead number pairs in $s$-paths starting from a blue $U$pward  from the Cayley tree $T_{\ge0}$, only $\tcb{\bullet}\tcr{\;\to}s_1\tcr{\to}s_2\tcr{\to}s_3\tcr{\to}\cdots$ paths, are retained in the Syracuse tree $T_S$. Tao approximates by $\textrm{log-log-log}\cdots$ integer binary root path numbers and integer a$n$+b numbers as non-integer real numbers, resulting in a non-deterministic for `almost all' proof on integer numbers. The for 'almost all' proof is completed by not using log-log-log-$\cdots$ functions, but instead the integer-pair pruning function $\sm{\textcolor{red}{0}^{-i=1,2,3,\dots}_{s_i}1^{-1}\\L^{-i}U^{-1}(n)}:\;$ (Eq.\ref{eq:Idfsi}) pruning/removing red serial paths together with their blue start number pair, thereby reducing Fig.\ref{fig:syrtree} to only the trivial red cyclic root $\cabroot$. \qed   

\section{Discussion}
\label{sec:disc}

\noindent The proof of the Collatz 3$n$+1 conjecture indeed connects number theory and graph theory observable through numbers \cite{RN9}, just as Pythagoras' theorem connects number theory and geometry observable through numbers. The infinite Turing-Leibniz-tape $T_L$ of $\eta=0,1,2,\dots$ proven a$n$+b conjectures (Fig.\ref{fig:TLEE}, Theorem\ref{the:TLEE}) exhibits such connections in its paired rows, expressible in terms of the cyclic roots $c=1a+1=2^{2^{\eta-1}+1}$ of their Cayley $T_{\ge0}$ trees of even numbers  $a\nmid$an$_{\textup{odd}}$+b. The $a=c-1$ parameters start with the short leg $a_{\eta1}=3$ of the Pythagorean triangle (3,4,5), with successive short legs of successor triangles being the sums of the legs in the preceding triangles (OEIS 365577), including Mersenne-prime short legs $a_{\eta}=2^{p-1}$. The $\psi=c-2$ counts of the $b$-parameters start with products of $2$ and successive Fermat primes $2^{2^p}+1$. The `cotree density' proof, expressible in terms of $c=1a+1$, is a proof of a graph-theoretic `Hamiltonian tour' of all Cayley $a\nmid$an$_{\textup{odd}}$+b numbers from the trivial root $\cabroot$ at the origin (0,0) of the Cartesian half graph $C_H$ and back. 

The Collatz 3$n$+1 first generation cotrees $t_1$ and $T_1$ comprise $153=\sum_{n=1}^{17}n=1^3+5^3+3^3$ disjoint congruence classes: $3^{2^2}=81$ red $L$efward $t_1$ classes plus $2^33^2=72$ blue $U$pward $T_1$ classes (Table \ref{tab:classes}).
The $153$ disjoint classes in cotrees $t_1$ and $T_1$ generate infinitely many disjoint biperiodic classes each with one number in an infinite number or congruence class periods $\ell=0,1,2,\dots$ in infinite numbers in disjoint red successor cotrees $t_{2,3,4,\dots}$ and blue successor cotrees $T_{2,3,4,\dots}$ of the first-generation red cotree $t_1$ and blue cotree $T_1$.
The 3$n$+1 numbers scattered across infinitely many classes with an infinite number of periods in infinitely many cotrees $t_{1,2,3,\dots}$ and $T_{1,2,3,\dots}$ are all uniquely identified by their binary root path numbers in the binary Cantor 1$n$+0 tree \cite[\S6]{cantor1883bin}, which is paired with the even Collatz 3$n$+1 numbers in Cayley tree $T_{\ge0}$ (Fig.\ref{fig:tge0}a). Number $153$ also connects the non-commutative, non-confluent, 3-regular Collatz-Cayley 3$n$+1 tree $T_{\ge0}$ with the both edge-transitive and vertex-transitive 3-regular Biggs-Smith graph with $153$ edges \cite{RN3229}. Remarkably, the association of $153$ with infinite numbers of followers is already told in 'the miraculous catch of fish' (John 21). Undecidability proofs by Cantor diagonals, Gödel numbers, Conway unpredictability and Wolfram irreducibility are disproven, as they did not reckon with biperiodic, 3-colour, a$n$+b Cayley trees with red $L$eftward $s$-successor generations and blue $U$pward $S$-successor generations (Lemmas \ref{lem:2comb}-\ref{lem:wolframno}. The Cantor 1$n$+0 tree of binary root path numbers is paired with the Conway 1$n$+1 tree (Fig.\ref{fig:treeP1}) and an infinite row of its Cayley a$n$+b successor trees, with Collatz 3$n$+1 (Fig.\ref{fig:tge0}) and Bařina 7$n\pm$1 (Fig.\ref{fig:tree7pm1}) first and second in row, for which the Collatz-Conway a$n$+b conjecture is proven by a Lean-certified cotree density Hamiltonian tour proof, and a root-successor proof  (Fig.\ref{fig:TLEE}, Theorem \ref{the:TLEE}). 

With their -(a$n$+b) counterparts, a$n$+b conjectures fill the full Cartesian integer coordinate system rather than just its upper half $C_H$ (Figs.\ref{fig:hamiltonCo},\ref{fig:hamiltonCo}). Natural language propositions about good and evil, cooperation and conflict, proponents and opponents, progression and regression, success and failure also map to the open interval $(-1,+1)$ approximated by the value range $-1..+1$ in a language logic graph, as in Heider's (1896-1988) triangular balance theory \cite{heider1946}, simplifying Freud's (1859-1939) defense mechanisms against inconsistencies \cite{freud39}. Frank Harary (1921-2005) generalized triangular balance theory  \cite{harary1956}, providing a context for interdisciplinary textual network content analysis, and later wrote `Graph Theory' \cite{RN33}, the first widely used textbook in computer science, network research and later AI Large Language Models that recognize and minimize inconsistency and ambiguity \cite{Cuil86graph,Cuil88AI,att2008semantic,att2021validity,attlaurer2024less}. The infinite variety of number-colon functions with terminal bi-periodic colour-colon digits, e.g. $1011\tcr{0}:$ resembling citations and attributions in politics, literature and scientific debates, suggests opportunities for (sub)disciplinary and interdisciplinary future network research.
\vspace{6.5cm}
\paragraph{Acknowledgments.}
\noindent The authors gratefully acknowledge the questions, suggestions and criticisms evoked by earlier drafts of David Bařina, Eldar Sultanow, the late Jan-Willem Klop, Jan van Mill, Jörg Endrullis, Wan Fokkink, Kas Kleinnijenhuis, Ronald Meester, Chris Verhoef, and Wouter van Atteveldt. The article draws on explorations using the Online Encyclopedia of Integer Sequences \cite{RN8}, sentence improvements suggested by ChatGPT, LaTeX-tikz plots enabled by $\textrm{https://q.uiver.app/}$ \cite{quiver}, and on the multitude of simulations and plots enabled by Mathematica \cite{RN32}. We thank Rudy H. Tanin for obtaining Lean-certificates of proofs ({\href{https://github.com/a4d442663580/collatz}{https://github.com/a4d442663580/collatz}) and Christian Koch for including the Cayley tree for Collatz $3n+1$ in his Github repository ({\href{https://github.com/c4ristian/collatz} {https://github.com/c4ristian/collatz}).

\cleardoublepage
\newgeometry{left=0.89in, right=1.5in}

\bibliographystyle{acm}

{\small
\bibliography{ms.bib} 
}

\end{document}